\documentclass[11pt,a4paper]{article}
\usepackage[margin=1.2in]{geometry}
\usepackage[utf8]{inputenc}
\usepackage[english]{babel}
\usepackage{amsmath}
\usepackage{graphicx}
\usepackage[dvipsnames]{xcolor}
\usepackage{mathtools}

\usepackage[linesnumbered,ruled,vlined]{algorithm2e}

\usepackage[charter, cal=cmcal]{mathdesign} 
\usepackage{natbib}

\usepackage{hyperref}
\usepackage{cleveref}

\usepackage{scalerel}
\usepackage{tikz}
\usetikzlibrary{svg.path}

\usepackage{booktabs}
\usepackage{pgfplotstable}
\usepackage{siunitx}
\usepackage{multicol}

\pgfplotsset{compat=1.18}
\usepackage{underscore}

\definecolor{orcidlogocol}{HTML}{A6CE39}
\tikzset{
	orcidlogo/.pic={
		\fill[orcidlogocol] svg{M256,128c0,70.7-57.3,128-128,128C57.3,256,0,198.7,0,128C0,57.3,57.3,0,128,0C198.7,0,256,57.3,256,128z};
		\fill[white] svg{M86.3,186.2H70.9V79.1h15.4v48.4V186.2z}
		svg{M108.9,79.1h41.6c39.6,0,57,28.3,57,53.6c0,27.5-21.5,53.6-56.8,53.6h-41.8V79.1z M124.3,172.4h24.5c34.9,0,42.9-26.5,42.9-39.7c0-21.5-13.7-39.7-43.7-39.7h-23.7V172.4z}
		svg{M88.7,56.8c0,5.5-4.5,10.1-10.1,10.1c-5.6,0-10.1-4.6-10.1-10.1c0-5.6,4.5-10.1,10.1-10.1C84.2,46.7,88.7,51.3,88.7,56.8z};
	}
}
\newcommand\orcidicon[1]{\href{https://orcid.org/#1}{\mbox{\scalerel*{
				\begin{tikzpicture}[yscale=-1,transform shape]
				\pic{orcidlogo};
				\end{tikzpicture}
			}{|}}}}

\DeclareMathOperator*{\argmin}{arg\,min}

\providecommand{\keywords}[1] {
  \small	
  \textbf{\textit{Keywords---}} #1
}

\title{The Augmented Mixing Method:\\ Computing High-Accuracy Primal-Dual Solutions\\ to Large-Scale SDPs via Column Updates}

\author{
Daniel Brosch\footnote{Institut f\"ur Mathematik,
  Alpen-Adria-Universit\"at Klagenfurt, Universit\"atstra{\ss}e 65--67,
  9020 Klagenfurt,
  \href{mailto:daniel.brosch@aau.at}{daniel.brosch@aau.at}, \,\href{mailto:jan.schwiddessen@aau.at}{jan.schwiddessen@aau.at}, \, \href{mailto:angelika.wiegele@aau.at}{angelika.wiegele@aau.at}} 
~\orcidicon{0000-0003-3988-4336}
\and Jan Schwiddessen$^{*}$\footnote{Corresponding Author: \href{mailto:jan.schwiddessen@aau.at}{jan.schwiddessen@aau.at}} \orcidicon{0000-0003-4380-2778}
\and Angelika Wiegele$^{*}$\orcidicon{0000-0003-1670-7951}
}

\begin{document}

\maketitle

\begin{abstract}
The Burer--Monteiro factorization has become a powerful tool for solving large-scale semidefinite programs (SDPs), enabling recently developed low-rank solvers to tackle problems previously beyond reach. However, existing methods are typically designed to prioritize scalability over solution accuracy.

We introduce the \emph{Augmented Mixing Method}, a new algorithm that combines the Burer--Monteiro factorization with an inexact augmented Lagrangian framework and a block coordinate descent scheme. Our method emphasizes solving low-dimensional subproblems efficiently and to high precision. Inequality constraints are handled directly, without explicitly maintaining slack variables in the algorithm. A novel dynamic update strategy for the penalty parameter ensures that primal and dual feasibility progress remain balanced. This approach enables our method to compute highly accurate primal-dual solutions, even for large-scale SDPs with over ten million inequality constraints.

Despite lacking theoretical convergence guarantees, the Augmented Mixing Method shows strong practical performance with default parameters across a wide range of SDP instances. It often produces more accurate primal-dual solutions than state-of-the-art interior-point methods and scales significantly better. Our open-source Julia implementation is memory-efficient, customizable, and supports arbitrary-precision arithmetic.
\end{abstract}



\section{Introduction}

\paragraph{Problem setting.}

In this paper, we consider linear semidefinite programs with arbitrary numbers of equality and inequality constraints, that is, optimization problems of the form
\begin{equation*}\label{eq:primal_sdp}
\begin{array}{rl}
\text{minimize} \quad & \langle C, X \rangle \\[1.2ex]
\text{subject to} \quad & \mathcal{A}(X) = a, \\
& \mathcal{B}(X) \geq b, \\
& X \in \mathcal{S}_+^n,
\end{array}
\tag{SDP}
\end{equation*}
where $C \in \mathcal{S}^n$, $a \in \mathbb{R}^{m_a}$, $b \in \mathbb{R}^{m_b}$, and $m_a, m_b \in \mathbb{N}_0$. We denote by $\langle A,B\rangle \coloneqq \operatorname{trace}(B^\top A)$ the trace inner product, and $\mathcal{S}_+^n$ is the cone of $n\times n$ symmetric positive semidefinite matrices. The linear operators $\mathcal{A} \colon \mathcal{S}^n \to \mathbb{R}^{m_a}$ and $\mathcal{B} \colon \mathcal{S}^n \to \mathbb{R}^{m_b}$ are defined componentwise by
\begin{equation*}
\mathcal{A}(X) =
\begin{pmatrix}
\langle A_1, X \rangle \\
\vdots \\
\langle A_{m_a}, X \rangle
\end{pmatrix},
\qquad
\mathcal{B}(X) =
\begin{pmatrix}
\langle B_1, X \rangle \\
\vdots \\
\langle B_{m_b}, X \rangle
\end{pmatrix},
\end{equation*}
for given matrices $A_1, \ldots, A_{m_a}$ and $B_1, \ldots, B_{m_b}$ in $\mathcal{S}^n$. 
 Problem~\eqref{eq:primal_sdp} can be extended to a semidefinite program involving several matrix variables. We will cover this generalization in Section~\ref{sec:multi-block}.

The dual problem corresponding to~\eqref{eq:primal_sdp} is
\begin{equation*}\label{eq:dual_sdp}
\begin{array}{rl}
\text{maximize} \quad & \langle a, y^a \rangle + \langle b, y^b \rangle \\[1.2ex]
\text{subject to} \quad & C - \mathcal{A}^\top(y^a) - \mathcal{B}^\top(y^b) = Z, \\
& y^a \in \mathbb{R}^{m_a}, ~ y^b \in \mathbb{R}_+^{m_b}, \\
& Z \in \mathcal{S}_+^n,
\end{array}
\tag{DSDP}
\end{equation*}
where the adjoint operators $\mathcal{A}^\top \colon \mathbb{R}^{m_a} \to \mathcal{S}^n$ and $\mathcal{B}^\top \colon \mathbb{R}^{m_b} \to \mathcal{S}^n$ are given by
\[
\mathcal{A}^\top(y^a) = \sum_{i=1}^{m_a} y^a_i A_i,
\quad
\mathcal{B}^\top(y^b) = \sum_{i=1}^{m_b} y^b_i B_i.
\]
Throughout the paper, we assume that the constraints in (SDP) are linearly independent and that both (SDP) and (DSDP) satisfy Slater's condition, ensuring strong duality and attainment of the optimal value for both problems.

\paragraph{Motivation.}

Semidefinite programming can be seen as a generalization of linear programming and is well known as a powerful modeling technique. Although semidefinite programs are convex optimization problems, they can be challenging to solve in practice. Under suitable regularity conditions, well-posed SDPs can be solved in polynomial time to arbitrary precision using interior-point methods~\cite{deklerk2016on}. However, despite their favorable theoretical complexity, primal-dual interior-point methods are known to scale poorly in practice: their computational cost typically grows at least cubically in both the order of the matrix variable and the number of constraints; see, e.g.,~\cite{borchers2007implementation}.

To address scalability, a variety of methods have been proposed in the literature, including eigenvalue optimization algorithms such as the spectral bundle method~\cite{helmberg2000spectral}, approaches based on the dual augmented Lagrangian~\cite{sun2020sdpnal,yang2015majorized}, related operator splitting methods like the alternating direction method of multipliers (ADMM)~\cite{povh2006boundary,wen2010alternating}, and pure penalty-based approaches for specifically structured SDPs~\cite{malick2013on}, among others. These methods typically offer lower per-iteration costs and better scalability to large-scale problems compared to interior-point methods. However, they often involve a trade-off between computational efficiency and solution accuracy. Some of these methods, in particular first-order approaches, are well-suited for obtaining medium-precision solutions efficiently, but achieving high-precision solutions remains challenging. Moreover, many existing approaches are tailored to specifically structured SDPs and are primarily used in combinatorial optimization, where the SDP serves as a convex relaxation of an \textnormal{NP}-hard problem. Extending these methods to broader classes of SDPs is nontrivial.

Another important class of methods seeks to exploit the low-rank structure that is often present in semidefinite programs. It has been shown in~\cite{barvinok1995problems,pataki1998rank,shapiro1982rank} that, for any SDP in standard form, there exists an optimal solution whose rank is bounded by the square root of twice the number of constraints. Motivated by this observation, some approaches avoid optimizing directly over the positive semidefinite matrix variable $X$ and instead apply a change of variables via the factorization $X = V^\top V$, where $V \in \mathbb{R}^{k \times n}$ with $k \ll n$. Optimization is then carried out with respect to $V$ rather than $X$. This idea, commonly referred to as the Burer--Monteiro factorization approach~(BM)~\cite{burer2003nonlinear,burer2005local}, significantly reduces the number of variables and enables more scalable algorithms. However, both the cost function and the affine constraints become quadratic in $V$, introducing nonconvexity into the problem. 

For generic cost matrices, despite the nonconvexity of the formulation, it was shown in~\cite{boumal2020deterministic} that the BM approach successfully solves semidefinite programs with equality constraints and a chosen rank $k$ of order $\sqrt{2m}$, where $m$ is the number of constraints. This result was extended to arbitrary semidefinite programs with inequality constraints and multiple positive semidefinite matrix variables in~\cite{cifuentes2021burer}. However, the method can fail in nongeneric settings~\cite{ocarroll2022theburermonteiro,waldspurger2020rank}.

It was shown in~\cite{cifuentes2022polynomial} that the BM method admits polynomial-time guarantees after slight perturbation, and can solve SDPs to any desired accuracy.

\paragraph{Our contribution.}
In this work, we propose a new method built on the Burer--Monteiro factorization, targeting SDPs in which the number of constraints far exceeds the matrix size. The method further decomposes the reformulated problem into smaller, lower-dimensional subproblems by applying a block coordinate descent scheme. As a result, the subproblems are very small and can be solved efficiently to high accuracy using a warm-started gradient-based method from previous iterates. Our main contributions are as follows:
\begin{itemize}
    \item We introduce a new algorithm, called the \emph{Augmented Mixing Method}, for solving SDPs with arbitrary numbers of equality and inequality constraints. The algorithm is particularly suited for problems with a large number of inequality constraints, as these are handled directly without reformulation or maintaining slack variables in the algorithm.
    \item We develop a novel dynamic update strategy for the penalty parameter that reliably balances primal and dual progress. This enables the computation of highly accurate solutions on many SDP instances, often reaching primal-dual errors close to the limits of the working precision.
    \item The algorithm has low per-iteration cost, relying primarily on (sparse) matrix-vector multiplications. It avoids eigenvalue computations, projections onto the cone of positive semidefinite matrices, and matrix factorizations of the SDP data. Moreover, in our experiments, stagnation in the primal-dual errors sets in later than for interior-point solvers.
    \item We present extensive numerical experiments across a variety of SDP types using the solver’s default parameters. These experiments demonstrate the algorithm’s effectiveness and scalability, especially on problems with a large number of constraints, and show competitive performance against state-of-the-art interior-point solvers in both runtime and accuracy.
    \item We provide an open-source Julia implementation. The solver is fully type-generic, supports arbitrary-precision arithmetic, and can be warm-started from any initial point. Numerical experiments confirm that it remains efficient and accurate when using extended floating-point precision.
\end{itemize}
To the best of our knowledge, this is the first solver based on the Burer--Monteiro factorization to combine these features in a unified and practical framework.

\paragraph{Outline of the paper.}
The remainder of this paper is organized as follows. In Section~\ref{sec:relation}, we introduce the Burer--Monteiro factorization and discuss how the Augmented Mixing Method relates to other approaches in the literature. Section~\ref{sec:augmented_mixing} presents our new algorithm, including the procedure for performing column updates, our novel strategy for dynamically updating the penalty parameter, and a high-level summary of the overall method. In Section~\ref{sec:further_ingredients}, we describe implementation details such as data scaling, stopping criteria, and user-configurable parameters. We also explain how to use our Julia implementation and describe the workflow for computations with extended floating-point precision. Section~\ref{sec:computational_experiments} presents extensive numerical experiments on various classes of SDPs, comparing the Augmented Mixing Method to other available SDP solvers. Finally, Section~\ref{sec:conclusion} summarizes our contributions and outlines possible directions for future work.

\paragraph{Notation.} 
We denote the Frobenius norm of a matrix $A$ by $\lVert A \rVert_F$. For vectors, we use $\lVert \cdot \rVert_2$ for the Euclidean norm and $\lVert \cdot \rVert_\infty$ for the maximum norm. Given a matrix $X \in \mathcal{S}^n$, we write $[X]_+$ for its metric projection onto $\mathcal{S}_+^n$, that is, the matrix obtained by zeroing out all negative eigenvalues in its spectral decomposition. For a vector $x \in \mathbb{R}^n$, the notation $[x]_+$ denotes the vector in $\mathbb{R}^n$ whose $i$th component is given by $\max\{x_i, 0\}$. The operator $\operatorname{diag} \colon \mathbb{R}^{n \times n} \rightarrow \mathbb{R}^n$ maps a square matrix to a vector consisting of its diagonal elements. The adjoint operator of $\operatorname{diag}$ is denoted by $\operatorname{Diag} \colon \mathbb{R}^n \rightarrow \mathcal{S}^{n}$. We denote the all-ones matrix of order $n$ by $J_n$ and the all-ones vector of dimension $n$ by $e_n$. For $n \in \mathbb{N}$, we write $[n]$ for the set $\{1,\ldots,n\}$. For a matrix $A \in \mathcal{S}^n$ and $i \in [n]$, we use the notation $A_{(i)}$ to refer to the $i$th column of $A$ and $A_{(ii)}$ to refer to the $i$th diagonal entry.

\section{Relation to Other Methods} \label{sec:relation}

The \emph{Augmented Mixing Method} combines and extends ideas from two prominent approaches for solving semidefinite programs: the Burer--Monteiro factorization framework~\cite{burer2003nonlinear,burer2005local} and its implementation \texttt{SDPLR}, as well as the \emph{mixing method} introduced in~\cite{wang2017mixing}. While the mixing method is highly specialized and limited to diagonally constrained SDPs, \texttt{SDPLR} is designed to solve equality-constrained semidefinite programs in standard form. Both approaches share the idea of employing a low-rank factorization. In this section, we briefly recall the core ideas behind \texttt{SDPLR} and the mixing method. 

\subsection{The Burer--Monteiro approach and \texttt{SDPLR}}

Burer and Monteiro~\cite{burer2003nonlinear,burer2005local} proposed solving SDPs in standard form
\begin{equation}\label{eq:bmsdp}
\begin{array}{rl}
\text{minimize} \quad & \langle C, X \rangle \\[1.2ex]
\text{subject to} \quad & \mathcal{A}(X) = a, \\
& X \in \mathcal{S}_+^n,
\end{array}
\end{equation}
by introducing a low-rank factorization $X = V^\top V$ using a rectangular matrix $V \in \mathbb{R}^{k \times n}$. This leads to the reformulated problem
\begin{equation}\label{eq:bmvec}
\begin{array}{rl}
\text{minimize} \quad & \langle C, V^\top V \rangle \\[1.2ex]
\text{subject to} \quad & \mathcal{A}(V^\top V) = a, \\
& V \in \mathbb{R}^{k \times n}.
\end{array}
\end{equation}
It is well known from~\cite{barvinok1995problems,pataki1998rank,shapiro1982rank} that~\eqref{eq:bmsdp} and~\eqref{eq:bmvec} are equivalent if $k \geq \lceil \sqrt{2m} \rceil$, where $m$ is the number of constraints.

The approach of~\cite{burer2003nonlinear,burer2005local} combines two ideas: first, they apply a standard augmented Lagrangian method to~\eqref{eq:bmvec}, solving the subproblems using a limited-memory BFGS method. Second, they use a dynamic strategy for updating the rank $k$, starting with a small value and increasing it at runtime only if necessary. Together, these ideas form the basis of the solver \texttt{SDPLR}. The solver proved successful on some large-scale instances of~\eqref{eq:bmsdp} that were difficult for other solvers at the time.

However, a known limitation of \texttt{SDPLR} is that it typically yields only moderately accurate primal feasible solutions, and the associated dual solutions are sometimes far from optimal. This shortcoming arises because \texttt{SDPLR} uses a standard augmented Lagrangian approach that does not explicitly target finding an accurate primal-dual pair. In contrast, the Augmented Mixing Method addresses this issue by adopting and enhancing the augmented Lagrangian scheme. In this scheme, subproblems are solved efficiently to any desired accuracy using a block coordinate descent strategy. Paired with a novel dynamic update scheme for the penalty parameter, the Augmented Mixing Method is able to compute high-accuracy primal-dual solutions. Unlike \texttt{SDPLR}, we deliberately keep the rank~$k$ fixed during runtime. Dynamic rank updates are not essential in our setting, particularly for SDPs with many constraints. The performance of our method does not depend strongly on~$k$, so little is gained by adjusting the rank adaptively. Moreover, the method behaves more stably for larger values of~$k$, and we therefore fix~$k$ from the beginning to a value that yields robust behavior.

Another solver that uses a low-rank factorization together with an (inexact) augmented Lagrangian approach is \texttt{HALLaR}~\cite{monteiro2024low}. It is a hybrid method based on an adaptive inexact proximal-point framework with inner acceleration and Frank--Wolfe steps to escape spurious local stationary points. The authors report numerical results showing that \texttt{HALLaR} achieves state-of-the-art performance on very large-scale SDP instances, with relative accuracy around $10^{-5}$. A GPU implementation, \texttt{cuHALLaR}, is described in~\cite{aguirre2025cuhallar}. Related approaches, such as those in~\cite{wang2025solving,wang2023decomposition}, combine the augmented Lagrangian framework with manifold optimization.

A further low-rank solver, called \texttt{LoRADS}, is proposed in~\cite{han2025low}. It uses a bilinear factorization $X = U^\top V$ with the additional constraint $U = V$. This constraint is incorporated into the objective function via a penalty term, and the resulting problem is solved using an ADMM framework. A related idea was proposed in~\cite{chen2023burermonteiro} for diagonally constrained SDPs; see Section~\ref{sec:mixing_method}. Like \texttt{SDPLR}, \texttt{LoRADS} employs a dynamic rank update strategy. The authors of~\cite{han2025low} report promising results on a variety of SDP instances, including some very large-scale problems. The numerical experiments indicate that \texttt{LoRADS} places more emphasis on achieving good primal feasibility, as the reported dual errors are typically larger than the corresponding primal errors. A GPU implementation of \texttt{LoRADS} is described in~\cite{han2024accelerating}.

\subsection{The Mixing Method}\label{sec:mixing_method}

The \emph{mixing method} was proposed in~\cite{wang2017mixing} to solve diagonally constrained SDPs of the form
\begin{equation}\label{eq:max_cut_sdp}
\begin{array}{rl}
\text{maximize} \quad & \langle C, X \rangle \\[1.2ex]
\text{subject to} \quad & X_{ii} = 1, \quad i = 1,\dots,n, \\
& X \in \mathcal{S}_+^n,
\end{array}
\end{equation}
which includes the well-known Goemans--Williamson SDP relaxation for the Max-Cut problem. Applying the Burer--Monteiro factorization $X = V^\top V$ with a column-wise stored matrix $V = (v_1 \mid \cdots \mid v_n) \in \mathbb{R}^{k \times n}$ leads to the reformulated problem
\begin{equation}\label{eq:mixingnonconvex}
\begin{array}{rl}
\text{maximize} \quad & \sum_{i,j=1}^n C_{ij} v_i^\top v_j \\[1.2ex]
\text{subject to} \quad & \lVert v_i \rVert_2 = 1, \quad i = 1,\dots,n, \\
& v_i \in \mathbb{R}^k, \quad i = 1,\dots, n.
\end{array}
\end{equation}
The mixing method then applies a block coordinate ascent approach with respect to the columns $v_i$, $i = 1,\ldots,n$. By selecting a fixed $i \in [n]$ and fixing all columns but $v_i$, Problem~\eqref{eq:mixingnonconvex} reduces to
\begin{equation} \label{eq:mixing_reduced}
\begin{array}{rl}
\text{maximize} \quad & g^\top v_i \\[1.2ex]
\text{subject to} \quad & \lVert v_i \rVert_2 = 1, \\
& v_i \in \mathbb{R}^k,
\end{array}
\end{equation}
where $g = \sum_{j \neq i} C_{ij} v_j$. Although this subproblem remains nonconvex, it has a closed-form solution for~$v_i$: specifically, $v_i^\ast = g / \lVert g \rVert_2$ when $g \neq 0$.

Under a mild nondegeneracy assumption and for $k > \sqrt{2n}$, Wang, Chang, and Kolter~\cite{wang2017mixing} establish local linear convergence of a block coordinate ascent scheme using this closed-form update in a cyclic fashion and starting from a random feasible point. A key feature of the mixing method is that it always maintains primal feasibility and guarantees strict monotonic improvement of the objective value. In addition, approximate dual multipliers can be constructed from any primal iterate, and their sequence also converges to the optimum~\cite{wang2017mixing}. Extensions and refinements of the mixing method for very large-scale problems have been proposed in~\cite{chen2023burermonteiro,erdogdu2022convergence,kim2021momentuminspired}.

As in the mixing method, our proposed algorithm updates only a single column of the BM factorization matrix $V$ at a time. However, since no closed-form solution is available in our setting, we use gradient-based methods to perform these column updates. Nevertheless, the computational effort required for each column update remains low, allowing us to solve the low-dimensional subproblems numerically to any desired precision.

\section{The Augmented Mixing Method} \label{sec:augmented_mixing}

In this section, we introduce our new algorithm, the \emph{Augmented Mixing Method}.
The Augmented Mixing Method shares with the mixing method the use of a block coordinate descent scheme, resulting in small subproblems that can be efficiently solved to high accuracy. At the same time, it shares with \texttt{SDPLR} the use of an augmented Lagrangian framework, which provides direct access to the dual variables. However, the Augmented Mixing Method is specifically designed to compute high-accuracy primal-dual solutions and extends the applicability of low-rank approaches to semidefinite programs with general affine constraints, including a large number of inequality constraints. 

\subsection{Burer--Monteiro factorization}

We reformulate~\eqref{eq:primal_sdp} via the Burer--Monteiro factorization, which allows us to eliminate the positive semidefiniteness constraint on the matrix variable $X$. Specifically, we introduce a rectangular matrix $V \in \mathbb{R}^{k \times n}$ for some $k \in \mathbb{N}$ and set $X = V^\top V$. For sufficiently large $k$,~\eqref{eq:primal_sdp} is then equivalent to
\begin{equation*}\label{eq:reformulated_sdp}
\begin{array}{rl}
\text{minimize} \quad & \langle C, V^\top V \rangle \\[1.2ex]
\text{subject to} \quad & \mathcal{A}(V^\top V) = a, \\
& \mathcal{B}(V^\top V) \geq b, \\
& V \in \mathbb{R}^{k \times n}.
\end{array}
\tag{BM-SDP}
\end{equation*}
While~\eqref{eq:primal_sdp} is a convex optimization problem, the reformulated problem~\eqref{eq:reformulated_sdp} is generally nonconvex due to the quadratic nature of both the objective and the constraints. In accordance with the \emph{Barvinok--Pataki bound}~\cite{barvinok1995problems,pataki1998rank}, we choose
\begin{equation}\label{eq:our_k}
k = \min\left(n, \left\lceil \sqrt{2(m_a + m_b)} \right\rceil\right)
\end{equation}
throughout the Augmented Mixing Method. This choice proved robust and efficient across all our numerical experiments. Note that for SDP instances with many constraints, it can yield $k=n$.

\subsection{Augmented Lagrangian and its gradient}

We adopt the framework of the augmented Lagrangian method, originally developed for equality-constrained problems by Hestenes~\cite{hestenes1969multiplier} and Powell~\cite{powell1969method}. To handle the inequality constraints in~\eqref{eq:reformulated_sdp}, we use Rockafellar's extension~\cite{rockafellar1973multiplier}. The first step in applying this extension to~\eqref{eq:reformulated_sdp} is to transform each inequality constraint $\langle B_i, V^\top V \rangle \geq b_i$ into an equality constraint of the form
\begin{equation*}
    b_i - \langle B_i, V^\top V \rangle + s_i^2 = 0
\end{equation*}
by introducing slack variables $s_i \in \mathbb{R}$ for $i = 1,\ldots,m_b$. With $s=(s_i)_{i=1}^{m_b}$, forming the standard augmented Lagrangian for the resulting equality-constrained problem yields
\begin{equation*}
\begin{aligned}
\bar{\mathcal{L}}(V, s, y^a, y^b; \mu) \coloneqq {} & \langle C, V^\top V \rangle 
+ \langle y^a, a - \mathcal{A}(V^\top V) \rangle
+ \frac{\mu}{2} \| a - \mathcal{A}(V^\top V) \|_2^2 \\
& + \sum_{i=1}^{m_b} \left[
y_i^b \left( b_i - \langle B_i, V^\top V \rangle + s_i^2 \right)
+ \frac{\mu}{2} \left( b_i - \langle B_i, V^\top V \rangle + s_i^2 \right)^2
\right].
\end{aligned}
\end{equation*}
Finally, the slack variables can be eliminated by minimizing $\bar{\mathcal{L}}$ with respect to them while keeping $V$ fixed; see, e.g., Section~5.4.2 in~\cite{geiger2002theorie}. This gives
\begin{equation*}
\begin{aligned}
\mathcal{L}(V, y^a, y^b; \mu) \coloneqq {}
& \min_s \bar{\mathcal{L}}(V, s, y^a, y^b; \mu) \\
={}& \langle C, V^\top V \rangle 
+ \langle y^a, a - \mathcal{A}(V^\top V) \rangle
+ \frac{\mu}{2} \| a - \mathcal{A}(V^\top V) \|_2^2 \\
& + \frac{1}{2\mu}
\left\lVert
\left[
y^b + \mu \left( b - \mathcal{B}(V^\top V) \right)
\right]_+
\right\rVert^2_2
- \frac{1}{2\mu} \left\lVert y^b \right\rVert_2^2,
\end{aligned}
\end{equation*}
which is equivalent to the standard Rockafellar augmented Lagrangian for inequality-constrained optimization problems. Here, $y^a \in \mathbb{R}^{m_a}$ and $y^b \in \mathbb{R}_+^{m_b}$ are the dual variables associated with the equality and inequality constraints, respectively. The scalar $\mu > 0$ is the penalty parameter.

For later use, note that
\begin{equation*}
\begin{split}
&\frac{1}{2\mu} \left\lVert \left[ y^b + \mu \left( b - \mathcal{B}(V^\top V) \right) \right]_+ \right\rVert^2_2 - \frac{1}{2\mu} \left\lVert y^b \right\rVert^2_2 \\
&= \frac{\mu}{2} \sum_{j \in I} (b_j - \langle B_j, V^\top V \rangle)^2 + \sum_{j \in I} y_j^b (b_j - \langle B_j, V^\top V \rangle) - \frac{1}{2\mu} \sum_{j \in I^c} (y_j^b)^2,
\end{split}
\end{equation*}
where the index sets $I$ and $I^c$ are defined as
\begin{equation*}
I \coloneqq \left\{ j \in [m_b] \colon y^b_j + \mu \left( b_j - \langle B_j, V^\top V \rangle \right) > 0 \right\}, \qquad
I^c \coloneqq [m_b] \setminus I.
\end{equation*}
Putting everything together, we obtain the expanded form
\begin{equation*}
\begin{aligned}
\mathcal{L}(V, y^a, y^b; \mu) = {} & \langle C, V^\top V \rangle 
+ \frac{\mu}{2} \sum_{j \in [m_a]} (a_j - \langle A_j, V^\top V \rangle)^2 + \sum_{j \in [m_a]} y_j^a (a_j - \langle A_j, V^\top V \rangle) \\
& + \frac{\mu}{2} \sum_{j \in I} (b_j - \langle B_j, V^\top V \rangle)^2 + \sum_{j \in I} y_j^b (b_j - \langle B_j, V^\top V \rangle) - \frac{1}{2\mu} \sum_{j \in I^c} (y_j^b)^2.
\end{aligned}
\end{equation*}
The augmented Lagrangian is continuously differentiable, and its gradient with respect to $V$ is given by
\begin{equation}\label{eq:gradientV}
\nabla_V \mathcal{L}(V, y^a, y^b; \mu) = 2 V \left( C - \sum_{j \in [m_a]} \left[ y^a_j + \mu \left( a_j - \langle A_j, V^\top V \rangle \right) \right] A_j - \sum_{j \in I} \left[ y^b_j + \mu \left( b_j - \langle B_j, V^\top V \rangle \right) \right] B_j \right).
\end{equation}
Note that not all inequality constraints are relevant for computing the gradient~\eqref{eq:gradientV}. Only those in the index set $I$ contribute and hence need to be considered.

Following the standard augmented Lagrangian approach, an approximate minimizer of~\eqref{eq:reformulated_sdp} can be computed by initializing $V$, $y^a$, $y^b$, and $\mu$, then minimizing the augmented Lagrangian for fixed $y^a$, $y^b$, and $\mu$ using the gradient~\eqref{eq:gradientV}, and subsequently updating $y^a$, $y^b$, and $\mu$ iteratively until a suitable stopping criterion is satisfied. However, using this approach would require solving a sequence of subproblems in $kn$ variables, for which it is challenging to obtain high-accuracy solutions. This limitation is essentially the bottleneck in \texttt{SDPLR} and, due to the resulting numerical challenges, it can typically attain only moderate accuracy.

\subsection{Column updates} \label{sec:column_updates}

To overcome this limitation, we employ a \emph{block coordinate descent} strategy based on column updates. This yields low-dimensional subproblems that can be solved to high precision, which is crucial for attaining high-accuracy primal-dual solutions.

At each step of the Augmented Mixing Method, all columns except one are fixed, and the augmented Lagrangian is minimized only with respect to the selected column, which is then updated accordingly. In this way, each subproblem involves only $k$ variables.

To formalize this, we write $\overline{V}(v_i)$ for the $k \times n$ matrix obtained by keeping all columns of $V$ fixed except for the $i$th column, which is treated as a variable. That is, $\overline{V}(v_i) = (v_1 \mid \cdots \mid v_i \mid \cdots \mid v_n)$. When updating column $v_i$ of $V$, its new value is obtained via
\begin{equation}\label{eq:subproblem}
v_i \leftarrow \argmin_{v_i \in \mathbb{R}^k} \mathcal{L}(\overline{V}(v_i), y^a, y^b; \mu),
\end{equation}
where the dual variables $y^a$, $y^b$, the penalty parameter $\mu$, and all other columns of $V$ are held fixed.

Subproblem~\eqref{eq:subproblem} is a nonlinear optimization problem in $k$ variables with a continuously differentiable objective function. According to~\eqref{eq:gradientV}, the partial derivative of $\mathcal{L}$ with respect to $v_i$ is given by
\begin{equation}\label{eq:gradient_vi}
\begin{split}
\nabla_{v_i} \mathcal{L}(\overline{V}(v_i), y^a, y^b; \mu) = &
2 \overline{V}(v_i) \Bigg(
    C_{(i)}
    - \sum_{j \in [m_a]} \left[ y^a_j + \mu \left( a_j - \langle A_j, \overline{V}(v_i)^\top \overline{V}(v_i) \rangle \right) \right] (A_j)_{(i)} \\
& \quad - \sum_{j \in I} \left[ y^b_j + \mu \left( b_j - \langle B_j, \overline{V}(v_i)^\top \overline{V}(v_i) \rangle \right) \right] (B_j)_{(i)}
\Bigg).
\end{split}
\end{equation}

To approximately solve subproblem~\eqref{eq:subproblem}, we employ a limited-memory BFGS method, warm-started from the previous iterate $v_i^{\textrm{start}}$. Throughout the algorithm, we use an inexact augmented Lagrangian approach in which the subproblems are solved only to a stopping tolerance that becomes more stringent as the iterates progress. Specifically, we terminate the subproblem solver and accept the current $v_i$ as an approximate solution to~\eqref{eq:subproblem} as soon as the following stopping criterion is satisfied:
\begin{equation}\label{eq:subproblemstop}
    \lVert \nabla_{v_i} \mathcal{L}(\overline{V}(v_i), y^a, y^b; \mu) \rVert_\infty < \min \left\{ \varepsilon,\, \delta \lVert \nabla_{v_i} \mathcal{L}(\overline{V}(v_i^{\textrm{start}}), y^a, y^b; \mu) \rVert_\infty \right\},
\end{equation}
where $\varepsilon,\delta > 0$ are predetermined parameters. We use the default values $\varepsilon = \delta = 0.01$.

\subsection{Dual and penalty parameter update} \label{sec:update}

After each full pass of column updates, the dual variables $y^a$ and $y^b$ are updated using the standard augmented Lagrangian rules
\begin{equation}\label{eq:dual_update}
y^a \leftarrow y^a + p \mu (a - \mathcal{A}(V^\top V)), \quad
y^b \leftarrow \left[ y^b + p \mu (b - \mathcal{B}(V^\top V)) \right]_+,
\end{equation}
where we use a default value of $p = 1$ for the dual step size. Although this choice is used in all experiments in Section~\ref{sec:computational_experiments}, our implementation allows users to specify a different value of $p > 0$ if desired (see Section~\ref{sec:parameters}). 

While the update of the dual variables is straightforward, updating the penalty parameter $\mu$ is crucial and has a strong impact on the performance of augmented Lagrangian methods in general. In many implementations of augmented Lagrangian methods, $\mu$ is either kept constant or increased monotonically, but rarely decreased. In contrast, we found that adjusting $\mu$ in both directions, depending on the progress in primal feasibility, is essential to reliably obtain high-accuracy primal-dual solutions. A different penalty parameter update rule that allows for reductions is used in~\cite{wang2025solving}.

\newcommand{\IActive}{\mathcal{I}}
\newcommand{\projActive}{P_\IActive}

After each full sweep over all columns of $V$, we define the active index set of the inequality constraints as
\begin{equation*}
\IActive \coloneqq \left\{ j \in [m_b] \;\colon\; b_j - \langle B_j, V_{\mathrm{new}}^\top V_{\mathrm{new}} \rangle \geq 0 \text{ or } y_j^b > 0 \right\},
\end{equation*}
and introduce the projection operator \( \projActive \colon \mathbb{R}^{m_b} \to \mathbb{R}^{\IActive} \), \( \projActive(u) \coloneqq (u_j)_{j\in\IActive}\), which removes all inactive components. We define the ratio
\begin{equation}\label{eq:penalty_ratio}
\texttt{ratio} \coloneqq
\frac{
    \left\lVert 
    \begin{pmatrix}
        a - \mathcal{A}(X_{\text{new}}) \\
        \projActive\left(b - \mathcal{B}(X_{\text{new}})\right)
    \end{pmatrix}
    \right\rVert_2
}{
    \mu \left\lVert 
    \begin{pmatrix}
        \mathcal{A}(X_{\text{new}} - X_{\text{old}}) \\
        \projActive\left( \mathcal{B}(X_{\text{new}} - X_{\text{old}}) \right)
    \end{pmatrix}
    \right\rVert_2
},
\end{equation}
where $X_{\text{old}} = V_{\text{old}}^\top V_{\text{old}}$ and $X_{\text{new}}=V_{\text{new}}^\top V_{\text{new}}$ denote the primal iterate before and after the most recent primal update, respectively. We aim to keep \texttt{ratio} close to $1$ throughout our algorithm.

Assuming that the iterates move directly towards feasibility, meaning that $\mathcal{A}(X_{\text{new}} - X_{\text{old}})$ is a scalar multiple of $\mathcal{A}(X_{\text{new}}) - a$ (and analogously for inequalities), the primal infeasibility decreases geometrically when \texttt{ratio} is close to $1$. In such cases, each iteration reduces the primal residual by a factor of $1 - (\mu+1)^{-1}$. This ensures we progress toward primal feasibility quickly enough to keep the number of iterations reasonable, while avoiding convergence to primal feasibility too early. This is important because once the constraints $\mathcal{A}(X) = a$ and $\mathcal{B}(X) \geq b$ are nearly satisfied, the dual variables are barely updated, limiting further improvement in dual feasibility.

If \texttt{ratio} exceeds the prescribed upper threshold \texttt{rat\_max}, we increase $\mu$ by the factor $\tau$; if it falls below the lower threshold \texttt{rat\_min}, we decrease $\mu$ by the factor $1/\tau$. Otherwise, $\mu$ remains unchanged. Our update rule is:
\begin{equation}\label{eq:penalty_update}
\mu \leftarrow 
\begin{cases}
\tau \mu, & \text{if } \texttt{ratio} > \texttt{rat\_max}, \\
\mu / \tau, & \text{if } \texttt{ratio} < \texttt{rat\_min}, \\
\mu, & \text{otherwise}.
\end{cases}
\end{equation}
We use the default values $\texttt{rat\_min} = 0.8$, $\texttt{rat\_max} = 1.2$, and $\tau = 1.03$.

\subsection{Algorithm and algorithmic complexity}

Algorithm~\ref{alg:augmented_mixing} summarizes the core steps of the Augmented Mixing Method, which alternates between column-wise updates of the matrix variable and updates of the dual variables together with the penalty parameter. Unless warm-started (see Section~\ref{sec:warmstart}), we initialize each column $v_i$ of $V$ uniformly at random on the unit sphere, set all dual variables to zero, and choose an initial penalty parameter of $\mu = \sqrt{n}$. For implementation details, including scaling strategies, stopping criteria, tolerances, and various extensions, we refer the reader to Section~\ref{sec:further_ingredients}.

Most of the computational effort in the Augmented Mixing Method is spent evaluating the restricted augmented Lagrangian in~\eqref{eq:subproblem} and its gradient~\eqref{eq:gradient_vi}. A particularly expensive step is evaluating the linear operators \( \mathcal{A}(\overline{V}(v_i)^\top \overline{V}(v_i)) \) and \( \mathcal{B}(\overline{V}(v_i)^\top \overline{V}(v_i)) \) at the current trial point \( v_i \). In the following, we describe how \( \mathcal{A}(\overline{V}(v_i)^\top \overline{V}(v_i)) \) can be computed efficiently using incremental updates that reuse already computed information. The update of \( \mathcal{B}(\cdot) \) works analogously.

Due to warm-starting from the point \( v_i^{\text{start}} \), we can assume the value
\begin{equation*}
\mathcal{A}^{\text{start}} \coloneqq \mathcal{A}(\overline{V}(v_i^{\text{start}})^\top \overline{V}(v_i^{\text{start}}))
\end{equation*}
is already available. Then, for any trial point $v_i$, the updated value can be computed as
\begin{equation}\label{eq:AX_incr}
\mathcal{A}(\overline{V}(v_i)^\top \overline{V}(v_i)) = \mathcal{A}^{\text{start}} + (\|v_i\|_2^2 - \|v_i^{\text{start}}\|_2^2)\, ((A_j)_{(ii)})_{j=1}^{m_a} + 2 \, \hat{A}_i^\top \overline{V}(0)^\top (v_i - v_i^{\text{start}}),
\end{equation}
where $\hat{A}_i \coloneqq ( (A_1)_{(i)} \mid \cdots \mid (A_{m_a})_{(i)} ) \in \mathbb{R}^{n \times m_a}$. Note that typically not all $m_a$ entries of $\mathcal{A}(\overline{V}(v_i)^\top \overline{V}(v_i))$ need to be computed due to sparsity. We evaluate~\eqref{eq:AX_incr} by first computing the dense matrix-vector product $\overline{V}(0)^\top (v_i - v_i^{\text{start}})$, which can be carried out in $\mathcal{O}(kn)$. After this, the most expensive step is to multiply the matrix $\hat{A}_i^\top$ from the left. For dense constraint matrices, this can be done in $\mathcal{O}(m_a n)$. However, if $\hat{A}_i^\top$ is a sparse matrix, the computational effort is proportional to $\mathrm{nnz}(\hat{A}_i)$, where $\mathrm{nnz}(\hat{A}_i)$ denotes the number of nonzeros in $\hat{A}_i$. Moreover, when evaluating the gradient~\eqref{eq:gradient_vi}, we have to compute
\begin{equation*}
    - \sum_{j \in [m_a]} \left[ y^a_j + \mu \left( a_j - \langle A_j, \overline{V}(v_i)^\top \overline{V}(v_i) \rangle \right) \right] (A_j)_{(i)} = - \hat{A}_i \left(y^a + \mu \left(a - \mathcal{A}(\overline{V}(v_i)^\top \overline{V}(v_i))\right)\right).
\end{equation*}
Again, multiplying by $\hat{A}_i$ can be done in $\mathcal{O}(m_a n)$ computational steps for dense matrices. For sparse constraint matrices, the computational effort is proportional to $\mathrm{nnz}(\hat{A}_i)$. Finally, multiplying with $2\overline{V}(v_i)$ from the left when evaluating the gradient~\eqref{eq:gradient_vi} can be carried out in $\mathcal{O}(kn)$.

To summarize, the computational effort for evaluating the objective in~\eqref{eq:subproblem} and its gradient~\eqref{eq:gradient_vi} is $\mathcal{O}(kn + m_a n)$ for dense constraint matrices and $\mathcal{O}(kn + \mathrm{nnz}(\hat{A}_i))$ for sparse constraint matrices. Therefore, in the sparse case with a large number of constraints (where $k \approx n$ holds), the effort to evaluate the objective and the gradient scales \emph{linearly} in the number of constraints. Note that according to~\eqref{eq:our_k}, the number of variables $k$ in each subproblem is at most $n$, which is independent of the number of constraints once $k = n$ is reached.

\begin{algorithm}[ht]
\caption{Augmented Mixing Method (high-level description)}
\label{alg:augmented_mixing}
Initialize each column $v_i$ of $V = (v_1 \mid \cdots \mid v_n) \in \mathbb{R}^{k \times n}$ uniformly at random on the unit sphere; set $y^a = 0$, $y^b = 0$, and $\mu = \sqrt{n}$\;
\Repeat{stopping criterion is met}{
    \For{$i \leftarrow 1$ \KwTo $n$}{
        Update column $v_i$ by approximately solving subproblem~\eqref{eq:subproblem} \\
        until the stopping criterion~\eqref{eq:subproblemstop} is satisfied\;
    }
    Update the dual variables $y^a$ and $y^b$ using~\eqref{eq:dual_update}\;
    Update the penalty parameter $\mu$ using~\eqref{eq:penalty_update}\;
}
\end{algorithm}

\section{Extensions, Implementation Details, and Software}\label{sec:further_ingredients}

This section outlines several extensions and implementation details of the Augmented Mixing Method. In particular, we discuss data scaling, update order strategies, stopping criteria, support for multi-block SDPs, high-precision computation via warm-starting, and configurable parameters available to the user. We also showcase two small examples that our method has difficulty solving because it can get stuck in a local maximum or exhibiting extremely slow convergence.

\subsection{Scaling}\label{sec:scaling}

In augmented Lagrangian methods, it is essential to scale the problem data to ensure that the algorithm balances progress between minimizing the objective and improving feasibility. We therefore apply the following automatic scaling procedure to the input data.

First, each constraint matrix, as well as the cost matrix $C$, is normalized by its Frobenius norm:
\begin{equation*}
\tilde{A}_i = \frac{A_i}{\lVert A_i \rVert_F}, \quad i = 1,\ldots,m_a, \qquad 
\tilde{B}_i = \frac{B_i}{\lVert B_i \rVert_F}, \quad i = 1,\ldots,m_b, \qquad 
\tilde{C} = \frac{C}{\lVert C \rVert_F}.
\end{equation*}
Next, the right-hand side vectors are scaled accordingly:
\begin{equation*}
\bar{a}_i = \frac{a_i}{\lVert A_i \rVert_F}, \quad i = 1,\ldots,m_a, \qquad 
\bar{b}_i = \frac{b_i}{\lVert B_i \rVert_F}, \quad i = 1,\ldots,m_b.
\end{equation*}
Finally, the resulting vectors are normalized to unit $\ell_2$-norm:
\begin{equation*}
\tilde{a} = \frac{\bar{a}}{\lVert \bar{a} \rVert_2}, \qquad 
\tilde{b} = \frac{\bar{b}}{\lVert \bar{b} \rVert_2}.
\end{equation*}

The Augmented Mixing Method then solves the scaled problem defined by the cost matrix $\tilde{C}$, the constraint matrices $\tilde{A}_1,\ldots,\tilde{A}_{m_a}, \tilde{B}_1,\ldots,\tilde{B}_{m_b}$, and the right-hand side vectors $\tilde{a}, \tilde{b}$. The computed solution is subsequently rescaled to yield an approximate solution to the original problem.

\subsection{Order of column updates} \label{sec:order_of_updates}

By default, the Augmented Mixing Method performs cyclic updates over the columns $v_1,\ldots,v_n$ of~$V$, as described in Section~\ref{sec:column_updates}. However, this update order can be modified. For instance, the paper~\cite{erdogdu2022convergence} on block coordinate descent for the Max-Cut relaxation~\eqref{eq:max_cut_sdp} explores several update strategies, including dynamic heuristics that select the next column at runtime based on recent progress.

In our setting, such approaches are less straightforward due to the augmented Lagrangian formulation, where the objective functions change over the course of the algorithm. Nevertheless, results from the literature on multi-block ADMM and block coordinate descent suggest that randomized update orders can improve convergence behavior, sometimes even ensuring convergence in expectation, whereas deterministic schemes may fail~\cite{cipolla2023alinear,mihic2021managing,sun2020ontheefficiency}. The paper~\cite{sun2020ontheefficiency} also discusses double-sweep strategies (i.e., forward and backward passes) which can improve practical convergence in some cases.

Our implementation supports these variants through two boolean parameters: \texttt{shuffling}, which randomizes the column order in each outer iteration, and \texttt{double\_sweep}, which enables forward and reverse sweeps. 
In practice, we observe that a fixed cyclic order tends to yield the best performance in terms of iteration count and runtime for most problem instances. 
Hence, both parameters are disabled by default and were not used in the numerical experiments reported in Section~\ref{sec:computational_experiments}.
However, our software allows users to activate these alternative strategies as needed.

\subsection{Stopping criteria} \label{sec:stopping_criteria}

In the following, we discuss the stopping criteria implemented in the Augmented Mixing Method. Note that the stopping criteria are applied to the scaled version of the SDP; see Section~\ref{sec:scaling}. This means that the specific values of the error measures presented below will differ slightly from those on the user’s original (unscaled) SDP. To simplify notation, we write $C$ instead of $\tilde{C}$ and so on. 

At any point during the Augmented Mixing Method, we have access to the current value of the primal matrix variable $X = V^\top V$, which is positive semidefinite by construction, as well as approximate dual variables $y^a, y^b$ with $y^b \geq 0$. An approximation of the dual slack matrix $Z$ is not directly available and must be computed by projecting $C - \mathcal{A}^\top(y^a) - \mathcal{B}^\top(y^b)$ onto the cone of positive semidefinite matrices:
\begin{equation}\label{eq:compute_Z}
    Z = [C - \mathcal{A}^\top(y^a) - \mathcal{B}^\top(y^b)]_+.
\end{equation}

Our stopping criteria are based on the KKT conditions for the primal-dual pair~\eqref{eq:primal_sdp} and~\eqref{eq:dual_sdp}, analogous to those used in \texttt{MOSEK}~\cite{mosek}. For given $X, Z \in \mathcal{S}_+^n$, $y^a \in \mathbb{R}^{m_a}$, and $y^b \in \mathbb{R}_+^{m_b}$, we define the error measures as
\begin{equation}\label{eq:stopping_criteria}
\begin{aligned}
\texttt{pinf} &\coloneqq \frac{
    \max\left\{ 
        \lVert a - \mathcal{A}(X) \rVert_\infty,\,
        \lVert [b - \mathcal{B}(X)]_+ \rVert_\infty
    \right\}
}{
    1 + \max\left\{ \lVert a \rVert_\infty,\, \lVert b \rVert_\infty \right\}
}, \\
\texttt{gap} &\coloneqq \frac{
    \left| \langle C, X \rangle - ( \langle a, y^a \rangle + \langle b, y^b \rangle ) \right|
}{
    1 + \left| \langle C, X \rangle \right| + \left| \langle a, y^a \rangle + \langle b, y^b \rangle \right|
}, \\
\texttt{dinf} &\coloneqq \frac{
    \left\lVert C - \mathcal{A}^\top(y^a) - \mathcal{B}^\top(y^b) - Z \right\rVert_\infty
}{
    1 + \lVert C \rVert_\infty
}, \\
\texttt{compl} &\coloneqq \frac{
    \langle X, Z \rangle
}{
    1 + \left| \langle C, X \rangle \right| + \left| \langle a, y^a \rangle + \langle b, y^b \rangle \right|
}.
\end{aligned}
\end{equation}
The Augmented Mixing Method terminates if
\begin{equation}\label{eq:stopping_tolerance}
    \max \{ \texttt{pinf}, \, \texttt{gap}, \, \texttt{dinf}, \, \texttt{compl} \} < \texttt{tol}
\end{equation}
holds, where \texttt{tol} is a user-specified parameter that defaults to~$10^{-12}$.

Note, however, that the error measures~\eqref{eq:stopping_criteria} are only required for termination and are not needed during the execution of the algorithm. For example, the dynamic update strategy for the penalty parameter is based on a different quantity, namely the ratio defined in~\eqref{eq:penalty_ratio}. Moreover, computing $Z$ can be very time-consuming when using extended precision arithmetic. Therefore, we avoid computing the expensive projection in~\eqref{eq:compute_Z} too frequently and only compute $Z$ when the solver appears close to convergence.

To implement this, the computation of $Z$ is triggered only every \texttt{iters\_Z} outer iterations of the Augmented Mixing Method and only if
\[
\max \{ \texttt{pinf}, \, \texttt{gap}, \, \texttt{compl}^* \} < \texttt{tol}
\]
holds, where $\texttt{compl}^*$ is an approximation of \texttt{compl} defined as
\[
\texttt{compl}^* \coloneqq \frac{
    \lvert \langle X, C - \mathcal{A}^\top(y^a) - \mathcal{B}^\top(y^b) \rangle \rvert
}{
    1 + \left| \langle C, X \rangle \right| + \left| \langle a, y^a \rangle + \langle b, y^b \rangle \right|
}.
\]

The parameter \texttt{iters\_Z} can be set by the user and has a default value of $50$. Moreover, users can also set a time limit or a maximum number of outer iterations to terminate the method early; see Section~\ref{sec:parameters} for a full list of available parameters.

We remark that later, in Section~\ref{sec:computational_experiments}, we compute these measures for the unscaled SDP when presenting results in figures and tables. This ensures a fair comparison across all tested solvers, independent of their internal stopping criteria or scaling strategies.

\subsection{Handling multi-block SDPs}\label{sec:multi-block}

The Augmented Mixing Method also supports SDPs with multiple positive semidefinite matrix variables, i.e., problems of the form
\begin{equation*}
\begin{array}{rl}
\text{minimize} \quad & \displaystyle \sum_{i=1}^q \langle C_i, X_i \rangle \\[1.2ex]
\text{subject to} \quad & \displaystyle \sum_{i=1}^q \mathcal{A}_i(X_i) = a, \\
& \displaystyle \sum_{i=1}^q \mathcal{B}_i(X_i) \geq b, \\
& X_i \in \mathcal{S}_+^{n_i}, \quad i = 1,\ldots, q. \\
\end{array}
\end{equation*}
Here, $q \in \mathbb{N}$ denotes the number of blocks, and $n_1,\ldots,n_q \in \mathbb{N}$ are the corresponding block sizes. For such SDPs, the Augmented Mixing Method maintains a factorization $X_i = V_i^\top V_i$ for each block $i \in \{1,\ldots,q\}$, where each rectangular matrix $V_i$ has dimensions $k_i \times n_i$. Similar to~\eqref{eq:our_k}, the rank parameter $k_i$ is chosen as
\begin{equation*}
k_i = \min\left(n_i, \left\lceil \sqrt{2(m_a + m_b)} \right\rceil\right)
\end{equation*}
for all $i = 1,\ldots,q$. Each outer iteration of the Augmented Mixing Method then consists of sequentially updating every column of every matrix $V_i$.

The initial penalty parameter is set to $\mu = \sqrt{\max\{ n_1,\ldots,n_q \}}$. Updating the penalty parameter and scaling the problem is handled analogously to the one-block case, and the error measures~\eqref{eq:stopping_criteria} are defined accordingly. Numerical results for multi-block SDPs are reported in Section~\ref{sec:rand}.

\subsection{Simple examples where the method fails or stagnates}\label{sec:fail}

Our method does not successfully solve all semidefinite programs. Even in cases where Slater’s condition holds (i.e. both the primal and dual problems are strictly feasible), convergence is not guaranteed. Below, we outline two primary obstacles we have identified.

\paragraph{Getting stuck in a local maximum.}
Consider the following strictly feasible SDP involving two blocks of size $1 \times 1$ and its Burer--Monteiro factorization counterpart with $x_i = v_i^2$:
\begin{equation}\label{eq:counter_example_multi_block}
\begin{array}{rlcrl}
\text{minimize} \quad & x_1 & = & \text{minimize} \quad & v_1^2\\[1.2ex]
\text{subject to} \quad & x_1+x_2=2,&& \text{subject to} \quad & v_1^2+v_2^2=2,\\
&x_2 = 1, &&&v_2^2=1,\\
& x_1,x_2\in\mathcal{S}_+^1,&&& v_1,v_2\in\mathbb{R}.
\end{array}
\end{equation}
The corresponding augmented Lagrangian is given by
\begin{equation*}
    \mathcal{L}(v_1,v_2,y;\mu) = v_1^2 + y_1(2-v_1^2-v_2^2)+y_2(1-v_2^2)+\frac{\mu}{2}(2-v_1^2-v_2^2)^2+\frac{\mu}{2}(1-v_2^2)^2.
\end{equation*}
Suppose the current iterate consists of $v_1 = 0$, $v_2 = \sqrt{1.5}$, and $y = (2,-2)^\top$. In this case, neither block update changes the values of $v_1$ or $v_2$:
\begin{align*}
    \mathcal{L}(v_1, \sqrt{1.5},y;\mu) &=  v_1^2 + y_1(0.5-v_1^2)-0.5y_2+\frac{\mu}{2}(0.5-v_1^2)^2+\frac{\mu}{2}(0.5)^2\\
    &= 2 - v_1^2 + \frac{\mu}{2}(0.5 - v_1^2)^2 + \frac{\mu}{2}(0.5)^2,\\
    \mathcal{L}(0,v_2,y;\mu) &= y_1(2-v_2^2)+y_2(1-v_2^2)+\frac{\mu}{2}(2-v_2^2)^2+\frac{\mu}{2}(1-v_2^2)^2\\
    &= 1 + \frac{\mu}{2}(2-v_2^2)^2 + \frac{\mu}{2}(1-v_2^2)^2.
\end{align*}
The first block remains at a local maximum of the augmented Lagrangian, and the second at a local minimum. The dual variable $y$ keeps getting updated in the direction $(+1, -1)^\top$, and the penalty parameter $\mu > 0$ has no effect on this behavior. Thus, this example shows that, without additional assumptions, no update rule for~$\mu$ can guarantee that the algorithm solves every (primal-dual strictly feasible) SDP instance. The augmented Lagrangian used in the Augmented Mixing Method can admit non-optimal first-order critical points even when~$n=k$, and the algorithm may get stuck at such points. While this example is initialized to trigger failure, we observe that similar stagnation often occurs, though not always, when starting from random initializations. Notably, this issue is not caused by Problem~\eqref{eq:counter_example_multi_block} having multiple blocks: combining the two blocks into a single $2\times 2$ matrix variable does not affect the outcome. This type of stagnation, where iterates do not improve in feasibility, is frequently observed in SDPs whose objective depends only on a single $1\times1$ block. Such structures commonly arise in polynomial optimization relaxations, particularly in the moment-sums-of-squares hierarchy. 

A simple workaround is to randomly perturb the current iterate slightly before solving the inner problem. While such a random perturbation can prevent complete stalling on a given instance, we found that it did not enable the solver to succeed on any additional instances from our test set (see Section~\ref{sec:computational_experiments}) that it could not already solve.

\paragraph{Slow convergence.}
The Augmented Mixing Method performs more than one million iterations to solve the following SDP to a relative accuracy of $10^{-8}$:
\begin{equation*}
\begin{array}{rl}
\text{minimize} \quad & \left\langle X_1, \left(\begin{smallmatrix}
  1 & 1 & 0\\ 1 & 1 & 0\\ 0&0&0
\end{smallmatrix} \right) \right\rangle + x_2\\[1.2ex]
\text{subject to} \quad & \operatorname{trace}(X_1) + x_2 = 4, \\
& \langle X_1, J_3\rangle + x_2 = 4,\\
& X_1\in \mathcal{S}^3_+, \quad x_2 \in \mathbb{R}_+.
\end{array}
\end{equation*}
In this example, the primal iterates initially converge rapidly toward feasibility, but progress slows down significantly in later iterations. Our penalty parameter update strategy targets a geometric reduction in the primal residual by maintaining \texttt{ratio} close to $1$. However, as primal feasibility improves, this strategy continues to increase the penalty parameter $\mu$, even when such increases are no longer beneficial. As a result, the primal residual is driven even smaller than necessary, leading to extremely small dual updates.

Note that the (original) mixing method~\cite{wang2017mixing} for the Goemans--Williamson SDP relaxation of the Max-Cut problem can already exhibit arbitrarily slow linear convergence near the optimum. Explicit examples are given in~\cite{Eelkema2022}, based on instances whose objective matrix has a nonzero eigenvalue close to zero. The Augmented Mixing Method exhibits the same behavior on these instances. This appears to be caused by the block coordinate descent approach, since \texttt{SDPLR}~\cite{burer2003nonlinear} converges quickly on them to a moderate-precision solution.

\subsection{Arbitrary-precision computation via warm-starting} \label{sec:warmstart}

Most semidefinite programming solvers rely on fixed-precision arithmetic, primarily because modern CPU architectures natively support certain floating-point types and because highly optimized BLAS and LAPACK routines are available for these types. Typically, double precision is used, which often limits the maximum achievable relative primal-dual accuracy of state-of-the-art interior-point methods to around~$10^{-10}$. However, some applications require significantly higher accuracy; see, e.g., the recent work~\cite{delaat2024optimalityuniquenessd4root} on the famous kissing number problem.

Our implementation of the Augmented Mixing Method is fully parametrized and can operate with any floating-point type supported by the Julia programming language. To efficiently compute arbitrary-precision solutions, we employ a two-stage warm-starting strategy:
\begin{enumerate}
    \item \textbf{Initial solve using double precision.} We first solve the SDP using the \texttt{Float64} data type with the tolerance parameter \texttt{tol} set to~$10^{-12}$.
    \item \textbf{Refinement using extended precision.} We promote the solution to the extended precision type. The solver resumes the optimization process from this warm-start point until the user-specified tolerance is reached.
\end{enumerate}
 
A key advantage of our implementation compared to other arbitrary-precision solvers is that the most computationally expensive operations, namely (sparse) matrix-vector multiplications, are relatively straightforward to optimize even for extended precision, unlike the more demanding computations required by interior-point methods.

In Section~\ref{sec:arbitrary_precision}, we present numerical results on randomly generated SDPs where the accuracy parameter \texttt{tol} in the refinement phase is set to~$10^{-20}$.

\subsection{Solver usage and parameters}\label{sec:parameters}

Our implementation of the Augmented Mixing Method is open-source and written entirely in the Julia programming language. It is designed to be easily configurable and extensible, enabling users to apply the solver to a broad class of SDPs with varying structure and precision requirements.

Using the solver in practice requires constructing an object of type \texttt{SdpData}, which encodes the problem data in the format expected by the Augmented Mixing Method. For a one-block SDP, this can be done via
\begin{center}
\texttt{sdp = SdpData(A, b, C, ineq\_start)},
\end{center}
where \texttt{A} is a list of symmetric constraint matrices, \texttt{b} is the right-hand side vector, and \texttt{C} is the symmetric cost matrix. The matrices in \texttt{A} must be ordered such that all equality constraints appear first, followed by the inequality constraints. The integer \texttt{ineq\_start} specifies the index of the first inequality constraint. If the SDP has no inequality constraints, the argument \texttt{ineq\_start} must be set to \texttt{length(b) + 1}. The implementation supports both dense and sparse matrices. All problem data must share the same element type (e.g., \texttt{Float64}). For multi-block SDPs, the lists \texttt{A} and \texttt{C} must include an additional leading index to indicate the corresponding block.

To solve the SDP with default parameters, call
\begin{center}
\texttt{X, y, Z, warm\_start = augmented\_mixing(sdp)}.
\end{center}
This returns the computed primal-dual triplet $(X, y, Z)$ along with a \texttt{WarmStart} object. To refine the solution using higher precision, first construct a new \texttt{SdpData} object \texttt{sdp\_T} with the desired floating-point type \texttt{T}, and continue the optimization using
\begin{center}
\texttt{augmented\_mixing(sdp\_T; tol=T(prec), warm\_start=WarmStart(warm\_start, T))},
\end{center}
where \texttt{prec} denotes the target tolerance.

An overview of all configurable parameters for \texttt{augmented\_mixing()} is listed in Table~\ref{table:parameters}.

\begin{table}[ht]
\centering
\resizebox{\textwidth}{!}{%
\setlength{\tabcolsep}{7pt}
\renewcommand{\arraystretch}{1.1}
\begin{tabular}{@{}llr p{7cm}@{}}
\toprule
\textbf{Parameter} & \textbf{Type} & \textbf{Default} & \textbf{Description} \\
\midrule
\texttt{tol}              & \texttt{T}                   & \texttt{1e-12}               & stopping tolerance; see~\eqref{eq:stopping_tolerance} \\
\texttt{mu_start}        & \texttt{T}                   & $\sqrt{\max\{n_1,\dots,n_q\}}$ & initial penalty parameter $\mu$; defaults to the square root of the largest block size\\
\texttt{time_limit}      & \texttt{Float64}             & \texttt{typemax(Float64)}    & maximum allowed wall-clock time in seconds \\
\texttt{max_iters}        & \texttt{Int}                 & \texttt{typemax(Int)}        & maximum number of outer iterations \\
\texttt{iters_Z}     & \texttt{Int}                 & \texttt{50}                  & $Z$ is computed every \texttt{iters_Z} iterations if all other errors are below \texttt{tol}; see Section~\ref{sec:stopping_criteria} \\
\texttt{scaling}          & \texttt{Bool}                & \texttt{true}                & whether automatic scaling is applied; see Section~\ref{sec:scaling} \\
\texttt{shuffling}        & \texttt{Bool}                & \texttt{false}               & whether column update order is randomized; see Section~\ref{sec:order_of_updates} \\
\texttt{double_sweep}    & \texttt{Bool}                & \texttt{false}               & whether columns are updated in both forward and reverse order; see Section~\ref{sec:order_of_updates} \\
\texttt{p}                & \texttt{T}                   & \texttt{1.0}                 & dual step size; see Section~\ref{sec:update} \\
\texttt{warm_start}      & \texttt{WarmStart\{\textit{T}\}} & \texttt{nothing}        & warm-start object; obtained from previous call to \texttt{augmented\_mixing()} \\
\texttt{delta}      & \texttt{T} & \texttt{0.01}        & relative tolerance for solving subproblems; see~\eqref{eq:subproblemstop} \\
\texttt{epsilon}      & \texttt{T} & \texttt{0.01}        & absolute tolerance for solving subproblems; see~\eqref{eq:subproblemstop} \\
\texttt{max_evals}      & \texttt{Int} & \texttt{1000}        & maximum number of function and gradient evaluations per column update \\
\texttt{tau}      & \texttt{T} & \texttt{1.03}        & factor for updating the penalty parameter; see~\eqref{eq:penalty_update} \\
\texttt{rat_min}      & \texttt{T} & \texttt{0.8}        & lower bound for \texttt{ratio}; see~\eqref{eq:penalty_update} \\
\texttt{rat_max}      & \texttt{T} & \texttt{1.2}        & upper bound for \texttt{ratio}; see~\eqref{eq:penalty_update} \\
\bottomrule
\end{tabular}
}
\caption{Overview of configurable parameters for \texttt{augmented_mixing()}. The type \texttt{T} must be a floating-point type, e.g., \texttt{Float64}.}
\label{table:parameters}
\end{table}

\section{Computational Experiments}\label{sec:computational_experiments}

In this section, we present extensive numerical results that demonstrate the performance of the Augmented Mixing Method using its default settings. We compare the Augmented Mixing Method with several state-of-the-art SDP solvers in terms of runtime and solution accuracy. We focus primarily on instances with a large number of inequality constraints, as our method is specifically designed for this setting.

\subsection{Computational setup}

We conducted all experiments on a machine running Debian~12, equipped with an AMD EPYC Genoa processor clocked at 3.59\,GHz and 1024\,GB of RAM. In all experiments, the number of computation threads was limited to four. Our implementation of the Augmented Mixing Method is written in the Julia programming language~\cite{bezanson2017julia}, and all experiments were performed using Julia version~1.11.5.

To solve the subproblems~\eqref{eq:subproblem}, we use the L-BFGS algorithm provided by the \texttt{Optim.jl} package~\cite{mogensen2018optim}. For computations with extended floating-point precision, as described in Section~\ref{sec:arbitrary_precision}, we use the \texttt{Double64} data type provided by the \texttt{DoubleFloats.jl} package~\cite{sarnoff2022doublefloats}. All optimization problems for solvers other than the Augmented Mixing Method were modeled using \texttt{JuMP}~\cite{lubin2023jump}. The source code of the Augmented Mixing Method is available as a Julia package at~\url{https://github.com/jschwiddessen/AugmentedMixing.jl}. All code and data required to reproduce the computational experiments presented in this paper are included in the supplementary material on the arXiv page of this paper and at~\url{https://zenodo.org/records/15976859}.

\paragraph{Solver selection and configuration.}
In our experiments using double-precision arithmetic, we compare against the state-of-the-art interior-point solver \texttt{MOSEK}~\cite{mosek} (version~11.0.20), the low-rank factorization method \texttt{SDPLR} (version~1.03-beta) introduced in~\cite{burer2003nonlinear,burer2005local}, and the solver \texttt{SCS}~\cite{brendan2016conic} (version~3.2.7) for large-scale convex cone programs. However, we exclude some solvers for specific problem classes due to limitations in handling large numbers of inequality constraints. These exclusions are noted in the corresponding subsections.

We use the default parameters of each solver with the following modifications. For every solver, including the Augmented Mixing Method, we impose a time limit of two hours. We restrict \texttt{MOSEK} to four threads and set the parameters \texttt{eps\_abs} and \texttt{eps\_rel} of \texttt{SCS} to~$10^{-6}$. For the Augmented Mixing Method, we additionally cap the number of outer iterations at \num{100000}. In addition to the default setting of \texttt{tol=1e-12}, we also run the Augmented Mixing Method with the relaxed stopping tolerance \texttt{tol=1e-6} to investigate the trade-off between accuracy and runtime.

For the extended-precision experiments in Section~\ref{sec:arbitrary_precision}, we compare the Augmented Mixing Method to the interior-point solver \texttt{Hypatia}~\cite{coey2022solving} (version~0.8.2), which is also implemented in Julia. As before, we impose a two-hour time limit and set a maximum of \num{100000} outer iterations for the Augmented Mixing Method. Additionally, we set the parameter \texttt{tol\_rel\_opt} of \texttt{Hypatia} to~$10^{-20}$, and the parameter \texttt{tol} of the Augmented Mixing Method to~$10^{-20}$ after warm-starting. We provide \texttt{Hypatia} with the dual formulation of the SDP, as this is the solver’s preferred input form.

\paragraph{Format and interpretation of results.}

For each class of SDP instances, we provide detailed tables summarizing the performance of the Augmented Mixing Method using the default tolerance \texttt{tol=1e-12}, including all four primal-dual error metrics from~\eqref{eq:stopping_criteria}. If only one additional solver is included in an experiment, we also report that solver’s computation time and the maximum of its four primal-dual errors in the same table. In other cases, we use performance profiles~\cite{dolan2002benchmarking} to compare runtime and solution accuracy. These plots also include results of the Augmented Mixing Method with \texttt{tol=1e-6}, which are not shown in the tables.

All primal-dual errors are computed from the solutions returned by the respective solvers, using the error metrics defined in~\eqref{eq:stopping_criteria} and evaluated on the original, unscaled SDP. This ensures a fair comparison across solvers, independent of their internal stopping criteria or scaling strategies.

In the tables, the first four columns list the instance name, the matrix size~$n$, and the number of equality and inequality constraints, $m_a$ and $m_b$. The next three columns show the Augmented Mixing Method's return status (\texttt{tol} if all primal-dual errors on the internally scaled problem fall below the specified tolerance \texttt{tol}; \texttt{iter} if the maximum number of outer iterations is reached; \texttt{time} if the time limit is exceeded), the elapsed wall-clock time in seconds, and the number of outer iterations performed. The final four columns report the primal-dual error metrics defined in~\eqref{eq:stopping_criteria}. In some cases, two additional columns report the runtime and maximum error of another solver for comparison. For randomly generated SDP instances, the columns for $n$, $m_a$, and $m_b$ are omitted, since these values are evident from the instance names. In Tables~\ref{table:max_cut_met} and~\ref{table:srflp_met}, we also summarize runtime results using the geometric mean.

The performance profile figures show two plots side by side. The left plot evaluates cumulative runtime performance: for each solver, it shows the fraction of instances for which the runtime is within a given factor ($x$-axis, $\log$ scale) of the best runtime achieved on that instance. The right plot evaluates solution accuracy: it shows the fraction of instances where the maximum primal-dual error is within a given factor of the smallest such error across all solvers. In both cases, curves that rise more steeply and reach $1$ more quickly indicate better performance. All performance profiles presented in this section were generated using the Julia package \texttt{BenchmarkProfiles.jl}~\cite{orbanbenchmarkprofiles2023}.

\subsection{Max-Cut problem}

The Max-Cut problem in a simple undirected graph~$G$ with $n$~vertices and weighted adjacency matrix~$A$ can be formulated as the quadratic unconstrained binary optimization problem
\begin{equation}\label{eq:max_cut}
\begin{array}{rl}
\text{maximize} \quad & \frac{1}{4} x^\top L x \\[1.2ex]
\text{subject to} \quad & x \in \{-1,1\}^n, \\
\end{array}
\end{equation}
where $L = \operatorname{Diag}(A e_n) - A$ denotes the Laplacian matrix of $G$. For $C = \frac{1}{4} L $, the basic SDP relaxation of~\eqref{eq:max_cut} is given by~\eqref{eq:max_cut_sdp}. This relaxation can be solved very efficiently by specialized low-rank methods such as the mixing method~\cite{wang2017mixing} or its extensions~\cite{chen2023burermonteiro,erdogdu2022convergence,kim2021momentuminspired}, as well as algorithms from the \texttt{Manopt} package~\cite{boumal2014manopt}. Since the Augmented Mixing Method is not competitive for this specific problem class, we omit computational results for the basic SDP relaxation of Max-Cut.

However, it is well known that the basic relaxation~\eqref{eq:max_cut_sdp} can be strengthened by adding the following $4 \binom{n}{3}$ triangle inequalities
\begin{equation}\label{eq:triangle_inequalities}
    \begin{split}
        X_{ij} + X_{ik} + X_{jk} &\geq -1, \quad 1 \leq i < j < k \leq n, \\
        X_{ij} - X_{ik} - X_{jk} &\geq -1, \quad 1 \leq i < j < k \leq n, \\
        -X_{ij} + X_{ik} - X_{jk} &\geq -1, \quad 1 \leq i < j < k \leq n, \\
        -X_{ij} - X_{ik} + X_{jk} &\geq -1, \quad 1 \leq i < j < k \leq n, \\
    \end{split}
\end{equation}
which are facet-defining for the Max-Cut polytope~\cite{deza1997geometry}. Table~\ref{table:max_cut_met} shows numerical results of the Augmented Mixing Method on various instances of the Max-Cut SDP relaxation with $n$ between~60 and~250 in which all triangle inequalities~\eqref{eq:triangle_inequalities} were included at once. These instances are either directly taken from the \emph{Biq~Mac Library}~\cite{wiegele2007biqmaclibrary} or are generated in the same way using the graph generator \texttt{rudy}~\cite{rinaldi1998rudy}. To enhance numerical stability and reduce the effect of scaling, we set the main diagonal in all SDPs to zero, as it contributes only a constant to the objective. We report results for the Augmented Mixing Method with \texttt{tol=1e-12} and the \texttt{SCS} solver in Table~\ref{table:max_cut_met}. We exclude \texttt{MOSEK} and \texttt{SDPLR} here because neither is suited to handle the very large number of inequalities in these instances.

As the results in Table~\ref{table:max_cut_met} demonstrate, the Augmented Mixing Method can handle SDPs with extremely large numbers of inequality constraints while still producing primal-dual solutions with high accuracy. It typically requires fewer than \num{2000} outer iterations, with most instances being solved in under 1000 iterations. On instances with $n \in \{60,80 \}$, the \texttt{SCS} solver has a shorter runtime than the Augmented Mixing Method, but the opposite holds for larger instances. However, \texttt{SCS} generally returns solutions with maximum primal-dual errors between $10^{-6}$ and $10^{-7}$. In contrast, the primal-dual errors of the Augmented Mixing Method are almost always below $10^{-9}$ and often even smaller. Moreover, the method successfully terminates on two out of five instances with more than ten million inequality constraints. It reaches the time limit of two hours on the remaining three such large instances but still yields reasonably accurate solutions. On these large instances, \texttt{SCS} terminates only well after the time limit has already been exceeded and produces solutions with maximum primal-dual errors ranging roughly from $10^{-3}$ to $10^{-6}$.

Nevertheless, the Augmented Mixing Method may fail to converge within the prescribed time limit and maximum number of iterations in rare cases. In particular, some smaller instances of type `pw09' seem to pose difficulties, as the method reaches the time limit without producing even moderately accurate solutions. We suspect that this behavior is due to our scaling strategy. Moreover, note that our choice of $k$ in~\eqref{eq:our_k} yields $k=n$ for all instances considered here. We found that, for some instances, this choice is indeed necessary to avoid early stagnation of the Augmented Mixing Method, even when the SDP relaxation is tight and admits a rank-one solution.

\begin{table}[ht]
\centering
\resizebox{1.0\textwidth}{!}{
\setlength{\tabcolsep}{7pt}
\renewcommand{\arraystretch}{1.1}
\begin{filecontents*}{maxCutMet.csv}
instance,type,seed,statusCode,time,iter,tol,timeLimit,maxIters,itersZ,scaling,warmStart,k,n,numBlocks,avgBlockSize,minBlockSize,maxBlockSize,m,equations,inequalities,nnzA,densityAPercent,timeScaling,timeData,primalObj,dualObj,gap,pinf,dinf,compl,complNotPsdAbs,mu,evals,nnzY,gapOriginal,pinfOriginal,dinfOriginal,complOriginal,timeSCS,maxErrorSCS
g05_60.0,maxCutMet,940269413,tol,112.2,1600,1.0e-12,7200.0,100000,50,true,false,60,60,1,60.0,60,60,136940,60,136880,821340,0.16660581276471448,0.1,0.34080004692077637,-94.73754346371078,-94.73754346354539,4.9e-14,8.3e-13,7.6e-13,3.8e-14,1.3e-14,2.1e+01,1080800,790,8.7e-13,1.5e-10,6.5e-12,6.8e-13,58.9,7.8e-07
g05_80.0,maxCutMet,699298087,tol,123.8,750,1.0e-12,7200.0,100000,50,true,false,80,80,1,80.0,80,80,328720,80,328640,1971920,0.09373098685811633,0.6,0.8967299461364746,-144.236873186514,-144.23687318640475,1.6e-14,3.5e-13,5.8e-13,2.8e-14,4.6e-15,2.7e+01,630760,997,3.8e-13,1.0e-10,6.6e-12,6.6e-13,135.4,1.8e-06
g05_100.0,maxCutMet,486432763,tol,366.0,1250,1.0e-12,7200.0,100000,50,true,false,100,100,1,100.0,100,100,646900,100,646800,3880900,0.05999227083011285,0.8,2.7275590896606445,-204.49152302726588,-204.49152302795028,5.7e-14,8.1e-13,9.6e-14,3.8e-14,3.2e-14,3.9e+00,1081976,1187,1.7e-12,3.3e-10,1.4e-12,1.1e-12,312.3,1.7e-06
g05_120.0,maxCutMet,983590172,tol,454.4,900,1.0e-12,7200.0,100000,50,true,false,120,120,1,120.0,120,120,1123480,120,1123360,6740280,0.04166295795207747,1.8,5.627233982086182,-272.3211250205871,-272.3211250196715,4.9e-14,7.3e-13,1.1e-13,7.5e-15,1.9e-14,4.8e+00,874467,1407,1.7e-12,3.9e-10,1.9e-12,2.6e-13,638.6,7.1e-07
g05_140.0,maxCutMet,530550972,tol,601.0,650,1.0e-12,7200.0,100000,50,true,false,140,140,1,140.0,140,140,1790460,140,1790320,10742060,0.03061025019907095,5.4,9.438868999481201,-343.8069796194461,-343.8069796195459,3.6e-15,3.0e-13,1.7e-13,2.1e-14,9.8e-16,2.7e+01,856436,1670,1.4e-13,2.0e-10,3.4e-12,8.6e-13,1145.4,9.8e-08
g05_160.0,maxCutMet,210881114,tol,2761.6,2200,1.0e-12,7200.0,100000,50,true,false,160,160,1,160.0,160,160,2679840,160,2679680,16078240,0.02343633388560511,6.8,14.269469976425171,-428.06054052025286,-428.06054052022506,7.2e-16,2.7e-15,5.3e-13,8.5e-16,2.7e-17,3.6e+01,3040182,2068,3.2e-14,2.2e-12,1.2e-11,3.8e-14,2571.8,6.2e-08
g05_180.0,maxCutMet,305402420,tol,1455.3,800,1.0e-12,7200.0,100000,50,true,false,180,180,1,180.0,180,180,3823620,180,3823440,22940820,0.018517792040003974,7.9,16.37859797477722,-520.0995584072468,-520.0995584074586,4.1e-15,4.1e-13,1.9e-13,1.5e-14,7.0e-16,3.2e+01,1292195,2344,2.0e-13,4.0e-10,4.9e-12,7.6e-13,3475.4,7.6e-07
g05_200.0,maxCutMet,704552125,time,7200.2,3786,1.0e-12,7200.0,100000,50,true,false,200,200,1,200.0,200,200,5253800,200,5253600,31521800,0.014999524153945716,10.9,24.76013684272766,-604.7402716580273,-604.7402725461668,1.3e-11,3.9e-11,Inf,Inf,1.8e-12,1.4e+01,4797977,2352,7.3e-10,3.7e-08,3.1e-09,1.7e-10,5272.6,9.4e-07
g05_250.0,maxCutMet,888356939,tol,4847.0,1000,1.0e-12,7200.0,100000,50,true,false,250,250,1,250.0,250,250,10292250,250,10292000,61752250,0.009599805679030338,23.3,52.663142919540405,-858.9095090692487,-858.9095090696006,3.0e-15,8.8e-13,2.5e-13,1.2e-14,7.0e-15,1.3e+01,1764072,2952,2.0e-13,1.4e-09,8.9e-12,8.1e-13,9235.1,1.7e-03
pw05_60.0,maxCutMet,965201257,tol,257.0,4050,1.0e-12,7200.0,100000,50,true,false,60,60,1,60.0,60,60,136940,60,136880,821340,0.16660581276471448,0.1,0.30840396881103516,-594.4999999998713,-594.5000000037514,1.9e-13,1.4e-13,0.0e+00,1.8e-13,1.8e-13,5.9e-01,2340761,101734,3.3e-12,2.1e-11,5.5e-15,3.2e-12,25.1,2.1e-08
pw05_80.0,maxCutMet,252327208,tol,428.0,4100,1.0e-12,7200.0,100000,50,true,false,80,80,1,80.0,80,80,328720,80,328640,1971920,0.09373098685811633,0.4,0.8189878463745117,-947.3489811860644,-947.3489811911984,1.2e-13,9.9e-13,8.3e-14,1.4e-14,4.2e-15,2.3e+00,2007840,991,2.7e-12,2.9e-10,2.1e-12,3.1e-13,118.4,2.7e-06
pw05_100.0,maxCutMet,553395641,tol,204.3,650,1.0e-12,7200.0,100000,50,true,false,100,100,1,100.0,100,100,646900,100,646800,3880900,0.05999227083011285,0.9,1.892416000366211,-1383.4254504118119,-1383.425450408855,3.9e-14,4.1e-13,4.1e-13,2.6e-14,1.7e-15,2.7e+01,652619,1179,1.1e-12,1.6e-10,1.3e-11,7.0e-13,295.1,1.1e-06
pw05_120.0,maxCutMet,494034228,tol,402.5,750,1.0e-12,7200.0,100000,50,true,false,120,120,1,120.0,120,120,1123480,120,1123360,6740280,0.04166295795207747,1.6,3.5167791843414307,-1828.743613963565,-1828.7436139629217,5.4e-15,2.7e-13,6.3e-13,1.2e-14,3.7e-16,3.1e+01,887978,1510,1.8e-13,1.5e-10,2.4e-11,4.0e-13,635.7,1.2e-06
pw05_140.0,maxCutMet,444031489,tol,637.7,750,1.0e-12,7200.0,100000,50,true,false,140,140,1,140.0,140,140,1790460,140,1790320,10742060,0.03061025019907095,3.0,6.6043479442596436,-2344.9834854687765,-2344.9834854606565,4.7e-14,7.4e-13,4.0e-13,4.8e-14,3.1e-15,3.0e+01,1000780,1804,1.7e-12,4.9e-10,1.8e-11,1.7e-12,1139.9,5.2e-07
pw05_160.0,maxCutMet,527747595,tol,788.6,600,1.0e-12,7200.0,100000,50,true,false,160,160,1,160.0,160,160,2679840,160,2679680,16078240,0.02343633388560511,5.8,9.699463844299316,-2865.957141402115,-2865.957141324884,3.2e-13,9.7e-13,8.1e-13,7.3e-14,6.4e-15,2.7e+01,917260,1913,1.3e-11,7.9e-10,4.1e-11,3.0e-12,1843.4,5.3e-07
pw05_180.0,maxCutMet,904693368,tol,2021.0,1450,1.0e-12,7200.0,100000,50,true,false,180,180,1,180.0,180,180,3823620,180,3823440,22940820,0.018517792040003974,7.0,16.952823877334595,-3529.533914360478,-3529.5339143555057,1.5e-14,8.0e-13,1.3e-13,5.9e-15,9.2e-15,6.4e+00,1686294,2181,7.0e-13,7.8e-10,7.7e-12,2.7e-13,2797.3,9.3e-07
pw05_200.0,maxCutMet,202878517,tol,1908.7,750,1.0e-12,7200.0,100000,50,true,false,200,200,1,200.0,200,200,5253800,200,5253600,31521800,0.014999524153945716,12.8,26.182953119277954,-4124.988456239088,-4124.988456216741,5.3e-14,3.0e-13,9.0e-13,1.8e-14,3.9e-16,3.1e+01,1368339,2593,2.7e-12,3.5e-10,5.7e-11,9.2e-13,4000.7,1.2e-06
pw05_250.0,maxCutMet,794531307,time,7201.4,1650,1.0e-12,7200.0,100000,50,true,false,250,250,1,250.0,250,250,10292250,250,10292000,61752250,0.009599805679030338,28.7,62.969878911972046,-5827.212296700072,-5827.212302156091,7.4e-12,1.4e-11,Inf,Inf,2.9e-13,9.6e+00,2555120,3046,4.7e-10,1.9e-08,1.4e-09,3.5e-11,9225.7,6.8e-06
pw09_60.0,maxCutMet,958159902,tol,228.2,3450,1.0e-12,7200.0,100000,50,true,false,60,60,1,60.0,60,60,136940,60,136880,821340,0.16660581276471448,0.1,0.3067319393157959,-613.9999999998955,-614.0000000217333,8.0e-13,2.0e-13,0.0e+00,8.0e-13,8.0e-13,5.7e-01,2059202,102614,1.8e-11,3.0e-11,1.8e-14,1.8e-11,15.4,5.0e-07
pw09_80.0,maxCutMet,480919819,time,7200.1,43139,1.0e-12,7200.0,100000,50,true,false,80,80,1,80.0,80,80,328720,80,328640,1971920,0.09373098685811633,0.4,0.8807730674743652,-936.2366968073504,-941.5316210114422,9.4e-05,1.9e-04,Inf,Inf,1.8e-05,3.0e+00,37828286,1598,2.8e-03,5.4e-02,3.2e-02,2.9e-03,140.1,1.9e-06
pw09_100.0,maxCutMet,689647233,time,7200.3,23553,1.0e-12,7200.0,100000,50,true,false,100,100,1,100.0,100,100,646900,100,646800,3880900,0.05999227083011285,1.2,1.8302559852600098,-1351.4522379960074,-1353.845062072924,2.4e-05,1.2e-04,Inf,Inf,9.4e-06,3.6e+00,23372497,1951,8.8e-04,5.0e-02,2.2e-02,2.4e-03,248.8,1.2e-06
pw09_120.0,maxCutMet,743510841,time,7200.3,14290,1.0e-12,7200.0,100000,50,true,false,120,120,1,120.0,120,120,1123480,120,1123360,6740280,0.04166295795207747,2.1,3.356719970703125,-1795.2533052589274,-1802.4314135592524,4.6e-05,9.8e-05,Inf,Inf,1.6e-05,4.1e+00,15717580,2160,2.0e-03,5.2e-02,1.3e-02,2.3e-03,590.0,6.0e-07
pw09_140.0,maxCutMet,819596173,time,7200.4,9444,1.0e-12,7200.0,100000,50,true,false,140,140,1,140.0,140,140,1790460,140,1790320,10742060,0.03061025019907095,3.7,6.47329306602478,-2204.998289114688,-2209.282009151372,1.9e-05,4.7e-05,Inf,Inf,4.0e-06,5.2e+00,11288612,2054,9.7e-04,3.2e-02,1.5e-02,1.1e-03,845.6,2.5e-06
pw09_160.0,maxCutMet,134184283,tol,3825.4,3650,1.0e-12,7200.0,100000,50,true,false,160,160,1,160.0,160,160,2679840,160,2679680,16078240,0.02343633388560511,6.0,10.398859024047852,-2761.40907104475,-2761.409071046006,3.9e-15,2.3e-13,6.7e-13,9.3e-14,6.3e-15,5.7e+00,4432974,1942,2.3e-13,1.9e-10,4.5e-11,5.4e-12,1541.4,2.2e-06
pw09_180.0,maxCutMet,275544110,tol,2441.4,1650,1.0e-12,7200.0,100000,50,true,false,180,180,1,180.0,180,180,3823620,180,3823440,22940820,0.018517792040003974,8.7,17.860620975494385,-3330.313908490378,-3330.313908478458,2.8e-14,5.2e-13,6.6e-13,1.3e-13,1.2e-14,4.9e+00,2161745,2113,1.8e-12,5.1e-10,5.0e-11,8.4e-12,2577.4,6.4e-07
pw09_200.0,maxCutMet,498298642,time,7201.5,3228,1.0e-12,7200.0,100000,50,true,false,200,200,1,200.0,200,200,5253800,200,5253600,31521800,0.014999524153945716,12.8,27.470561027526855,-3939.342183820449,-3943.550091754988,7.5e-06,5.3e-05,Inf,Inf,1.1e-05,4.6e+00,5239369,4569,5.3e-04,5.5e-02,2.1e-02,2.5e-03,4423.9,6.6e-07
pw09_250.0,maxCutMet,994907493,time,7201.5,1464,1.0e-12,7200.0,100000,50,true,false,250,250,1,250.0,250,250,10292250,250,10292000,61752250,0.009599805679030338,29.7,66.78122210502625,-5581.990360121399,-5581.990012127163,3.6e-10,3.1e-09,Inf,Inf,1.2e-10,7.6e+00,2700590,3034,3.1e-08,4.9e-06,9.5e-07,1.1e-07,9210.2,4.4e-06
w05_60.0,maxCutMet,153501654,tol,107.8,1550,1.0e-12,7200.0,100000,50,true,false,60,60,1,60.0,60,60,136940,60,136880,821340,0.16660581276471448,0.1,0.34146904945373535,-702.1802753837262,-702.1802753816472,1.0e-13,7.9e-13,3.2e-13,3.6e-14,4.3e-15,2.2e+01,1044732,809,1.5e-12,1.5e-10,6.0e-12,5.2e-13,61.4,6.7e-07
w05_80.0,maxCutMet,277453054,tol,148.0,900,1.0e-12,7200.0,100000,50,true,false,80,80,1,80.0,80,80,328720,80,328640,1971920,0.09373098685811633,0.4,0.967034101486206,-1109.4647381886002,-1109.46473818629,5.6e-14,5.8e-13,3.9e-13,3.1e-14,3.4e-15,2.9e+01,790480,1065,1.0e-12,1.7e-10,9.6e-12,5.8e-13,145.3,2.1e-06
w05_100.0,maxCutMet,382926493,tol,234.8,750,1.0e-12,7200.0,100000,50,true,false,100,100,1,100.0,100,100,646900,100,646800,3880900,0.05999227083011285,1.6,2.8497068881988525,-1629.145144061904,-1629.1451440619583,7.4e-16,5.4e-13,4.6e-13,2.4e-14,1.6e-15,3.1e+01,736163,1213,1.7e-14,2.2e-10,1.4e-11,5.4e-13,296.7,2.2e-06
w05_120.0,maxCutMet,772547491,tol,491.3,900,1.0e-12,7200.0,100000,50,true,false,120,120,1,120.0,120,120,1123480,120,1123360,6740280,0.04166295795207747,2.6,4.560924053192139,-2162.013373289861,-2162.013373289718,1.2e-15,1.0e-14,1.9e-14,1.2e-15,2.0e-16,3.0e+01,1048410,1524,3.3e-14,5.4e-12,7.1e-13,3.3e-14,613.3,1.3e-06
w05_140.0,maxCutMet,327118557,tol,751.0,900,1.0e-12,7200.0,100000,50,true,false,140,140,1,140.0,140,140,1790460,140,1790320,10742060,0.03061025019907095,4.5,8.876321077346802,-2832.056742783418,-2832.056742782865,3.3e-15,8.9e-15,2.7e-13,8.4e-16,8.2e-18,3.3e+01,1152687,1784,9.8e-14,6.0e-12,1.2e-11,2.5e-14,1186.5,2.2e-06
w05_160.0,maxCutMet,889816777,tol,986.3,750,1.0e-12,7200.0,100000,50,true,false,160,160,1,160.0,160,160,2679840,160,2679680,16078240,0.02343633388560511,7.5,12.936612129211426,-3466.1639313237356,-3466.1639313074547,7.0e-14,6.3e-13,6.1e-13,3.6e-14,9.4e-16,3.4e+01,1125239,2063,2.3e-12,5.2e-10,3.0e-11,1.2e-12,1946.6,6.6e-07
w05_180.0,maxCutMet,122208585,tol,6702.1,4150,1.0e-12,7200.0,100000,50,true,false,180,180,1,180.0,180,180,3823620,180,3823440,22940820,0.018517792040003974,11.7,20.825779914855957,-4162.977390624843,-4162.977390624685,5.0e-16,1.0e-15,9.9e-13,2.2e-16,7.4e-18,4.9e+01,5971165,2404,1.9e-14,8.3e-13,5.5e-11,8.4e-15,2911.4,5.5e-07
w05_200.0,maxCutMet,992550569,tol,1967.4,750,1.0e-12,7200.0,100000,50,true,false,200,200,1,200.0,200,200,5253800,200,5253600,31521800,0.014999524153945716,15.5,30.955934047698975,-4933.558050991182,-4933.558050929306,1.5e-13,6.3e-13,4.8e-13,4.7e-14,1.5e-15,3.5e+01,1383487,2529,6.3e-12,7.2e-10,3.0e-11,1.9e-12,3676.9,1.2e-06
w05_250.0,maxCutMet,322566175,time,7203.2,1592,1.0e-12,7200.0,100000,50,true,false,250,250,1,250.0,250,250,10292250,250,10292000,61752250,0.009599805679030338,31.9,73.95794010162354,-7062.424447358648,-7062.42448606822,5.4e-11,8.9e-11,Inf,Inf,5.1e-13,1.5e+01,2556493,3103,2.7e-09,1.4e-07,5.6e-09,3.1e-10,9089.0,4.1e-05
w09_60.0,maxCutMet,30176789,tol,63.6,900,1.0e-12,7200.0,100000,50,true,false,60,60,1,60.0,60,60,136940,60,136880,821340,0.16660581276471448,0.1,0.3123490810394287,-956.4800343784425,-956.4800343727038,2.1e-13,9.4e-13,8.1e-13,7.7e-14,7.4e-15,2.5e+01,619357,760,3.0e-12,1.7e-10,2.0e-11,1.1e-12,58.6,1.3e-06
w09_80.0,maxCutMet,207282611,tol,121.8,800,1.0e-12,7200.0,100000,50,true,false,80,80,1,80.0,80,80,328720,80,328640,1971920,0.09373098685811633,0.6,0.9128589630126953,-1503.410280985087,-1503.4102809837998,2.3e-14,2.6e-13,3.5e-13,1.5e-14,6.5e-16,3.1e+01,641965,946,4.3e-13,7.5e-11,1.2e-11,2.8e-13,187.8,8.5e-07
w09_100.0,maxCutMet,403676044,tol,224.1,750,1.0e-12,7200.0,100000,50,true,false,100,100,1,100.0,100,100,646900,100,646800,3880900,0.05999227083011285,1.0,2.1624419689178467,-2236.8855823826443,-2236.8855823790605,3.6e-14,8.8e-14,1.2e-13,1.2e-14,1.6e-15,3.0e+01,727975,1268,8.0e-13,3.5e-11,4.9e-12,2.7e-13,344.9,1.1e-06
w09_120.0,maxCutMet,821219842,tol,383.3,750,1.0e-12,7200.0,100000,50,true,false,120,120,1,120.0,120,120,1123480,120,1123360,6740280,0.04166295795207747,1.8,4.188122034072876,-2920.447975206845,-2920.447975209414,1.7e-14,6.9e-13,3.4e-13,2.6e-14,2.9e-15,3.3e+01,834022,1486,4.4e-13,3.6e-10,1.7e-11,7.0e-13,697.8,1.2e-06
w09_140.0,maxCutMet,661364480,tol,675.8,800,1.0e-12,7200.0,100000,50,true,false,140,140,1,140.0,140,140,1790460,140,1790320,10742060,0.03061025019907095,4.7,8.062062978744507,-3837.3260323620616,-3837.3260323588524,1.4e-14,6.2e-13,5.7e-13,3.2e-14,5.0e-17,3.5e+01,1054295,1880,4.2e-13,4.2e-10,3.3e-11,9.4e-13,1466.7,1.1e-06
w09_160.0,maxCutMet,981172035,tol,945.3,750,1.0e-12,7200.0,100000,50,true,false,160,160,1,160.0,160,160,2679840,160,2679680,16078240,0.02343633388560511,6.9,12.381854057312012,-4652.989255133014,-4652.98925507755,1.8e-13,4.1e-13,4.2e-13,3.8e-14,3.9e-15,3.4e+01,1110407,2088,6.0e-12,3.4e-10,2.8e-11,1.3e-12,2137.3,8.1e-07
w09_180.0,maxCutMet,67995657,tol,1670.6,1150,1.0e-12,7200.0,100000,50,true,false,180,180,1,180.0,180,180,3823620,180,3823440,22940820,0.018517792040003974,9.3,17.85946798324585,-5608.910623496189,-5608.910623433593,1.5e-13,7.5e-13,2.1e-13,1.5e-14,3.8e-15,1.2e+01,1453602,2315,5.6e-12,7.4e-10,1.6e-11,5.8e-13,3687.1,7.7e-07
w09_200.0,maxCutMet,6347352,tol,3768.5,1700,1.0e-12,7200.0,100000,50,true,false,200,200,1,200.0,200,200,5253800,200,5253600,31521800,0.014999524153945716,13.9,29.219916820526123,-6639.404513889022,-6639.404513842797,8.4e-14,4.3e-13,1.5e-13,1.2e-14,3.2e-15,9.6e+00,2672696,2592,3.5e-12,5.0e-10,1.2e-11,4.9e-13,6243.1,9.2e-07
w09_250.0,maxCutMet,479307064,tol,4537.2,800,1.0e-12,7200.0,100000,50,true,false,250,250,1,250.0,250,250,10292250,250,10292000,61752250,0.009599805679030338,29.7,65.94290900230408,-9542.865814855077,-9542.865814605902,2.6e-13,1.0e-12,3.3e-13,3.4e-14,1.3e-15,3.6e+01,1736186,3030,1.3e-11,1.3e-09,3.4e-11,1.7e-12,9697.3,1.6e-04
Geom. mean,maxCutMet,,,1055.8,,,,,,,,,,,,,,,,,,,,,,,,,,,,,,,,,,,877.6,
\end{filecontents*}
\pgfplotstabletypeset[
    col sep=comma,
    string type,
    columns={instance,n,equations,inequalities,statusCode,time,iter,gapOriginal,pinfOriginal,dinfOriginal,complOriginal,timeSCS,maxErrorSCS},
    columns/instance/.style={column name=\textbf{Instance}, column type=l},
    columns/n/.style={column name=$n$, column type=r},
    columns/equations/.style={column name=$m_a$, column type=r},
    columns/inequalities/.style={column name=$m_b$, column type=r},
    columns/statusCode/.style={column name=\textbf{Status}, column type=|c},
    columns/time/.style={column name=\textbf{Time~[s]}, column type=r},
    columns/iter/.style={column name=\textbf{Iter}, column type=r},
    columns/gapOriginal/.style={column name=\texttt{gap}, column type=c},
    columns/pinfOriginal/.style={column name=\texttt{pinf}, column type=c},
    columns/dinfOriginal/.style={column name=\texttt{dinf}, column type=c},
    columns/complOriginal/.style={column name=\texttt{compl}, column type=c|},
    columns/timeSCS/.style={column name=\textbf{Time~[s]}, column type=r},
    columns/maxErrorSCS/.style={column name=\textbf{Max. error}, column type=c},
    every head row/.style={
        before row={
            \toprule
            & & & & \multicolumn{7}{c|}{\textbf{Augmented Mixing Method}} & \multicolumn{2}{c}{\textbf{SCS}} \\
        },
        after row=\midrule
    },
    every last row/.style={before row=\midrule, after row=\bottomrule}
]{maxCutMet.csv}
}
\caption{Results of the Augmented Mixing Method with \texttt{tol=1e-12} and \texttt{SCS}~\cite{brendan2016conic} on various instances of the Max-Cut SDP relaxation~\eqref{eq:max_cut_sdp}, strengthened by all triangle inequalities~\eqref{eq:triangle_inequalities}. The last row reports the geometric mean of runtimes over all instances.}
\label{table:max_cut_met}
\end{table}

\subsection{Bounding the stability number of a graph}

Given a simple undirected graph $G=(V,E)$ with $\lvert V \rvert = n$, the \emph{stability number} of $G$ can be computed as the optimal value of
\begin{equation*}
\begin{array}{rl}
\text{maximize} \quad & e_n^\top x \\[1.2ex]
\text{subject to} \quad & x_i x_j = 0, \quad \forall \{ i,j \} \in  E, \\
& x \in \{0,1\}^n. \\
\end{array}
\end{equation*}
A well-known SDP relaxation for obtaining an upper bound on the stability number is given by
\begin{equation}\label{eq:theta}
\begin{array}{rl}
\text{maximize} \quad & \langle J_n , X \rangle \\[1.2ex]
\text{subject to} \quad & X_{ij} = 0, \quad \forall \{ i,j \} \in  E, \\
& \operatorname{trace}(X) = 1, \\
& X \in \mathcal{S}_+^n. \\
\end{array}
\end{equation}
The optimal value of this SDP relaxation is also known as the Lovász number (or Lovász theta function) of $G$. This SDP relaxation can be strengthened by adding nonnegativity constraints on entries of $X$, resulting in the doubly nonnegative (DNN) relaxation
\begin{equation}\label{eq:theta_prime}
\begin{array}{rl}
\text{maximize} \quad & \langle J_n , X \rangle \\[1.2ex]
\text{subject to} \quad & X_{ij} = 0, \quad \forall \{ i,j \} \in  E, \\
& X_{ij} \geq 0, \quad \forall \{ i,j \} \notin  E, \\
& \operatorname{trace}(X) = 1, \\
& X \in \mathcal{S}_+^n. \\
\end{array}
\end{equation}
Tables~\ref{table:theta} and~\ref{table:thetaPrime} show numerical results of the Augmented Mixing Method applied to the basic SDP relaxation~\eqref{eq:theta} and the strengthened DNN relaxation~\eqref{eq:theta_prime} on selected instances from the second DIMACS implementation challenge~\cite{DIMACS1992}. For the DNN relaxation, we exclude instances with $n$ larger than $300$ because none of the general-purpose solvers in our benchmark set could solve these larger instances in reasonable time. We note that other solvers, such as \texttt{RiNNAL}~\cite{hou2025lowrank}, \texttt{SDPNAL+}~\cite{yang2015majorized}, and \texttt{HALLaR}~\cite{monteiro2024low}, can handle graphs with substantially larger numbers of vertices. In contrast to the Augmented Mixing Method, which treats the nonnegativity constraints as ordinary affine inequality constraints, these methods are primarily designed for scalability on this particular problem class. The purpose of our experiments, however, is to assess high-accuracy primal-dual solution quality in a general-purpose setting.

With a few exceptions in Table~\ref{table:theta}, the Augmented Mixing Method is able to compute high-accuracy primal-dual solutions of the basic SDP relaxation~\eqref{eq:theta} in reasonable time. The primal-dual errors are almost always smaller than $10^{-10}$, and in some cases significantly smaller. Five instances are not solved by the method: the time limit is reached on three, and the maximum number of outer iterations is reached on two.

For the DNN relaxation~\eqref{eq:theta_prime}, shown in Table~\ref{table:thetaPrime}, the primal-dual errors are again typically below $10^{-10}$ and occasionally even below $10^{-12}$. As with the basic SDP relaxation, five instances are not solved by the Augmented Mixing Method, although on one of them the maximum primal-dual error is still smaller than $10^{-9}$. Interestingly, the solving times for the DNN relaxation~\eqref{eq:theta_prime} are not much higher than for the basic SDP relaxation~\eqref{eq:theta}, despite the presence of significantly more affine constraints. In some cases, the DNN relaxation is even solved faster and with higher accuracy than its basic counterpart.

\begin{figure}[ht]
    \centering
    \includegraphics[width=\textwidth]{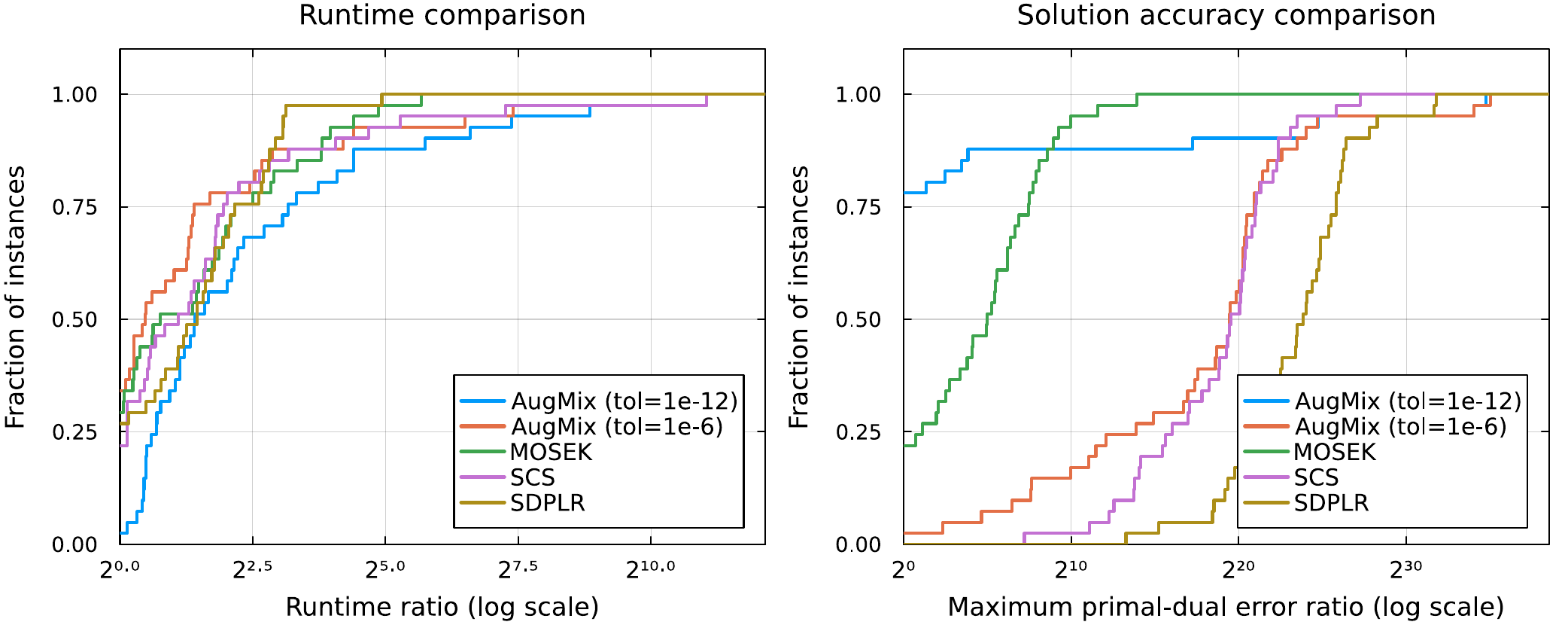}
    \caption{Performance profiles comparing the Augmented Mixing Method (with \texttt{tol=1e-12} and \texttt{tol=1e-6}) against \texttt{MOSEK}, \texttt{SCS}, and \texttt{SDPLR} on the basic SDP relaxation~\eqref{eq:theta} for bounding the stability number of a graph. Left: Runtime comparison. Right: Maximum primal-dual error comparison.}
    \label{fig:perfplot_theta}
\end{figure}

The left plot in Figure~\ref{fig:perfplot_theta}, corresponding to the basic SDP relaxation~\eqref{eq:theta}, shows that all solvers exhibit similar runtimes on the majority of instances. However, the version of the Augmented Mixing Method with \texttt{tol=1e-6} is the fastest solver on about 75\% of all instances, while the version with \texttt{tol=1e-12} is consistently the slowest. Nonetheless, the right plot in Figure~\ref{fig:perfplot_theta} reveals that the Augmented Mixing Method with the default setting \texttt{tol=1e-12} produces the most accurate solutions on over 75\% of all instances. Only \texttt{MOSEK} achieves comparable accuracy and successfully solves all instances considered. In contrast, the accuracy of the Augmented Mixing Method with \texttt{tol=1e-6} is similar to \texttt{SCS}, and \texttt{SDPLR} clearly yields the least accurate solutions overall.

\begin{figure}[ht]
    \centering
    \includegraphics[width=\textwidth]{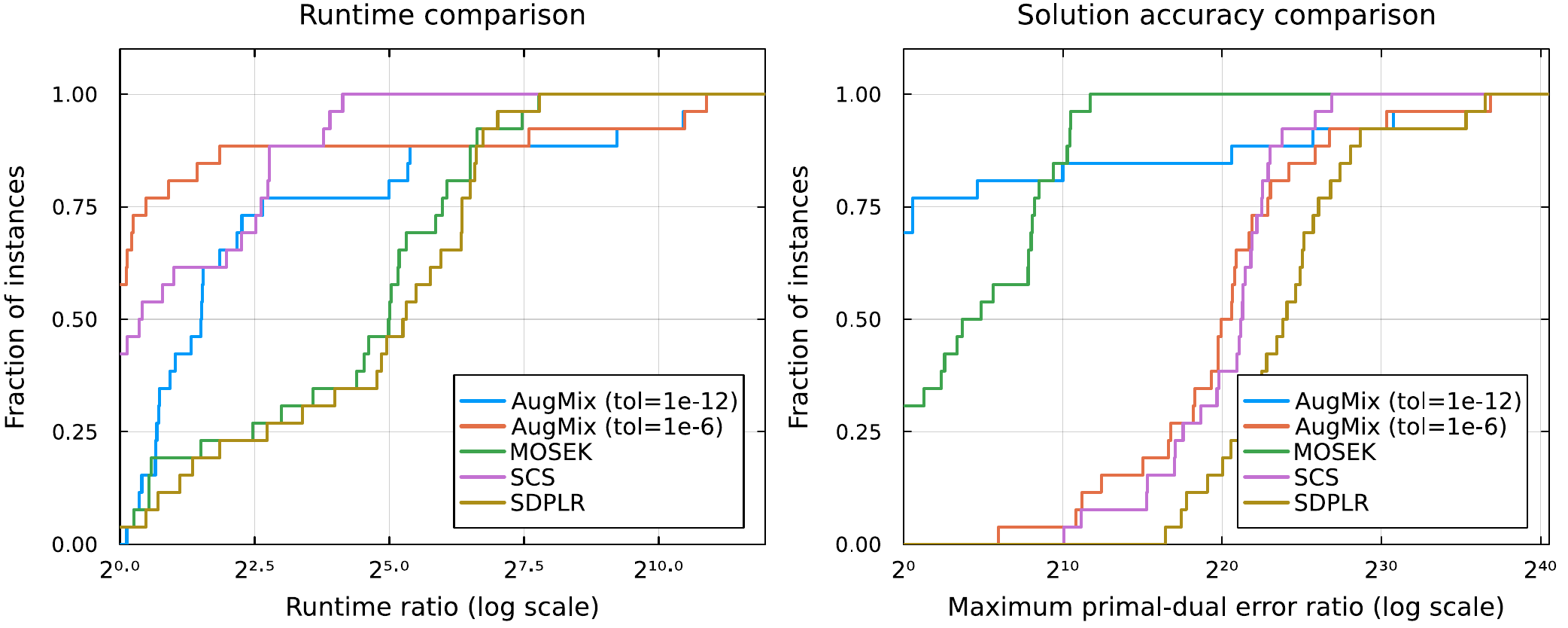}
    \caption{Performance profiles comparing the Augmented Mixing Method (with \texttt{tol=1e-12} and \texttt{tol=1e-6}) against \texttt{MOSEK}, \texttt{SCS}, and \texttt{SDPLR} on the DNN relaxation~\eqref{eq:theta_prime} for bounding the stability number of a graph. Left: Runtime comparison. Right: Maximum primal-dual error comparison.}
    \label{fig:perfplot_thetaPrime}
\end{figure}

Regarding the DNN relaxation~\eqref{eq:theta_prime}, the solution accuracy comparison in Figure~\ref{fig:perfplot_thetaPrime} is nearly identical to that of the basic SDP relaxation~\eqref{eq:theta} in Figure~\ref{fig:perfplot_theta}. The runtime plot, however, reveals significant differences. \texttt{MOSEK} and \texttt{SDPLR} are by far the slowest solvers, often by orders of magnitude. On almost all instances, the Augmented Mixing Method with \texttt{tol=1e-6} is the fastest solver, only occasionally outperformed by \texttt{SCS}. Notably, the Augmented Mixing Method with \texttt{tol=1e-12} is substantially faster than \texttt{MOSEK} on the instances it solves, while producing solutions of similar or even higher accuracy. This highlights the method’s favorable trade-off between runtime and accuracy, as demanding higher precision leads to only a modest increase in runtime.

\begin{table}[ht]
\centering
\resizebox{1.0\textwidth}{!}{
\setlength{\tabcolsep}{7pt}
\renewcommand{\arraystretch}{1.1}
\begin{filecontents*}{theta.csv}
instance,type,seed,statusCode,time,iter,tol,timeLimit,maxIters,itersZ,scaling,warmStart,k,n,numBlocks,avgBlockSize,minBlockSize,maxBlockSize,m,equations,inequalities,nnzA,densityAPercent,timeScaling,timeData,primalObj,dualObj,gap,pinf,dinf,compl,complNotPsdAbs,mu,evals,nnzY,gapOriginal,pinfOriginal,dinfOriginal,complOriginal
brock200_1,theta,103797663,tol,28.7,900,1.0e-12,7200.0,100000,50,true,false,101,200,1,200.0,200,200,5067,5067,0,10332,0.005097690941385435,0.0,0.013742208480834961,-27.456641221425084,-27.45664122141382,1.6e-13,2.4e-13,4.6e-13,3.6e-13,1.6e-14,2.0e+00,3590134,5067,2.0e-13,1.5e-14,4.7e-11,4.4e-13
brock200_2,theta,877428585,tol,13.8,300,1.0e-12,7200.0,100000,50,true,false,142,200,1,200.0,200,200,10025,10025,0,20248,0.005049376558603491,0.0,0.028522014617919922,-14.227205727552716,-14.22720572755061,4.9e-14,4.7e-13,1.2e-13,3.8e-13,4.3e-15,7.1e-01,1387085,10025,7.1e-14,3.3e-14,1.2e-11,5.5e-13
brock200_3,theta,569821980,tol,27.8,1250,1.0e-12,7200.0,100000,50,true,false,126,200,1,200.0,200,200,7853,7853,0,15904,0.005063033235706099,0.0,0.019808053970336914,-18.820535919235876,-18.820535919225087,2.1e-13,6.0e-13,1.5e-13,8.2e-13,5.1e-15,6.3e-01,2575018,7853,2.8e-13,5.0e-14,1.5e-11,1.1e-12
brock200_4,theta,444422726,tol,21.0,500,1.0e-12,7200.0,100000,50,true,false,117,200,1,200.0,200,200,6812,6812,0,13822,0.005072665883734586,0.0,0.016002178192138672,-21.29347571404538,-21.2934757140457,5.6e-15,7.2e-14,1.3e-13,6.6e-14,2.7e-17,9.6e-01,2437040,6812,8.2e-15,3.6e-15,1.3e-11,8.6e-14
brock400_1,theta,894079699,tol,173.6,1200,1.0e-12,7200.0,100000,50,true,false,201,400,1,400.0,400,400,20078,20078,0,40554,0.001262389182189461,0.0,0.1127769947052002,-39.701898706154964,-39.70189870615602,1.1e-14,1.7e-14,5.4e-13,2.8e-14,5.5e-16,1.5e+00,8441490,20078,1.3e-14,1.6e-15,1.1e-10,3.4e-14
brock400_2,theta,716028077,tol,378.0,4200,1.0e-12,7200.0,100000,50,true,false,201,400,1,400.0,400,400,20015,20015,0,40428,0.0012624281788658506,0.0,0.1505110263824463,-39.560523575971104,-39.560523575905044,6.7e-13,6.9e-13,2.0e-13,9.5e-13,1.8e-14,6.3e-01,15982776,20015,8.2e-13,2.4e-14,4.1e-11,1.2e-12
brock400_3,theta,871352809,tol,124.2,900,1.0e-12,7200.0,100000,50,true,false,201,400,1,400.0,400,400,20120,20120,0,40638,0.0012623633200795229,0.1,0.17025995254516602,-39.48055249664963,-39.48055249662987,2.0e-13,5.0e-13,1.9e-13,7.7e-13,1.5e-14,6.5e-01,5329978,20120,2.5e-13,4.5e-14,3.7e-11,9.6e-13
brock400_4,theta,135823888,tol,135.5,750,1.0e-12,7200.0,100000,50,true,false,201,400,1,400.0,400,400,20036,20036,0,40470,0.0012624151527250947,0.1,0.1672530174255371,-39.59962551127508,-39.59962551124717,2.8e-13,4.6e-13,2.0e-13,7.5e-13,2.1e-14,1.8e+00,6151887,20036,3.5e-13,1.6e-14,4.0e-11,9.3e-13
c_fat200_1,theta,497340725,tol,43.6,950,1.0e-12,7200.0,100000,50,true,false,192,200,1,200.0,200,200,18367,18367,0,36932,0.005026950509065172,0.1,0.04589700698852539,-12.000000000010356,-11.999999999999956,2.7e-13,6.6e-13,2.5e-15,6.2e-14,1.2e-15,6.0e-01,3473372,18367,4.2e-13,3.3e-14,5.5e-13,9.5e-14
c_fat200_2,theta,690319336,tol,5.0,150,1.0e-12,7200.0,100000,50,true,false,183,200,1,200.0,200,200,16666,16666,0,33530,0.005029701188047522,0.0,0.038815975189208984,-23.99999999997673,-23.999999999999638,3.7e-13,7.4e-13,7.4e-17,6.1e-16,5.3e-15,9.6e-01,286043,16666,4.7e-13,3.7e-14,1.4e-12,3.7e-15
c_fat200_5,theta,676266835,tol,5.9,200,1.0e-12,7200.0,100000,50,true,false,152,200,1,200.0,200,200,11428,11428,0,23054,0.005043314665733286,0.0,0.026495933532714844,-60.345275298509584,-60.345275298509485,7.5e-16,5.6e-16,1.8e-17,4.9e-16,1.6e-16,6.9e-01,530705,11428,5.8e-16,1.7e-15,1.5e-12,1.3e-15
hamming6_2,theta,885829393,tol,0.1,100,1.0e-12,7200.0,100000,50,true,false,20,64,1,64.0,64,64,193,193,0,448,0.05667098445595855,0.0,0.00023818016052246094,-32.00000000000004,-31.99999999999999,7.4e-16,3.4e-15,4.7e-17,1.7e-17,1.8e-16,3.1e+00,41275,193,6.0e-16,7.8e-16,3.3e-13,3.1e-16
hamming6_4,theta,248854592,tol,0.5,150,1.0e-12,7200.0,100000,50,true,false,52,64,1,64.0,64,64,1313,1313,0,2688,0.04998095963442498,0.0,0.0008819103240966797,-5.33333333333322,-5.333333333333362,7.6e-15,1.2e-14,2.9e-15,1.4e-14,1.5e-15,5.8e-01,167566,1313,1.1e-14,2.9e-15,8.8e-14,2.2e-14
hamming8_2,theta,808980131,tol,1.9,150,1.0e-12,7200.0,100000,50,true,false,46,256,1,256.0,256,256,1025,1025,0,2304,0.003429878048780488,0.0,0.0029380321502685547,-127.9999999999995,-128.00000000000006,2.0e-15,2.6e-14,0.0e+00,2.3e-16,2.2e-16,5.5e+00,224454,1025,2.4e-15,1.8e-15,3.5e-12,5.3e-16
hamming8_4,theta,547964795,tol,39.7,550,1.0e-12,7200.0,100000,50,true,false,154,256,1,256.0,256,256,11777,11777,0,23808,0.003084667147830517,0.0,0.053542137145996094,-16.000000000001222,-16.00000000000043,1.7e-14,4.1e-13,2.6e-13,8.5e-13,9.0e-15,4.2e+00,3250053,11777,2.4e-14,1.8e-14,3.4e-11,1.2e-12
hamming10_2,theta,744123520,tol,84.2,400,1.0e-12,7200.0,100000,50,true,false,102,1024,1,1024.0,1024,1024,5121,5121,0,11264,0.0002097673794180824,0.0,0.08901095390319824,-511.9999999999852,-512.0000000000116,2.5e-14,1.9e-13,0.0e+00,1.1e-14,1.1e-14,2.1e+01,2102547,5121,2.6e-14,9.7e-15,8.7e-12,1.6e-14
johnson8_2_4,theta,844880713,tol,0.7,900,1.0e-12,7200.0,100000,50,true,false,19,28,1,28.0,28,28,169,169,0,364,0.27472527472527475,0.0,8.702278137207031e-5,-3.9999999999999387,-3.999999999999724,1.6e-14,3.2e-13,4.9e-13,1.0e-13,2.1e-14,8.7e+00,413135,169,2.4e-14,4.3e-14,7.1e-12,1.5e-13
johnson8_4_4,theta,575588055,tol,0.4,150,1.0e-12,7200.0,100000,50,true,false,34,70,1,70.0,70,70,561,561,0,1190,0.04329004329004329,0.0,0.00041794776916503906,-13.999999999999718,-14.000000000000192,1.3e-14,1.9e-13,2.2e-14,6.0e-14,5.2e-15,3.7e+00,157454,561,1.7e-14,1.6e-14,7.8e-13,7.6e-14
johnson16_2_4,theta,305938573,tol,12.3,1300,1.0e-12,7200.0,100000,50,true,false,58,120,1,120.0,120,120,1681,1681,0,3480,0.014376363275827881,0.0,0.0017299652099609375,-7.999999999999743,-8.000000000000929,4.4e-14,4.8e-13,5.9e-13,3.7e-13,3.4e-14,9.7e+00,2807145,1681,6.9e-14,3.1e-14,3.6e-11,5.9e-13
johnson32_2_4,theta,307346117,time,7200.1,31104,1.0e-12,7200.0,100000,50,true,false,173,496,1,496.0,496,496,14881,14881,0,30256,0.0008264489682665273,0.1,0.15181708335876465,-17.352100303387203,-15.581842253621407,3.2e-02,8.2e-03,Inf,Inf,1.5e-02,1.2e+00,349848977,14881,5.2e-02,8.2e-03,2.0e+00,2.8e-02
keller4,theta,2836000,tol,91.8,4250,1.0e-12,7200.0,100000,50,true,false,102,171,1,171.0,171,171,5101,5101,0,10371,0.006953013777486302,0.0,0.007573127746582031,-14.01224167801701,-14.012241678016863,3.5e-15,5.1e-13,2.8e-13,9.6e-13,1.2e-14,2.1e+00,13107585,5101,5.2e-15,2.8e-14,2.4e-11,1.4e-12
MANN_a9,theta,864038,tol,0.4,550,1.0e-12,7200.0,100000,50,true,false,13,45,1,45.0,45,45,73,73,0,189,0.1278538812785388,0.0,9.989738464355469e-5,-17.475031613081008,-17.475031613081523,1.2e-14,1.8e-13,2.2e-13,5.2e-13,2.6e-14,3.6e+00,253044,73,1.6e-14,1.9e-14,5.1e-12,6.1e-13
MANN_a27,theta,428379973,tol,122.9,5600,1.0e-12,7200.0,100000,50,true,false,38,378,1,378.0,378,378,703,703,0,1782,0.0017740619244378643,0.0,0.006869792938232422,-132.76289108572172,-132.7628910854475,9.6e-13,7.2e-14,2.4e-13,7.2e-13,9.4e-13,3.6e+00,17621680,703,1.0e-12,2.1e-14,6.0e-11,7.7e-13
p_hat300_1,theta,655809194,tol,3190.7,28550,1.0e-12,7200.0,100000,50,true,false,261,300,1,300.0,300,300,33918,33918,0,68134,0.0022319843282164174,0.1,0.16380906105041504,-10.067964586910115,-10.067964586912453,6.2e-14,7.1e-13,1.9e-13,9.7e-13,2.5e-14,6.9e-01,146196900,33918,1.1e-13,3.8e-14,2.9e-11,1.7e-12
p_hat300_2,theta,454656574,time,7200.2,52742,1.0e-12,7200.0,100000,50,true,false,215,300,1,300.0,300,300,22923,22923,0,46144,0.0022366667151381193,0.0,0.06670713424682617,-26.96603479042605,-26.966034812820787,3.1e-10,4.3e-08,Inf,Inf,8.0e-10,8.3e+00,421732631,22923,4.1e-10,1.8e-09,1.1e-03,2.0e-08
p_hat300_3,theta,982868641,tol,259.0,3300,1.0e-12,7200.0,100000,50,true,false,152,300,1,300.0,300,300,11461,11461,0,23220,0.0022511124683709974,0.0,0.05006217956542969,-41.169930018811066,-41.16993001881169,6.3e-15,9.3e-13,3.2e-13,2.0e-13,2.7e-15,1.0e+00,21745258,11461,7.5e-15,3.8e-14,4.8e-11,2.4e-13
p_hat500_3,theta,633978584,tol,741.9,2100,1.0e-12,7200.0,100000,50,true,false,249,500,1,500.0,500,500,30951,30951,0,62400,0.0008064359794513909,0.1,0.4075901508331299,-58.567905885527075,-58.56790588553061,2.5e-14,8.1e-13,5.4e-13,1.8e-13,2.7e-15,7.1e+00,22540743,30951,3.0e-14,2.6e-14,1.3e-10,2.2e-13
san200_0_7_1,theta,249294285,tol,2.0,100,1.0e-12,7200.0,100000,50,true,false,110,200,1,200.0,200,200,5971,5971,0,12140,0.005082900686652152,0.0,0.013659954071044922,-30.000000000003862,-29.999999999999716,5.6e-14,1.8e-13,4.8e-17,1.3e-17,3.9e-15,7.6e-01,179309,5971,6.8e-14,3.3e-14,3.9e-13,5.4e-16
san200_0_7_2,theta,772154870,iter,2223.3,100000,1.0e-12,7200.0,100000,50,true,false,110,200,1,200.0,200,200,5971,5971,0,12140,0.005082900686652152,0.0,0.012682914733886719,-17.993417517887927,-18.00280916813754,1.9e-04,1.1e-03,Inf,Inf,4.9e-05,2.1e+00,240257510,5971,2.5e-04,5.3e-05,3.0e-01,9.8e-04
san200_0_9_1,theta,13572006,tol,1.4,150,1.0e-12,7200.0,100000,50,true,false,64,200,1,200.0,200,200,1991,1991,0,4180,0.005248618784530387,0.0,0.003031015396118164,-70.0000000000002,-70.0,1.3e-15,2.4e-15,3.4e-19,8.3e-17,1.2e-16,6.0e-01,166865,1991,1.8e-15,1.7e-15,3.2e-12,8.4e-16
san200_0_9_2,theta,610464760,tol,1.0,100,1.0e-12,7200.0,100000,50,true,false,64,200,1,200.0,200,200,1991,1991,0,4180,0.005248618784530387,0.0,0.003068208694458008,-59.999999999992674,-59.99999999999992,5.4e-14,1.7e-13,6.5e-18,8.3e-17,6.8e-16,7.6e-01,133701,1991,6.0e-14,8.5e-15,1.1e-12,1.7e-15
san200_0_9_3,theta,876189547,iter,1105.8,100000,1.0e-12,7200.0,100000,50,true,false,64,200,1,200.0,200,200,1991,1991,0,4180,0.005248618784530387,0.0,0.003696918487548828,-43.99498952363475,-43.99938669726905,4.3e-05,1.7e-04,Inf,Inf,5.9e-06,4.6e+00,176415049,1991,4.9e-05,8.4e-06,1.9e-01,2.7e-05
san400_0_5_1,theta,416452836,tol,71.2,300,1.0e-12,7200.0,100000,50,true,false,283,400,1,400.0,400,400,39901,39901,0,80200,0.0012562341795944963,0.1,0.34542298316955566,-13.000000000000927,-13.00000000000012,1.8e-14,1.7e-13,0.0e+00,3.1e-15,3.1e-15,8.2e-01,2382160,39901,2.7e-14,1.8e-14,9.4e-14,6.3e-15
san400_0_7_1,theta,84999349,tol,14.0,100,1.0e-12,7200.0,100000,50,true,false,219,400,1,400.0,400,400,23941,23941,0,48280,0.0012603901257257424,0.0,0.20465707778930664,-40.00000000000904,-40.000000000000604,8.4e-14,4.1e-13,0.0e+00,6.5e-15,6.4e-15,1.1e+00,570325,23941,1.0e-13,1.6e-14,8.2e-13,8.8e-15
san400_0_7_2,theta,877017282,tol,32.4,250,1.0e-12,7200.0,100000,50,true,false,219,400,1,400.0,400,400,23941,23941,0,48280,0.0012603901257257424,0.0,0.24290990829467773,-30.00000000002274,-29.9999999999972,3.2e-13,3.8e-13,2.3e-16,2.4e-15,3.3e-14,8.5e-01,1366156,23941,4.2e-13,1.1e-13,3.7e-13,8.2e-15
san400_0_7_3,theta,235377152,time,7200.1,51315,1.0e-12,7200.0,100000,50,true,false,219,400,1,400.0,400,400,23941,23941,0,48280,0.0012603901257257424,0.0,0.2506868839263916,-22.000000008424124,-21.99999998882462,3.1e-10,1.7e-08,Inf,Inf,1.7e-10,1.6e+00,310003813,23941,4.4e-10,6.0e-10,8.1e-05,1.5e-08
san400_0_9_1,theta,827345035,tol,5.3,100,1.0e-12,7200.0,100000,50,true,false,127,400,1,400.0,400,400,7981,7981,0,16360,0.0012811677734619721,0.0,0.037034034729003906,-99.99999999999957,-100.0,1.9e-15,5.0e-15,6.3e-19,1.3e-16,9.7e-17,1.1e+00,294626,7981,1.8e-15,3.3e-16,1.5e-12,5.5e-16
sanr200_0_7,theta,506670270,tol,46.9,1500,1.0e-12,7200.0,100000,50,true,false,110,200,1,200.0,200,200,6033,6033,0,12264,0.005082048731974142,0.0,0.010937929153442383,-23.836157572458376,-23.836157572459012,1.0e-14,2.2e-14,7.2e-13,2.9e-14,4.4e-16,1.5e+00,5832692,6033,1.3e-14,2.3e-15,7.3e-11,3.7e-14
sanr200_0_9,theta,336841823,tol,12.8,1300,1.0e-12,7200.0,100000,50,true,false,64,200,1,200.0,200,200,2038,2038,0,4274,0.0052428851815505394,0.0,0.0030510425567626953,-49.27351787933943,-49.27351787937168,2.9e-13,8.0e-13,7.1e-14,2.4e-13,1.6e-14,6.0e-01,2110823,2038,3.2e-13,1.6e-13,7.1e-12,2.7e-13
sanr400_0_5,theta,62849739,tol,501.3,2550,1.0e-12,7200.0,100000,50,true,false,283,400,1,400.0,400,400,39817,39817,0,80032,0.0012562473315418038,0.1,0.3933908939361572,-20.317656713148757,-20.317656713149155,6.6e-15,1.7e-14,8.7e-13,5.6e-14,3.2e-16,8.5e-01,16267178,39817,1.0e-14,1.5e-15,1.7e-10,8.2e-14
sanr400_0_7,theta,6514372,tol,148.9,950,1.0e-12,7200.0,100000,50,true,false,219,400,1,400.0,400,400,23932,23932,0,48262,0.0012603940330937658,0.1,0.07798910140991211,-34.278346079067305,-34.278346079069614,2.6e-14,2.5e-14,5.7e-13,5.2e-14,1.9e-16,9.8e-01,6517809,23932,3.3e-14,8.8e-16,1.1e-10,6.6e-14
\end{filecontents*}
\pgfplotstabletypeset[
    col sep=comma,
    string type,
    columns={instance,n,equations,inequalities,statusCode,time,iter,gapOriginal,pinfOriginal,dinfOriginal,complOriginal},
    columns/instance/.style={column name=\textbf{Instance}, column type=l},
    columns/n/.style={column name=$n$, column type=r},
    columns/equations/.style={column name=$m_a$, column type=r},
    columns/inequalities/.style={column name=$m_b$, column type=r},
    columns/statusCode/.style={column name=\textbf{Status}, column type=c},
    columns/time/.style={column name=\textbf{Time~[s]}, column type=r},
    columns/iter/.style={column name=\textbf{Iter}, column type=r},
    columns/gapOriginal/.style={column name=\texttt{gap}, column type=c},
    columns/pinfOriginal/.style={column name=\texttt{pinf}, column type=c},
    columns/dinfOriginal/.style={column name=\texttt{dinf}, column type=c},
    columns/complOriginal/.style={column name=\texttt{compl}, column type=c},
    every head row/.style={before row=\toprule, after row=\midrule},
    every last row/.style={after row=\bottomrule}
]{theta.csv}
}
\caption{Results of the Augmented Mixing Method with \texttt{tol=1e-12} on selected instances of the basic SDP relaxation~\eqref{eq:theta} for computing upper bounds on the stability number of a graph.}
\label{table:theta}
\end{table}

\begin{table}[ht]
\centering
\resizebox{1.0\textwidth}{!}{
\setlength{\tabcolsep}{7pt}
\renewcommand{\arraystretch}{1.1}
\begin{filecontents*}{thetaPrime.csv}
instance,type,seed,statusCode,time,iter,tol,timeLimit,maxIters,itersZ,scaling,warmStart,k,n,numBlocks,avgBlockSize,minBlockSize,maxBlockSize,m,equations,inequalities,nnzA,densityAPercent,timeScaling,timeData,primalObj,dualObj,gap,pinf,dinf,compl,complNotPsdAbs,mu,evals,nnzY,gapOriginal,pinfOriginal,dinfOriginal,complOriginal
brock200_1,thetaPrime,452490270,tol,26.5,750,1.0e-12,7200.0,100000,50,true,false,200,200,1,200.0,200,200,19901,5067,14834,40000,0.00502487312195367,0.0,0.07382583618164062,-27.196716094132412,-27.19671609409841,5.0e-13,8.7e-13,2.8e-13,9.8e-13,7.0e-16,6.3e-01,1868602,5907,6.1e-13,4.4e-14,2.8e-11,1.2e-12
brock200_2,thetaPrime,917291261,tol,18.6,300,1.0e-12,7200.0,100000,50,true,false,200,200,1,200.0,200,200,19901,10025,9876,40000,0.00502487312195367,0.0,0.047833919525146484,-14.131007419881275,-14.131007419874576,1.6e-13,5.3e-13,1.6e-13,6.1e-13,1.3e-14,8.0e-01,1505609,10746,2.3e-13,3.2e-14,1.6e-11,8.9e-13
brock200_3,thetaPrime,296298601,tol,16.9,300,1.0e-12,7200.0,100000,50,true,false,200,200,1,200.0,200,200,19901,7853,12048,40000,0.00502487312195367,0.0,0.04720616340637207,-18.67179791085259,-18.671797910845015,1.5e-13,7.1e-13,9.0e-14,4.9e-13,9.0e-15,7.4e-01,1338627,8640,2.0e-13,3.5e-14,9.0e-12,6.5e-13
brock200_4,thetaPrime,725759525,tol,15.4,250,1.0e-12,7200.0,100000,50,true,false,200,200,1,200.0,200,200,19901,6812,13089,40000,0.00502487312195367,0.0,0.05020904541015625,-21.12107382369795,-21.1210738236967,2.2e-14,4.4e-14,7.0e-14,7.7e-14,1.8e-15,1.6e+00,1235653,7569,2.9e-14,2.2e-15,7.0e-12,1.0e-13
c_fat200_1,thetaPrime,414303423,tol,69.7,1650,1.0e-12,7200.0,100000,50,true,false,200,200,1,200.0,200,200,19901,18367,1534,40000,0.00502487312195367,0.0,0.04942798614501953,-12.000000000005649,-11.99999999999996,1.5e-13,9.5e-13,8.1e-15,2.1e-13,1.0e-15,5.8e-01,5401752,18623,2.3e-13,4.7e-14,8.7e-13,3.3e-13
c_fat200_2,thetaPrime,374659006,tol,5.9,200,1.0e-12,7200.0,100000,50,true,false,200,200,1,200.0,200,200,19901,16666,3235,40000,0.00502487312195367,0.0,0.05363607406616211,-23.999999999996422,-23.999999999999993,5.7e-14,1.5e-13,0.0e+00,8.5e-17,1.4e-16,6.0e-01,343592,18107,7.3e-14,7.3e-15,3.2e-13,8.7e-16
c_fat200_5,thetaPrime,681033121,tol,7.8,200,1.0e-12,7200.0,100000,50,true,false,200,200,1,200.0,200,200,19901,11428,8473,40000,0.00502487312195367,0.0,0.06545495986938477,-60.34527529850957,-60.345275298509485,6.5e-16,1.0e-16,2.4e-17,5.7e-16,1.9e-16,7.1e-01,526375,11428,8.8e-16,1.0e-15,3.4e-13,1.4e-15
hamming6_2,thetaPrime,351276835,tol,0.2,100,1.0e-12,7200.0,100000,50,true,false,64,64,1,64.0,64,64,2017,193,1824,4096,0.0495785820525533,0.0,0.0021791458129882812,-31.999999999999226,-32.00000000000015,1.3e-14,1.3e-13,1.1e-15,3.4e-14,2.1e-15,6.1e-01,44115,1025,1.4e-14,1.2e-14,5.3e-13,3.7e-14
hamming6_4,thetaPrime,565228904,iter,657.0,100000,1.0e-12,7200.0,100000,50,true,false,64,64,1,64.0,64,64,2017,1313,704,4096,0.0495785820525533,0.0,0.0015189647674560547,-3.999999999256203,-4.000000002011574,1.7e-10,1.4e-08,Inf,Inf,1.3e-10,7.5e+00,196421548,1537,3.1e-10,1.3e-09,4.2e-05,2.7e-08
hamming8_2,thetaPrime,396003040,tol,12.4,250,1.0e-12,7200.0,100000,50,true,false,256,256,1,256.0,256,256,32641,1025,31616,65536,0.0030636316289329372,0.1,0.20135092735290527,-128.00000000001216,-128.00000000000006,4.5e-14,4.1e-13,9.7e-16,1.2e-13,2.0e-16,3.8e+00,498466,14223,4.7e-14,1.8e-14,5.7e-12,1.2e-13
hamming8_4,thetaPrime,455487179,time,7200.0,44739,1.0e-12,7200.0,100000,50,true,false,256,256,1,256.0,256,256,32641,11777,20864,65536,0.0030636316289329372,0.0,0.11245298385620117,-16.000027330689406,-15.99999910701726,5.9e-07,4.1e-07,Inf,Inf,1.8e-08,1.5e+01,415521427,24208,8.6e-07,1.8e-08,2.2e-02,9.0e-08
johnson8_2_4,thetaPrime,467996275,tol,0.2,200,1.0e-12,7200.0,100000,50,true,false,28,28,1,28.0,28,28,379,169,210,784,0.2638522427440633,0.0,0.00017690658569335938,-4.000000000000696,-4.000000000003235,1.9e-13,7.8e-13,3.9e-13,4.1e-13,2.4e-13,4.3e+00,100454,176,2.8e-13,1.0e-13,5.7e-12,6.0e-13
johnson8_4_4,thetaPrime,493771539,tol,0.9,250,1.0e-12,7200.0,100000,50,true,false,70,70,1,70.0,70,70,2416,561,1855,4900,0.041390728476821195,0.0,0.0019240379333496094,-14.000000000000483,-13.99999999999991,1.6e-14,5.0e-14,3.4e-15,1.2e-14,2.5e-15,1.1e+00,196430,1119,1.9e-14,6.3e-15,1.5e-13,1.5e-14
johnson16_2_4,thetaPrime,529580799,iter,1831.6,100000,1.0e-12,7200.0,100000,50,true,false,120,120,1,120.0,120,120,7141,1681,5460,14400,0.014003640946646127,0.0,0.00819087028503418,-8.092174592646348,-7.986660630736709,3.9e-03,2.5e-03,Inf,Inf,1.2e-03,1.3e+00,287006886,2246,6.2e-03,1.6e-04,3.7e-01,5.4e-03
keller4,thetaPrime,951297849,tol,983.4,36100,1.0e-12,7200.0,100000,50,true,false,171,171,1,171.0,171,171,14536,5101,9435,29241,0.006879471656576775,0.0,0.0530550479888916,-13.465895608401873,-13.465895608400288,4.0e-14,7.6e-14,1.0e-12,1.1e-13,1.1e-15,7.4e-01,78106792,7735,5.7e-14,4.1e-15,8.6e-11,1.5e-13
MANN_a9,thetaPrime,285557868,tol,0.9,550,1.0e-12,7200.0,100000,50,true,false,45,45,1,45.0,45,45,991,73,918,2025,0.10090817356205853,0.0,0.0005099773406982422,-17.475031613080226,-17.475031613082862,6.3e-14,1.7e-13,1.9e-13,4.2e-13,5.8e-14,3.3e+00,250486,73,7.4e-14,1.8e-14,4.4e-12,4.8e-13
p_hat300_1,thetaPrime,68306418,time,7200.0,62029,1.0e-12,7200.0,100000,50,true,false,300,300,1,300.0,300,300,44851,33918,10933,90000,0.00222960469108827,0.1,0.18824505805969238,-10.02020665807274,-10.020206657931855,3.8e-12,1.8e-11,Inf,Inf,6.5e-13,6.5e-01,299582135,35254,6.7e-12,2.1e-12,9.4e-10,5.5e-11
p_hat300_2,thetaPrime,752743828,time,7200.1,29402,1.0e-12,7200.0,100000,50,true,false,300,300,1,300.0,300,300,44851,22923,21928,90000,0.00222960469108827,0.1,0.19840693473815918,-26.713734228912173,-26.71373395155652,3.9e-09,6.9e-07,Inf,Inf,2.8e-08,7.3e+00,342145383,30487,5.1e-09,2.8e-08,6.0e-03,3.5e-07
p_hat300_3,thetaPrime,186442687,tol,555.3,6100,1.0e-12,7200.0,100000,50,true,false,300,300,1,300.0,300,300,44851,11461,33390,90000,0.00222960469108827,0.0,0.2992699146270752,-40.700271656225624,-40.70027165623156,6.0e-14,9.5e-13,1.5e-13,2.3e-13,8.7e-16,1.1e+00,22266207,14830,7.2e-14,3.9e-14,2.3e-11,2.7e-13
san200_0_7_1,thetaPrime,765013372,tol,5.3,150,1.0e-12,7200.0,100000,50,true,false,200,200,1,200.0,200,200,19901,5971,13930,40000,0.00502487312195367,0.0,0.03995800018310547,-29.999999999999993,-30.0,8.5e-17,1.1e-16,0.0e+00,6.6e-18,6.6e-19,5.6e-01,355353,18109,5.8e-17,3.1e-17,3.7e-13,3.8e-16
san200_0_7_2,thetaPrime,904007051,tol,46.5,1400,1.0e-12,7200.0,100000,50,true,false,200,200,1,200.0,200,200,19901,5971,13930,40000,0.00502487312195367,0.0,0.04275012016296387,-18.00000000000372,-17.999999999998206,1.1e-13,5.0e-13,1.4e-13,7.8e-18,3.6e-14,7.6e-01,3510527,17263,1.5e-13,1.5e-13,1.4e-11,4.4e-16
san200_0_9_1,thetaPrime,590039376,tol,4.3,150,1.0e-12,7200.0,100000,50,true,false,200,200,1,200.0,200,200,19901,1991,17910,40000,0.00502487312195367,0.0,0.04682803153991699,-70.00000000000004,-70.0,2.4e-16,7.9e-16,1.6e-18,7.6e-18,4.9e-17,6.0e-01,243671,15623,0.0e+00,1.0e-15,1.6e-12,1.1e-15
san200_0_9_2,thetaPrime,617034027,tol,3.5,100,1.0e-12,7200.0,100000,50,true,false,200,200,1,200.0,200,200,19901,1991,17910,40000,0.00502487312195367,0.0,0.041285037994384766,-60.00000000000512,-59.99999999999849,4.9e-14,4.0e-14,1.3e-16,5.6e-17,1.1e-14,8.0e-01,212442,14458,5.5e-14,4.0e-14,7.1e-13,7.3e-16
san200_0_9_3,thetaPrime,232076948,tol,42.6,1550,1.0e-12,7200.0,100000,50,true,false,200,200,1,200.0,200,200,19901,1991,17910,40000,0.00502487312195367,0.0,0.0436551570892334,-43.99999999999631,-43.99999999999992,3.5e-14,1.2e-13,9.9e-18,3.0e-16,5.5e-16,5.8e-01,2901086,9419,3.9e-14,7.5e-15,1.7e-12,1.5e-15
sanr200_0_7,thetaPrime,119179408,tol,408.7,17350,1.0e-12,7200.0,100000,50,true,false,200,200,1,200.0,200,200,19901,6033,13868,40000,0.00502487312195367,0.0,0.042222023010253906,-23.63328637868803,-23.633286378667613,3.3e-13,8.9e-13,2.0e-13,9.8e-13,1.9e-14,6.3e-01,27146448,6852,4.2e-13,7.3e-14,2.0e-11,1.2e-12
sanr200_0_9,thetaPrime,146230334,tol,39.1,1650,1.0e-12,7200.0,100000,50,true,false,200,200,1,200.0,200,200,19901,2038,17863,40000,0.00502487312195367,0.0,0.0449068546295166,-48.90455389762063,-48.90455389751111,9.8e-13,8.7e-13,6.5e-14,3.0e-13,2.1e-15,5.8e-01,2519494,2664,1.1e-12,4.3e-14,6.3e-12,3.4e-13
\end{filecontents*}
\pgfplotstabletypeset[
    col sep=comma,
    string type,
    columns={instance,n,equations,inequalities,statusCode,time,iter,gapOriginal,pinfOriginal,dinfOriginal,complOriginal},
    columns/instance/.style={column name=\textbf{Instance}, column type=l},
    columns/n/.style={column name=$n$, column type=r},
    columns/equations/.style={column name=$m_a$, column type=r},
    columns/inequalities/.style={column name=$m_b$, column type=r},
    columns/statusCode/.style={column name=\textbf{Status}, column type=c},
    columns/time/.style={column name=\textbf{Time~[s]}, column type=r},
    columns/iter/.style={column name=\textbf{Iter}, column type=r},
    columns/gapOriginal/.style={column name=\texttt{gap}, column type=c},
    columns/pinfOriginal/.style={column name=\texttt{pinf}, column type=c},
    columns/dinfOriginal/.style={column name=\texttt{dinf}, column type=c},
    columns/complOriginal/.style={column name=\texttt{compl}, column type=c},
    every head row/.style={before row=\toprule, after row=\midrule},
    every last row/.style={after row=\bottomrule}
]{thetaPrime.csv}
}
\caption{Results of the Augmented Mixing Method with \texttt{tol=1e-12} on selected instances of the DNN relaxation~\eqref{eq:theta_prime} for computing upper bounds on the stability number of a graph.}
\label{table:thetaPrime}
\end{table}


\subsection{Single-row facility layout problem}

The single-row facility layout problem (SRFLP) was first introduced in~\cite{simmons1969onedimensional} and asks to arrange $n$ facilities of given lengths on a straight line so that a weighted sum of distances between all pairs of facilities is minimized. It belongs to the class of linear quadratic ordering problems and can be formulated as a binary quadratic problem with both a quadratic objective and quadratic equality constraints, as first shown in~\cite{anjos2005asemidefinite}. From this formulation and a suitable choice of the cost matrix $C \in \mathcal{S}^{\binom{n}{2}}$, the SDP relaxation
\begin{equation}\label{eq:sdp_srflp}
\begin{array}{rl}
\text{minimize} \quad & \langle C , X \rangle \\[1.2ex]
\text{subject to} \quad & X_{ij,jk} - X_{ij,ik} - X_{ik,jk} = -1, \quad 1 \leq i < j < k \leq n, \\
& X_{ij,ij} = 1, \quad 1 \leq i < j \leq n, \\
& X \in \mathcal{S}_+^{\binom{n}{2}}, \\
\end{array}
\end{equation}
was derived in~\cite{anjos2005asemidefinite}. Note that each row (and column) of $X$ corresponds to an ordered pair $(i,j)$ of facilities. Thus, for an SRFLP instance with $n$ facilities, the SDP relaxation~\eqref{eq:sdp_srflp} involves a positive semidefinite matrix variable of order $\binom{n}{2}$ and has $\binom{n}{3} + \binom{n}{2}$ equality constraints. To strengthen the relaxation, it was proposed in~\cite{anjos2008computing} to add the triangle inequalities~\eqref{eq:triangle_inequalities}. Note that the number of triangle inequalities grows as $\mathcal{O}(n^6)$ in the SDP relaxation~\eqref{eq:sdp_srflp} and some of them are already implied by the equality constraints.

For the strengthened SDP relaxation with triangle inequalities, we compare the Augmented Mixing Method with \texttt{tol=1e-12} only to the \texttt{SCS} solver, as the number of inequality constraints exceeds the capacity of the other solvers. Consequently, performance profiles are provided only for the basic SDP relaxation~\eqref{eq:sdp_srflp}, where multiple solvers could be applied.

Table~\ref{table:srflp_basic} shows numerical results of the Augmented Mixing Method on the basic SDP relaxation~\eqref{eq:sdp_srflp}, and Table~\ref{table:srflp_met} reports results on the strengthened relaxation including all nonredundant triangle inequalities. The `O' instances in these tables are taken from~\cite{obata1979quadratic}, the `S' instances from~\cite{sarker1989theamoebic}, the `Y' instances from~\cite{yu2003directional}, the `ste' instances from~\cite{anjos2009provably}, the `N40' instances from~\cite{hungerlander2013semidefinite}, and the remaining `N' instances from~\cite{nugent1968anexperimental}. For the strengthened relaxation, we consider only instances with up to $n = 300$.

As shown in Table~\ref{table:srflp_basic}, the Augmented Mixing Method successfully computes high-accuracy solutions on all instances of the basic SDP relaxation~\eqref{eq:sdp_srflp}, with all primal-dual errors below $10^{-10}$ on nearly all instances. Although the time limit is reached on some of the largest instances, the resulting solutions still exhibit small primal-dual errors (at most $10^{-7}$, typically much smaller). The number of outer iterations ranges from roughly~\num{1000} on smaller instances to around~\num{45000} on larger ones, generally increasing with problem size.

\begin{figure}[ht]
    \centering
    \includegraphics[width=\textwidth]{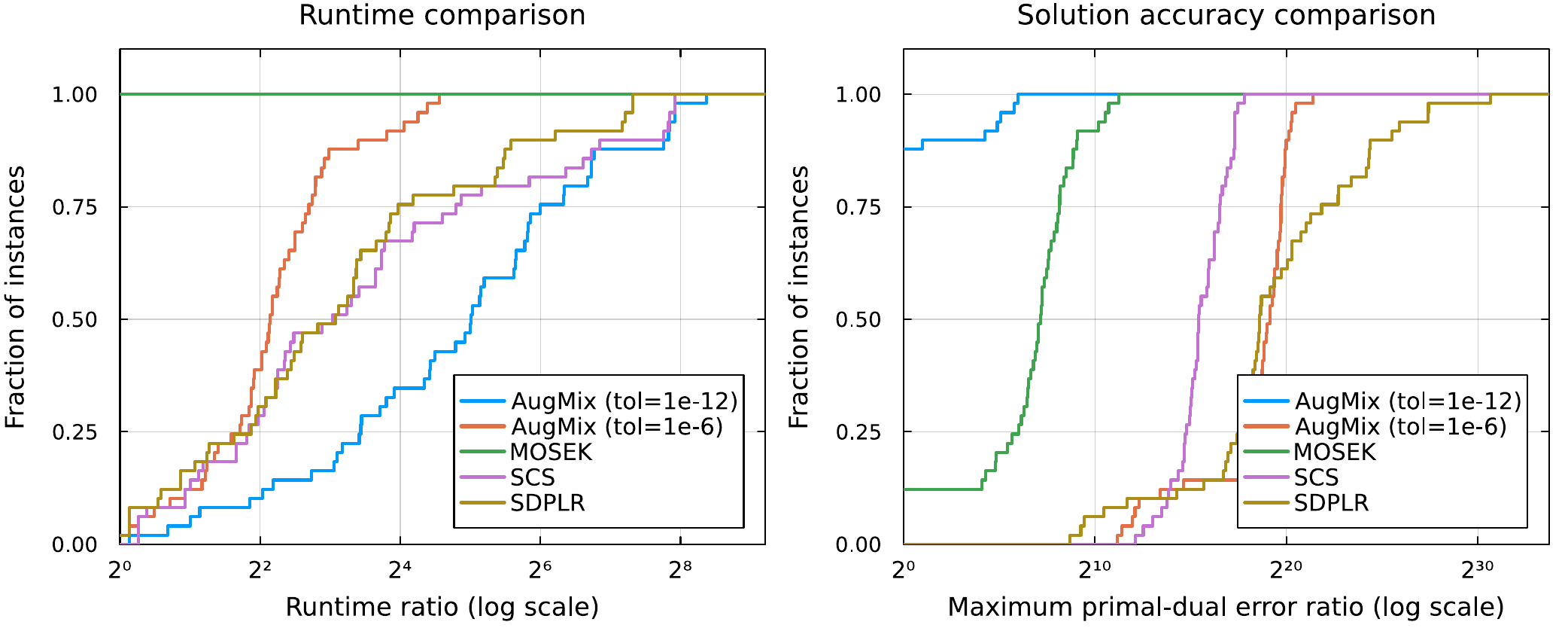}
    \caption{Performance profiles comparing the Augmented Mixing Method (with \texttt{tol=1e-12} and \texttt{tol=1e-6}) against \texttt{MOSEK}, \texttt{SCS}, and \texttt{SDPLR} on the basic SDP relaxation~\eqref{eq:sdp_srflp} of the single-row facility layout problem. Left: Runtime comparison. Right: Maximum primal-dual error comparison.}
    \label{fig:perfplot_srflpBasic}
\end{figure}

The performance profiles in Figure~\ref{fig:perfplot_srflpBasic} indicate that \texttt{MOSEK} is the fastest solver on most instances of the basic relaxation, while the Augmented Mixing Method with \texttt{tol=1e-12} is the slowest. However, the version with \texttt{tol=1e-6} is the second fastest on the majority of instances. In terms of solution accuracy, the Augmented Mixing Method with \texttt{tol=1e-12} outperforms all other solvers, including \texttt{MOSEK}. All other solvers yield solutions of significantly lower accuracy. These results again demonstrate the flexibility of the Augmented Mixing Method to trade runtime for precision.

For the SDP relaxation with triangle inequalities, Table~\ref{table:srflp_met} shows that \texttt{SCS} has shorter runtime on smaller instances, but often exceeds the time limit on larger ones. The Augmented Mixing Method performs fewer than \num{2000} outer iterations in total and often far fewer than for the basic relaxation. However, on the largest instances with $n = 300$ and nearly 18 million inequality constraints, only about $200$ iterations are performed before hitting the time limit, leading to lower-accuracy solutions. The same behavior is observed for \texttt{SCS}. On all other instances, the Augmented Mixing Method achieves substantially better primal-dual accuracy than \texttt{SCS}, often by several orders of magnitude. Most notably, it solves an instance with $n = 276$ and nearly $14$ million inequality constraints to a relative precision of almost $10^{-11}$.

\begin{table}[ht]
\centering
\resizebox{1.0\textwidth}{!}{
\setlength{\tabcolsep}{7pt}
\renewcommand{\arraystretch}{1.1}
\begin{filecontents*}{srflpBasic.csv}
instance,type,seed,statusCode,time,iter,tol,timeLimit,maxIters,itersZ,scaling,warmStart,k,n,numBlocks,avgBlockSize,minBlockSize,maxBlockSize,m,equations,inequalities,nnzA,densityAPercent,timeScaling,timeData,primalObj,dualObj,gap,pinf,dinf,compl,complNotPsdAbs,mu,evals,nnzY,gapOriginal,pinfOriginal,dinfOriginal,complOriginal
N-15_t,srflpBasic,332857678,tol,38.3,8950,1.0e-12,7200.0,100000,50,true,false,34,105,1,105.0,105,105,560,560,0,2835,0.04591836734693878,0.0,0.0008800029754638672,-2332.433785312531,-2332.4337853130774,8.1e-14,9.9e-13,1.1e-13,3.5e-14,2.5e-14,1.2e+01,9661070,560,1.2e-13,1.3e-11,2.0e-12,5.1e-14
N-16a_t,srflpBasic,847097585,tol,26.0,4600,1.0e-12,7200.0,100000,50,true,false,37,120,1,120.0,120,120,680,680,0,3480,0.03553921568627451,0.0,0.001135110855102539,-3043.5588698624983,-3043.5588698624806,2.0e-15,9.3e-13,1.3e-13,2.7e-14,6.3e-15,1.1e+01,5814058,680,3.2e-15,1.3e-11,2.8e-12,3.9e-14
N-16b_t,srflpBasic,608328147,tol,41.1,7350,1.0e-12,7200.0,100000,50,true,false,37,120,1,120.0,120,120,680,680,0,3480,0.03553921568627451,0.0,0.0011670589447021484,-2808.218390066028,-2808.2183900656632,4.5e-14,9.3e-13,1.3e-13,2.2e-14,1.2e-14,1.4e+01,9177704,680,6.6e-14,1.3e-11,2.7e-12,3.2e-14
N-17_t,srflpBasic,199453294,tol,429.5,72350,1.0e-12,7200.0,100000,50,true,false,41,136,1,136.0,136,136,816,816,0,4216,0.02793396770472895,0.0,0.0014841556549072266,-3525.1716496749304,-3525.1716496749304,0.0e+00,3.1e-15,1.0e-12,5.4e-17,9.6e-17,1.7e+00,76336514,816,5.8e-16,4.7e-14,2.3e-11,2.4e-16
N-18_t,srflpBasic,647053474,tol,31.3,3500,1.0e-12,7200.0,100000,50,true,false,45,153,1,153.0,153,153,969,969,0,5049,0.022258645460248083,0.0,0.0019230842590332031,-4075.0109505344985,-4075.010950533741,6.4e-14,9.9e-13,6.0e-14,2.8e-14,1.4e-14,9.5e+00,5705232,969,9.3e-14,1.7e-11,1.5e-12,4.1e-14
N-20_t,srflpBasic,108253285,tol,57.1,4550,1.0e-12,7200.0,100000,50,true,false,52,190,1,190.0,190,190,1330,1330,0,7030,0.014641867827463396,0.0,0.0029277801513671875,-5861.867456417995,-5861.867456419225,7.2e-14,9.5e-13,7.3e-14,1.4e-14,1.9e-15,9.1e+00,9443325,1330,1.1e-13,1.8e-11,2.2e-12,2.0e-14
N-21_t,srflpBasic,733960829,tol,111.8,7500,1.0e-12,7200.0,100000,50,true,false,56,210,1,210.0,210,210,1540,1540,0,8190,0.012059369202226345,0.0,0.003634929656982422,-5818.105309920511,-5818.105309919756,4.4e-14,9.4e-13,1.5e-13,1.2e-14,7.5e-16,1.1e+01,16749334,1540,6.5e-14,2.0e-11,4.4e-12,1.7e-14
N-22_t,srflpBasic,651811358,tol,158.7,8900,1.0e-12,7200.0,100000,50,true,false,60,231,1,231.0,231,231,1771,1771,0,9471,0.010021975013505222,0.0,0.004158973693847656,-7037.089240341744,-7037.0892403418875,6.9e-15,9.5e-13,4.1e-14,3.1e-14,2.0e-14,1.3e+01,21641394,1771,8.3e-15,2.1e-11,1.3e-12,4.6e-14
N-24_t,srflpBasic,956426398,tol,338.3,14350,1.0e-12,7200.0,100000,50,true,false,68,276,1,276.0,276,276,2300,2300,0,12420,0.007088846880907372,0.0,0.006179094314575195,-9125.612727454052,-9125.612727454134,3.0e-15,9.9e-13,1.3e-14,1.8e-14,1.1e-14,6.3e+00,39019579,2300,4.3e-15,2.5e-11,4.7e-13,2.6e-14
N40_1,srflpBasic,278026475,tol,6400.4,22200,1.0e-12,7200.0,100000,50,true,false,147,780,1,780.0,780,780,10660,10660,0,60060,0.0009260595564535527,0.0,0.12987804412841797,-133938.406500141,-133938.40650012656,3.5e-14,9.9e-13,1.1e-14,1.6e-14,5.1e-15,6.0e+00,181380702,10660,5.7e-14,4.3e-11,7.0e-13,2.5e-14
N40_2,srflpBasic,893109454,tol,6651.3,22950,1.0e-12,7200.0,100000,50,true,false,147,780,1,780.0,780,780,10660,10660,0,60060,0.0009260595564535527,0.0,0.2680318355560303,-134819.64817323352,-134819.64817324484,2.8e-14,9.9e-13,1.2e-14,1.0e-14,3.7e-15,1.1e+01,192452850,10660,4.9e-14,4.3e-11,7.4e-13,1.5e-14
N40_3,srflpBasic,566265073,tol,3127.4,10550,1.0e-12,7200.0,100000,50,true,false,147,780,1,780.0,780,780,10660,10660,0,60060,0.0009260595564535527,0.0,0.256619930267334,-123856.82833991759,-123856.82833991696,1.7e-15,9.8e-13,1.4e-13,7.4e-15,3.2e-15,9.9e+00,90960992,10660,2.6e-15,5.2e-11,8.3e-12,1.1e-14
N40_4,srflpBasic,329567457,time,7200.1,24442,1.0e-12,7200.0,100000,50,true,false,147,780,1,780.0,780,780,10660,10660,0,60060,0.0009260595564535527,0.0,0.2601289749145508,-120048.18388792001,-120048.18388775115,4.6e-13,3.8e-11,Inf,Inf,1.1e-13,1.3e+01,212232977,10660,7.0e-13,1.7e-09,3.5e-11,5.8e-13
N40_5,srflpBasic,886122421,tol,3318.4,11200,1.0e-12,7200.0,100000,50,true,false,147,780,1,780.0,780,780,10660,10660,0,60060,0.0009260595564535527,0.0,0.1957409381866455,-140415.01643523836,-140415.01643523076,1.8e-14,9.7e-13,1.7e-13,1.4e-14,7.5e-15,9.4e+00,98102719,10660,3.1e-14,4.2e-11,1.1e-11,2.2e-14
O-9_t,srflpBasic,742573675,tol,1.0,1450,1.0e-12,7200.0,100000,50,true,false,16,36,1,36.0,36,36,120,120,0,540,0.3472222222222222,0.0,0.00011110305786132812,-676.2271960586024,-676.2271960582893,1.6e-13,9.2e-13,1.0e-13,2.4e-14,1.9e-14,8.8e+00,554021,120,2.3e-13,5.7e-12,1.2e-12,3.6e-14
O-10_t,srflpBasic,609398784,tol,1.2,1250,1.0e-12,7200.0,100000,50,true,false,19,45,1,45.0,45,45,165,165,0,765,0.22895622895622897,0.0,0.00015401840209960938,-952.2579921957737,-952.2579921956709,3.7e-14,6.5e-13,1.6e-13,9.5e-15,3.6e-15,6.5e+00,609408,165,5.4e-14,4.9e-12,2.3e-12,1.4e-14
O-15_t,srflpBasic,98521087,tol,9.3,2200,1.0e-12,7200.0,100000,50,true,false,34,105,1,105.0,105,105,560,560,0,2835,0.04591836734693878,0.0,0.0007150173187255859,-3611.154804792492,-3611.154804792764,2.6e-14,8.8e-13,1.2e-13,3.3e-14,1.9e-14,5.2e+00,2588182,560,3.9e-14,9.3e-12,3.4e-12,4.8e-14
O-20_t,srflpBasic,934582034,tol,152.9,15600,1.0e-12,7200.0,100000,50,true,false,52,190,1,190.0,190,190,1330,1330,0,7030,0.014641867827463396,0.0,0.0023429393768310547,-8797.526658671693,-8797.526658671624,2.7e-15,9.9e-13,1.3e-14,3.0e-14,1.8e-14,8.8e-01,23312720,1330,5.9e-15,1.9e-11,6.0e-13,4.4e-14
S-12_t,srflpBasic,143649825,tol,4.0,2150,1.0e-12,7200.0,100000,50,true,false,24,66,1,66.0,66,66,286,286,0,1386,0.11125238397965671,0.0,0.00034308433532714844,-2321.745503772384,-2321.7455037725804,2.8e-14,9.4e-13,5.3e-14,3.4e-14,1.4e-14,6.4e+00,1673944,286,4.2e-14,9.0e-12,1.1e-12,5.1e-14
S-13_t,srflpBasic,50638069,tol,6.3,2250,1.0e-12,7200.0,100000,50,true,false,27,78,1,78.0,78,78,364,364,0,1794,0.0810087348548887,0.0,0.0004420280456542969,-3117.650931314446,-3117.650931314322,1.3e-14,7.4e-13,3.0e-13,2.3e-14,9.4e-15,4.7e+00,2203952,364,1.8e-14,7.9e-12,6.8e-12,3.5e-14
S-14_t,srflpBasic,871874611,tol,10.9,3100,1.0e-12,7200.0,100000,50,true,false,31,91,1,91.0,91,91,455,455,0,2275,0.06037918125830214,0.0,0.0005340576171875,-3947.330731806439,-3947.330731806699,2.2e-14,1.0e-12,3.9e-14,2.9e-14,1.8e-14,6.3e+00,3311523,455,3.3e-14,1.2e-11,9.8e-13,4.3e-14
S-15_t,srflpBasic,505015922,tol,14.7,3300,1.0e-12,7200.0,100000,50,true,false,34,105,1,105.0,105,105,560,560,0,2835,0.04591836734693878,0.0,0.0007159709930419922,-4873.579395170685,-4873.579395170616,4.7e-15,8.9e-13,4.4e-14,2.0e-14,2.7e-15,5.7e+00,4071934,560,5.5e-15,1.2e-11,1.2e-12,3.0e-14
S-16_t,srflpBasic,162807843,tol,11.8,1950,1.0e-12,7200.0,100000,50,true,false,37,120,1,120.0,120,120,680,680,0,3480,0.03553921568627451,0.0,0.0009291172027587891,-6024.864761102412,-6024.864761102534,6.9e-15,9.1e-13,7.4e-14,2.6e-14,1.1e-14,5.4e+00,2908291,680,1.1e-14,1.3e-11,2.3e-12,3.9e-14
S-17_t,srflpBasic,843375701,tol,18.6,2550,1.0e-12,7200.0,100000,50,true,false,41,136,1,136.0,136,136,816,816,0,4216,0.02793396770472895,0.0,0.001177072525024414,-7257.7203099886565,-7257.720309988784,5.9e-15,7.9e-13,2.8e-14,1.4e-14,7.4e-15,6.3e+00,4170354,816,8.3e-15,1.2e-11,9.7e-13,2.0e-14
S-18_t,srflpBasic,991247657,tol,29.5,3200,1.0e-12,7200.0,100000,50,true,false,45,153,1,153.0,153,153,969,969,0,5049,0.022258645460248083,0.0,0.0015060901641845703,-8704.038818267449,-8704.038818266377,4.1e-14,9.3e-13,6.6e-13,7.9e-15,2.8e-14,4.0e+00,5884751,969,6.2e-14,1.6e-11,2.5e-11,1.2e-14
S-19_t,srflpBasic,156478564,tol,22.7,2150,1.0e-12,7200.0,100000,50,true,false,48,171,1,171.0,171,171,1140,1140,0,5985,0.017954242330973633,0.0,0.0018610954284667969,-10260.894097267485,-10260.89409726841,3.0e-14,8.4e-13,1.6e-14,2.0e-14,1.5e-14,3.2e+00,4632445,1140,4.4e-14,1.5e-11,6.4e-13,3.0e-14
S-20_t,srflpBasic,319730356,tol,65.0,5500,1.0e-12,7200.0,100000,50,true,false,52,190,1,190.0,190,190,1330,1330,0,7030,0.014641867827463396,0.0,0.0023589134216308594,-11969.202081609537,-11969.202081608857,1.9e-14,9.6e-13,5.0e-14,4.7e-14,3.9e-15,1.3e+00,11278065,1330,2.6e-14,1.9e-11,2.2e-12,7.1e-14
S-21_t,srflpBasic,236440055,tol,39.7,2450,1.0e-12,7200.0,100000,50,true,false,56,210,1,210.0,210,210,1540,1540,0,8190,0.012059369202226345,0.0,0.002925872802734375,-13767.893577137207,-13767.893577136903,7.4e-15,8.3e-13,2.1e-14,8.6e-15,1.3e-14,5.8e+00,6425751,1540,1.2e-14,1.7e-11,9.8e-13,1.3e-14
S-22_t,srflpBasic,82659411,tol,125.1,7000,1.0e-12,7200.0,100000,50,true,false,60,231,1,231.0,231,231,1771,1771,0,9471,0.010021975013505222,0.0,0.003688812255859375,-15884.72625271469,-15884.726252714168,1.1e-14,1.0e-12,2.0e-14,1.6e-14,3.1e-15,2.9e+00,17923655,1771,1.6e-14,2.2e-11,1.0e-12,2.3e-14
S-23_t,srflpBasic,63015030,tol,76.4,3400,1.0e-12,7200.0,100000,50,true,false,64,253,1,253.0,253,253,2024,2024,0,10879,0.008397256635785593,0.0,0.0044329166412353516,-18247.42046731529,-18247.420467314798,9.1e-15,8.4e-13,2.3e-14,1.0e-14,2.5e-15,5.7e+00,10506272,2024,1.5e-14,2.0e-11,1.2e-12,1.5e-14
S-24_t,srflpBasic,399145894,tol,93.3,3450,1.0e-12,7200.0,100000,50,true,false,68,276,1,276.0,276,276,2300,2300,0,12420,0.007088846880907372,0.0,0.0054361820220947266,-20794.552401250712,-20794.55240124947,2.0e-14,9.1e-13,7.1e-14,6.8e-15,3.6e-15,6.5e+00,11767759,2300,2.8e-14,2.3e-11,4.1e-12,1.0e-14
S-25_t,srflpBasic,334933276,tol,126.5,3400,1.0e-12,7200.0,100000,50,true,false,73,300,1,300.0,300,300,2600,2600,0,14100,0.006025641025641026,0.0,0.006510019302368164,-23478.21460909972,-23478.214609099145,8.3e-15,8.2e-13,2.9e-14,1.4e-14,9.3e-15,5.6e+00,12560465,2600,1.6e-14,2.2e-11,1.8e-12,2.1e-14
ste36-1,srflpBasic,744155305,time,7200.2,43342,1.0e-12,7200.0,100000,50,true,false,125,630,1,630.0,630,630,7770,7770,0,43470,0.0014095728381442667,0.0,0.06640791893005371,-37162.85998367685,-37162.85996375688,1.4e-10,2.3e-09,Inf,Inf,2.6e-11,4.3e+01,301093076,7770,2.7e-10,1.1e-07,9.8e-09,5.6e-11
ste36-2,srflpBasic,108410693,time,7200.1,42929,1.0e-12,7200.0,100000,50,true,false,125,630,1,630.0,630,630,7770,7770,0,43470,0.0014095728381442667,0.0,0.050924062728881836,-660675.8296075613,-660675.829318479,1.1e-10,3.9e-09,Inf,Inf,6.8e-11,2.7e+01,300028398,7770,2.2e-10,1.5e-07,4.1e-09,1.8e-10
ste36-3,srflpBasic,189712974,time,7200.1,42567,1.0e-12,7200.0,100000,50,true,false,125,630,1,630.0,630,630,7770,7770,0,43470,0.0014095728381442667,0.0,0.052643775939941406,-415761.12833763164,-415761.12837458943,2.3e-11,2.9e-09,Inf,Inf,3.2e-11,4.8e+01,302784608,7770,4.4e-11,1.1e-07,3.6e-10,4.5e-12
ste36-4,srflpBasic,823618862,time,7200.2,44472,1.0e-12,7200.0,100000,50,true,false,125,630,1,630.0,630,630,7770,7770,0,43470,0.0014095728381442667,0.0,0.05269289016723633,-372599.6854585339,-372599.6851910553,1.8e-10,2.1e-09,Inf,Inf,3.4e-11,1.8e+01,303649980,7770,3.6e-10,9.8e-08,1.3e-09,8.8e-12
ste36-5,srflpBasic,308963030,time,7200.1,44406,1.0e-12,7200.0,100000,50,true,false,125,630,1,630.0,630,630,7770,7770,0,43470,0.0014095728381442667,0.0,0.06960701942443848,-359342.2140655078,-359342.21396064066,7.5e-11,2.7e-09,Inf,Inf,1.4e-11,3.8e+01,302401353,7770,1.5e-10,1.2e-07,2.8e-08,9.1e-11
Y-9_t,srflpBasic,25564899,tol,0.1,100,1.0e-12,7200.0,100000,50,true,false,16,36,1,36.0,36,36,120,120,0,540,0.3472222222222222,0.0,0.00011801719665527344,-1336.500000000024,-1336.4999999999916,8.0e-15,2.3e-13,4.1e-17,7.8e-18,2.0e-15,5.8e-01,33440,120,1.2e-14,1.2e-12,4.0e-15,8.3e-17
Y-10_t,srflpBasic,748596922,tol,0.6,550,1.0e-12,7200.0,100000,50,true,false,19,45,1,45.0,45,45,165,165,0,765,0.22895622895622897,0.0,0.00015592575073242188,-1732.5380839029785,-1732.538083903193,4.1e-14,5.3e-13,2.1e-13,3.7e-14,2.7e-14,3.0e+00,321757,165,6.2e-14,4.0e-12,3.4e-12,5.6e-14
Y-11_t,srflpBasic,394107148,tol,9.9,8000,1.0e-12,7200.0,100000,50,true,false,21,55,1,55.0,55,55,220,220,0,1045,0.15702479338842976,0.0,0.000225067138671875,-2025.253596054061,-2025.2535960538357,3.7e-14,8.8e-13,6.6e-14,3.1e-14,6.8e-15,5.7e+00,4208928,220,5.6e-14,6.9e-12,1.3e-12,4.6e-14
Y-12_t,srflpBasic,509513661,tol,2.6,1350,1.0e-12,7200.0,100000,50,true,false,24,66,1,66.0,66,66,286,286,0,1386,0.11125238397965671,0.0,0.0003199577331542969,-2440.675416205119,-2440.6754162051075,1.5e-15,8.1e-13,3.3e-14,2.1e-14,5.7e-15,4.4e+00,1131177,286,1.6e-15,7.8e-12,7.2e-13,3.1e-14
Y-13_t,srflpBasic,272940807,tol,3.5,1300,1.0e-12,7200.0,100000,50,true,false,27,78,1,78.0,78,78,364,364,0,1794,0.0810087348548887,0.0,0.00043010711669921875,-2950.218681613301,-2950.2186816129038,4.5e-14,8.1e-13,7.5e-14,1.9e-14,1.2e-14,5.7e+00,1201382,364,6.8e-14,8.6e-12,1.8e-12,2.9e-14
Y-14_t,srflpBasic,973254063,tol,8.9,2600,1.0e-12,7200.0,100000,50,true,false,31,91,1,91.0,91,91,455,455,0,2275,0.06037918125830214,0.0,0.0005831718444824219,-3365.9089672500872,-3365.908967250199,1.1e-14,2.9e-13,8.6e-13,1.1e-14,1.6e-15,2.9e+00,2678868,455,1.7e-14,3.4e-12,2.2e-11,1.6e-14
Y-15_t,srflpBasic,95155869,tol,9.0,2050,1.0e-12,7200.0,100000,50,true,false,34,105,1,105.0,105,105,560,560,0,2835,0.04591836734693878,0.0,0.0007190704345703125,-3939.0560174983316,-3939.0560174983775,3.9e-15,9.0e-13,9.6e-14,1.9e-14,9.6e-15,6.2e+00,2490178,560,7.1e-15,1.2e-11,2.7e-12,2.8e-14
Y-20_t,srflpBasic,845781852,tol,103.3,8350,1.0e-12,7200.0,100000,50,true,false,52,190,1,190.0,190,190,1330,1330,0,7030,0.014641867827463396,0.0,0.0023589134216308594,-7468.7970479665555,-7468.797047965854,3.2e-14,9.7e-13,1.2e-13,1.8e-14,9.9e-15,1.2e+01,18262106,1330,4.8e-14,1.9e-11,4.4e-12,2.6e-14
Y-25_t,srflpBasic,622848081,tol,186.2,5200,1.0e-12,7200.0,100000,50,true,false,73,300,1,300.0,300,300,2600,2600,0,14100,0.006025641025641026,0.0,0.006701946258544922,-12193.636920678286,-12193.636920677904,1.1e-14,3.1e-13,1.0e-12,2.4e-15,4.4e-15,2.3e+00,16788761,2600,1.6e-14,8.3e-12,4.0e-11,3.6e-15
Y-30_t,srflpBasic,544680419,tol,778.6,10050,1.0e-12,7200.0,100000,50,true,false,95,435,1,435.0,435,435,4495,4495,0,24795,0.0029151164128725406,0.0,0.02110910415649414,-17114.39681214244,-17114.396812141553,1.8e-14,9.9e-13,2.6e-14,1.1e-14,8.0e-15,2.4e+00,46751853,4495,2.9e-14,3.5e-11,1.6e-12,1.5e-14
Y-35_t,srflpBasic,315344506,tol,701.5,4800,1.0e-12,7200.0,100000,50,true,false,120,595,1,595.0,595,595,7140,7140,0,39865,0.001577101428806817,0.0,0.057541847229003906,-24306.321164231376,-24306.321164229954,2.0e-14,9.4e-13,1.1e-13,8.8e-15,3.5e-15,4.7e+00,34047937,7140,2.9e-14,4.1e-11,7.4e-12,1.3e-14
Y-40_t,srflpBasic,775044439,tol,3029.1,11800,1.0e-12,7200.0,100000,50,true,false,147,780,1,780.0,780,780,10660,10660,0,60060,0.0009260595564535527,0.0,0.29473209381103516,-30169.022895703016,-30169.022895707123,4.6e-14,9.0e-13,1.4e-14,6.3e-14,2.3e-14,6.2e-01,80179231,10660,6.9e-14,4.8e-11,1.1e-12,9.3e-14
\end{filecontents*}
\pgfplotstabletypeset[
    col sep=comma,
    string type,
    columns={instance,n,equations,inequalities,statusCode,time,iter,gapOriginal,pinfOriginal,dinfOriginal,complOriginal},
    columns/instance/.style={column name=\textbf{Instance}, column type=l},
    columns/n/.style={column name=$n$, column type=r},
    columns/equations/.style={column name=$m_a$, column type=r},
    columns/inequalities/.style={column name=$m_b$, column type=r},
    columns/statusCode/.style={column name=\textbf{Status}, column type=c},
    columns/time/.style={column name=\textbf{Time~[s]}, column type=r},
    columns/iter/.style={column name=\textbf{Iter}, column type=r},
    columns/gapOriginal/.style={column name=\texttt{gap}, column type=c},
    columns/pinfOriginal/.style={column name=\texttt{pinf}, column type=c},
    columns/dinfOriginal/.style={column name=\texttt{dinf}, column type=c},
    columns/complOriginal/.style={column name=\texttt{compl}, column type=c},
    every head row/.style={before row=\toprule, after row=\midrule},
    every last row/.style={after row=\bottomrule}
]{srflpBasic.csv}
}
\caption{Results of the Augmented Mixing Method with \texttt{tol=1e-12} on instances of the basic SDP relaxation~\eqref{eq:sdp_srflp} for the single-row facility layout problem.}
\label{table:srflp_basic}
\end{table}

\begin{table}[ht]
\centering
\resizebox{1.0\textwidth}{!}{
\setlength{\tabcolsep}{7pt}
\renewcommand{\arraystretch}{1.1}
\begin{filecontents*}{srflpMet.csv}
instance,type,seed,statusCode,time,iter,tol,timeLimit,maxIters,itersZ,scaling,warmStart,k,n,numBlocks,avgBlockSize,minBlockSize,maxBlockSize,m,equations,inequalities,nnzA,densityAPercent,timeScaling,timeData,primalObj,dualObj,gap,pinf,dinf,compl,complNotPsdAbs,mu,evals,nnzY,gapOriginal,pinfOriginal,dinfOriginal,complOriginal,timeSCS,maxErrorSCS
N-15_t,srflpMet,317668309,tol,972.7,1800,1.0e-12,7200.0,100000,50,true,false,105,105,1,105.0,105,105,749945,560,749385,4499145,0.05441541903549117,0.8,2.3403258323669434,-2268.999999999584,-2269.0000000366,4.7e-13,8.2e-13,7.8e-15,4.7e-13,4.7e-13,6.2e-01,2965516,556918,8.2e-12,3.6e-10,1.4e-13,8.1e-12,471.6,3.0e-08
N-16a_t,srflpMet,90854801,tol,410.8,500,1.0e-12,7200.0,100000,50,true,false,120,120,1,120.0,120,120,1123480,680,1122800,6740280,0.04166295795207747,2.4,3.682466983795166,-3013.9999999990473,-3014.0000000064147,6.4e-14,2.9e-13,0.0e+00,5.6e-14,5.6e-14,5.5e-01,920351,828875,1.2e-12,1.6e-10,3.5e-15,1.1e-12,175.6,6.6e-10
N-16b_t,srflpMet,314577408,tol,333.1,350,1.0e-12,7200.0,100000,50,true,false,120,120,1,120.0,120,120,1123480,680,1122800,6740280,0.04166295795207747,2.2,3.6110260486602783,-2784.0000000002565,-2783.9999999977294,2.4e-14,1.6e-13,1.3e-15,1.6e-19,2.1e-14,5.9e-01,758474,835764,4.5e-13,8.7e-11,3.3e-14,3.1e-17,229.4,2.0e-08
N-17_t,srflpMet,900624326,tol,662.8,500,1.0e-12,7200.0,100000,50,true,false,136,136,1,136.0,136,136,1640296,816,1639480,9841096,0.03243720502095891,2.3,6.114405155181885,-3463.000000000351,-3463.0000000280356,1.9e-13,7.0e-14,0.0e+00,1.9e-13,1.9e-13,5.7e-01,1100869,1230420,4.0e-12,3.6e-11,3.5e-15,4.0e-12,660.3,8.7e-08
N-18_t,srflpMet,265520487,tol,1223.9,700,1.0e-12,7200.0,100000,50,true,false,153,153,1,153.0,153,153,2341257,969,2340288,14046777,0.02562977167946158,3.4,8.713390111923218,-3987.999999996074,-3988.000000155993,8.6e-13,7.1e-14,0.0e+00,8.4e-13,8.4e-13,6.1e-01,1580231,1755664,2.0e-11,4.4e-11,5.5e-15,2.0e-11,1355.6,3.9e-07
N-20_t,srflpMet,81113020,time,7201.1,1864,1.0e-12,7200.0,100000,50,true,false,190,190,1,190.0,190,190,4500910,1330,4499580,27004510,0.01661991393787151,8.7,19.253477096557617,-5718.000000868788,-5718.000040050183,1.3e-10,3.7e-11,Inf,Inf,1.3e-10,6.4e-01,6012922,3195671,3.4e-09,3.9e-08,8.2e-10,3.5e-09,7646.7,3.6e-04
N-21_t,srflpMet,121811359,time,7203.1,1329,1.0e-12,7200.0,100000,50,true,false,210,210,1,210.0,210,210,6086290,1540,6084750,36516690,0.01360505097757184,14.2,27.71642017364502,-5689.000000039179,-5689.000000834113,2.3e-12,1.2e-12,Inf,Inf,2.4e-12,6.3e-01,4778940,4389208,7.0e-11,1.4e-09,1.4e-10,7.5e-11,7845.9,3.6e-04
N-22_t,srflpMet,830128640,time,7201.4,770,1.0e-12,7200.0,100000,50,true,false,231,231,1,231.0,231,231,8111411,1771,8109640,48667311,0.011243900241757705,19.6,42.13645100593567,-6884.000000548688,-6884.000037894782,8.5e-11,2.2e-11,Inf,Inf,8.6e-11,6.8e-01,3611374,6041928,2.7e-09,3.1e-08,7.0e-13,2.8e-09,8296.4,7.3e-04
N-24_t,srflpMet,726278586,time,7200.6,329,1.0e-12,7200.0,100000,50,true,false,276,276,1,276.0,276,276,13864676,2300,13862376,83186676,0.007876365871771198,33.9,76.41849398612976,-8879.201381161112,-8809.091528514937,1.1e-04,8.8e-06,Inf,Inf,1.1e-04,1.1e+00,2342713,1182677,4.0e-03,1.3e-02,1.1e-03,2.0e-04,10667.4,1.3e-03
O-9_t,srflpMet,169458414,tol,5.5,350,1.0e-12,7200.0,100000,50,true,false,36,36,1,36.0,36,36,28596,120,28476,171396,0.4624772695481886,0.0,0.02368307113647461,-669.0000000000001,-668.9999999999973,2.7e-16,3.1e-16,2.9e-17,4.6e-19,2.6e-16,6.2e-01,144636,21487,2.2e-15,2.1e-14,2.5e-15,2.9e-17,0.6,1.4e-08
O-10_t,srflpMet,325389824,tol,11.6,350,1.0e-12,7200.0,100000,50,true,false,45,45,1,45.0,45,45,56805,165,56640,340605,0.29610069536132383,0.1,0.03731513023376465,-937.9999999999998,-938.0000000000001,2.5e-17,4.3e-18,0.0e+00,7.5e-18,8.3e-18,5.4e-01,199471,42651,1.8e-16,1.8e-14,1.3e-15,1.1e-16,1.4,6.6e-08
O-15_t,srflpMet,565449994,tol,211.1,450,1.0e-12,7200.0,100000,50,true,false,105,105,1,105.0,105,105,749945,560,749385,4499145,0.05441541903549117,1.1,2.5864322185516357,-3581.0000000017703,-3581.000000107253,8.6e-13,1.7e-13,0.0e+00,8.8e-13,8.8e-13,5.8e-01,631400,561583,1.5e-11,6.0e-11,6.0e-15,1.5e-11,112.2,2.7e-07
O-20_t,srflpMet,817412986,tol,1301.9,350,1.0e-12,7200.0,100000,50,true,false,190,190,1,190.0,190,190,4500910,1330,4499580,27004510,0.01661991393787151,10.0,23.5545711517334,-8675.99999999839,-8676.000000035392,7.8e-14,2.1e-14,0.0e+00,7.5e-14,7.5e-14,5.5e-01,1069877,3375208,2.1e-12,1.8e-11,1.4e-14,2.0e-12,3455.4,9.2e-08
S-12_t,srflpMet,967214167,tol,54.5,400,1.0e-12,7200.0,100000,50,true,false,66,66,1,66.0,66,66,183106,286,182820,1098306,0.13769967321462084,0.3,0.42228102684020996,-2312.9999999999773,-2312.99999999997,1.2e-16,1.2e-15,4.0e-17,1.7e-20,4.8e-16,5.7e-01,432724,137403,1.5e-15,2.1e-13,3.6e-15,1.6e-17,9.1,1.1e-07
S-13_t,srflpMet,891148347,tol,129.2,600,1.0e-12,7200.0,100000,50,true,false,78,78,1,78.0,78,78,304382,364,304018,1825902,0.09859826948222451,0.4,0.7939350605010986,-3105.4999999999804,-3105.4999999979514,2.2e-14,3.8e-15,1.7e-15,2.8e-19,2.2e-14,5.7e-01,722607,228395,3.3e-13,8.5e-13,4.1e-14,1.9e-17,20.6,6.1e-08
S-14_t,srflpMet,890599987,tol,314.5,900,1.0e-12,7200.0,100000,50,true,false,91,91,1,91.0,91,91,486031,455,485576,2915731,0.07244371266419243,0.9,1.3655359745025635,-3925.9999999991173,-3926.0000000090768,7.7e-14,3.5e-13,0.0e+00,7.0e-14,7.0e-14,5.6e-01,1291815,356760,1.3e-12,1.2e-10,4.4e-15,1.2e-12,100.8,1.5e-07
S-15_t,srflpMet,901605809,tol,219.1,350,1.0e-12,7200.0,100000,50,true,false,105,105,1,105.0,105,105,749945,560,749385,4499145,0.05441541903549117,1.4,2.22979998588562,-4858.000000004969,-4858.0000000825385,4.4e-13,2.2e-13,0.0e+00,4.7e-13,4.7e-13,5.5e-01,664153,562621,8.0e-12,7.8e-11,6.4e-15,8.5e-12,197.6,7.6e-07
S-16_t,srflpMet,326330585,tol,660.0,800,1.0e-12,7200.0,100000,50,true,false,120,120,1,120.0,120,120,1123480,680,1122800,6740280,0.04166295795207747,2.2,4.100729942321777,-5996.999999995249,-5997.000000144435,6.2e-13,6.8e-14,0.0e+00,6.0e-13,6.0e-13,5.9e-01,1491214,842422,1.2e-11,2.9e-11,7.1e-15,1.2e-11,913.5,1.9e-09
S-17_t,srflpMet,461948314,tol,882.0,700,1.0e-12,7200.0,100000,50,true,false,136,136,1,136.0,136,136,1640296,816,1639480,9841096,0.03243720502095891,3.4,6.8429789543151855,-7219.499999999534,-7219.50000002847,9.1e-14,2.0e-14,0.0e+00,8.9e-14,8.9e-14,5.7e-01,1493503,1230443,2.0e-12,1.0e-11,1.9e-14,2.0e-12,893.4,7.9e-08
S-18_t,srflpMet,553918524,tol,1258.9,650,1.0e-12,7200.0,100000,50,true,false,153,153,1,153.0,153,153,2341257,969,2340288,14046777,0.02562977167946158,4.2,10.15283203125,-8664.000000000367,-8664.000000038151,9.0e-14,4.0e-13,0.0e+00,9.1e-14,9.1e-14,5.6e-01,1726640,1744789,2.2e-12,2.5e-10,1.1e-14,2.2e-12,1262.6,3.6e-08
S-19_t,srflpMet,106170185,tol,3508.7,1550,1.0e-12,7200.0,100000,50,true,false,171,171,1,171.0,171,171,3275391,1140,3274251,19651491,0.020518241381836468,7.4,15.289210796356201,-10189.500000002367,-10189.50000006596,1.2e-13,1.2e-14,0.0e+00,1.2e-13,1.2e-13,5.9e-01,3616837,2451846,3.1e-12,8.9e-12,2.1e-14,3.2e-12,3996.2,4.3e-09
S-20_t,srflpMet,346013764,tol,5792.4,1650,1.0e-12,7200.0,100000,50,true,false,190,190,1,190.0,190,190,4500910,1330,4499580,27004510,0.01661991393787151,11.7,22.52543807029724,-11874.999999999549,-11875.000000447122,6.6e-13,8.8e-13,0.0e+00,6.6e-13,6.6e-13,5.8e-01,4847800,3354899,1.9e-11,7.6e-10,1.1e-14,1.9e-11,7207.8,2.3e-08
S-21_t,srflpMet,832875909,time,7204.1,936,1.0e-12,7200.0,100000,50,true,false,210,210,1,210.0,210,210,6086290,1540,6084750,36516690,0.01360505097757184,16.6,33.79444193840027,-13675.499975108862,-13675.500409849088,5.2e-10,1.5e-10,Inf,Inf,4.9e-10,7.8e-01,4799785,2050953,1.6e-08,1.8e-07,3.0e-10,1.5e-08,7856.4,3.1e-04
S-22_t,srflpMet,105726206,time,7204.8,758,1.0e-12,7200.0,100000,50,true,false,231,231,1,231.0,231,231,8111411,1771,8109640,48667311,0.011243900241757705,23.0,47.087352991104126,-15756.000012397963,-15756.001098932227,1.0e-09,8.0e-10,Inf,Inf,1.1e-09,7.5e-01,3748061,4679937,3.4e-08,9.3e-07,4.0e-14,3.5e-08,8320.3,1.3e+01
S-23_t,srflpMet,229029560,time,7206.1,495,1.0e-12,7200.0,100000,50,true,false,253,253,1,253.0,253,253,10668757,2024,10666733,64011277,0.009373496585614703,22.1,58.34909009933472,-18128.99993985597,-18128.99968182446,2.0e-10,3.3e-10,Inf,Inf,2.5e-10,8.1e+00,2894851,3321757,7.1e-09,4.4e-07,5.7e-10,1.2e-10,9205.5,8.7e-04
S-24_t,srflpMet,262629756,tol,7124.1,500,1.0e-12,7200.0,100000,50,true,false,276,276,1,276.0,276,276,13864676,2300,13862376,83186676,0.007876365871771198,33.5,75.55810403823853,-20678.000000002994,-20678.00000106791,6.8e-13,2.1e-14,0.0e+00,6.8e-13,6.8e-13,6.6e-01,2105295,10398727,2.6e-11,3.1e-11,1.0e-13,2.6e-11,10824.6,1.4e-03
S-25_t,srflpMet,680420659,time,7223.8,205,1.0e-12,7200.0,100000,50,true,false,300,300,1,300.0,300,300,17820700,2600,17818100,106922700,0.006666573142469151,45.1,106.89261794090271,-23285.891973227597,-23262.375326102494,1.3e-05,4.2e-06,Inf,Inf,1.3e-05,1.4e+00,1847130,914754,5.1e-04,7.2e-03,3.4e-03,7.8e-05,13454.8,2.5e-03
Y-9_t,srflpMet,385653447,tol,2.9,150,1.0e-12,7200.0,100000,50,true,false,36,36,1,36.0,36,36,28596,120,28476,171396,0.4624772695481886,0.0,0.012660026550292969,-1336.4999999999989,-1336.499999999999,1.2e-17,2.3e-16,4.5e-18,3.3e-19,4.0e-17,6.0e-01,84353,21495,1.7e-16,1.7e-14,2.7e-15,1.8e-17,0.2,4.4e-08
Y-10_t,srflpMet,332292333,tol,13.9,350,1.0e-12,7200.0,100000,50,true,false,45,45,1,45.0,45,45,56805,165,56640,340605,0.29610069536132383,0.0,0.03267097473144531,-1712.0000000000061,-1711.9999999999816,7.1e-16,1.4e-15,5.3e-17,9.6e-20,5.2e-16,5.8e-01,253154,42657,6.6e-15,1.4e-13,2.0e-15,1.5e-17,1.3,1.2e-08
Y-11_t,srflpMet,293024053,tol,16.0,200,1.0e-12,7200.0,100000,50,true,false,55,55,1,55.0,55,55,104995,220,104775,629695,0.19826052322837884,0.1,0.15878701210021973,-2022.4999999999861,-2022.499999999555,9.2e-15,8.8e-16,8.7e-16,3.5e-19,9.5e-15,7.2e-01,187196,78825,1.1e-13,1.2e-13,2.2e-14,2.4e-17,3.3,6.6e-08
Y-12_t,srflpMet,105754475,tol,170.4,1450,1.0e-12,7200.0,100000,50,true,false,66,66,1,66.0,66,66,183106,286,182820,1098306,0.13769967321462084,0.1,0.39978885650634766,-2424.9999999953925,-2425.000000059152,9.9e-13,8.4e-13,0.0e+00,9.2e-13,9.2e-13,5.5e-01,1391558,135233,1.3e-11,1.8e-10,1.7e-15,1.2e-11,43.5,3.8e-08
Y-13_t,srflpMet,522673167,tol,124.1,550,1.0e-12,7200.0,100000,50,true,false,78,78,1,78.0,78,78,304382,364,304018,1825902,0.09859826948222451,0.4,0.8398549556732178,-2920.499999998183,-2920.4999999958495,2.7e-14,3.9e-13,3.8e-15,2.8e-20,4.8e-14,5.7e-01,718903,225484,4.0e-13,1.1e-10,1.0e-13,1.7e-17,32.9,4.6e-08
Y-14_t,srflpMet,613771449,tol,431.7,1400,1.0e-12,7200.0,100000,50,true,false,91,91,1,91.0,91,91,486031,455,485576,2915731,0.07244371266419243,0.7,1.447192907333374,-3340.9999999999914,-3341.0000000020877,1.9e-14,3.0e-15,0.0e+00,1.9e-14,1.9e-14,5.4e-01,1721880,364568,3.1e-13,8.5e-13,5.0e-15,3.1e-13,237.2,1.2e-07
Y-15_t,srflpMet,200947313,tol,288.4,500,1.0e-12,7200.0,100000,50,true,false,105,105,1,105.0,105,105,749945,560,749385,4499145,0.05441541903549117,1.6,2.359494924545288,-3921.000000001825,-3920.9999999958095,4.2e-14,4.4e-13,1.9e-15,6.7e-20,2.9e-14,5.8e-01,883248,557290,7.7e-13,1.9e-10,6.0e-14,1.6e-17,75.1,1.1e-08
Y-20_t,srflpMet,874667527,tol,3677.9,950,1.0e-12,7200.0,100000,50,true,false,190,190,1,190.0,190,190,4500910,1330,4499580,27004510,0.01661991393787151,11.5,26.24228000640869,-7435.000000015401,-7435.000000369153,8.7e-13,9.5e-13,0.0e+00,9.1e-13,9.1e-13,5.8e-01,3092086,3322841,2.4e-11,1.0e-09,1.4e-14,2.5e-11,4614.6,8.8e-07
Y-25_t,srflpMet,565349202,time,7207.5,192,1.0e-12,7200.0,100000,50,true,false,300,300,1,300.0,300,300,17820700,2600,17818100,106922700,0.006666573142469151,45.2,106.37434315681458,-12030.71762238126,-11968.029587408455,6.6e-05,7.6e-06,Inf,Inf,6.6e-05,2.0e+00,1864150,280288,2.6e-03,1.3e-02,1.4e-03,1.6e-04,12692.3,2.0e-03
Geom. mean,srflpMet,,,619.1,,,,,,,,,,,,,,,,,,,,,,,,,,,,,,,,,,,372.6,
\end{filecontents*}
\pgfplotstabletypeset[
    col sep=comma,
    string type,
    columns={instance,n,equations,inequalities,statusCode,time,iter,gapOriginal,pinfOriginal,dinfOriginal,complOriginal,timeSCS,maxErrorSCS},
    columns/instance/.style={column name=\textbf{Instance}, column type=l},
    columns/n/.style={column name=$n$, column type=r},
    columns/equations/.style={column name=$m_a$, column type=r},
    columns/inequalities/.style={column name=$m_b$, column type=r},
    columns/statusCode/.style={column name=\textbf{Status}, column type=|c},
    columns/time/.style={column name=\textbf{Time~[s]}, column type=r},
    columns/iter/.style={column name=\textbf{Iter}, column type=r},
    columns/gapOriginal/.style={column name=\texttt{gap}, column type=c},
    columns/pinfOriginal/.style={column name=\texttt{pinf}, column type=c},
    columns/dinfOriginal/.style={column name=\texttt{dinf}, column type=c},
    columns/complOriginal/.style={column name=\texttt{compl}, column type=c|},
    columns/timeSCS/.style={column name=\textbf{Time~[s]}, column type=r},
    columns/maxErrorSCS/.style={column name=\textbf{Max. error}, column type=c},
    every head row/.style={
        before row={
            \toprule
            & & & & \multicolumn{7}{c|}{\textbf{Augmented Mixing Method}} & \multicolumn{2}{c}{\textbf{SCS}} \\
        },
        after row=\midrule
    },
    every last row/.style={before row=\midrule, after row=\bottomrule}
]{srflpMet.csv}
}
\caption{Results of the Augmented Mixing Method with \texttt{tol=1e-12} and \texttt{SCS}~\cite{brendan2016conic} on various instances of the SRFLP relaxation~\eqref{eq:sdp_srflp}, strengthened by all triangle inequalities~\eqref{eq:triangle_inequalities}. The last row reports the geometric mean of runtimes over all instances.}
\label{table:srflp_met}
\end{table}

\subsection{Edge expansion problem}

For a given simple, undirected graph $G = (V,E)$ with $\lvert V \rvert = n \geq 3$, the edge expansion of $G$ is defined as the optimal value of
\begin{equation*}
    \min_{\emptyset \neq S \subsetneq V} \frac{\lvert \{ ij \in E \, \colon \, i \in S, \, j \in V \setminus S \} \rvert}{\min \{ \lvert S \rvert , \lvert V \setminus S \rvert \} }.
\end{equation*}
An SDP relaxation for this problem is given by
\begin{equation}\label{eq:edge_expansion_basic}
\begin{array}{rl}
\text{minimize} \quad & \langle L , X \rangle \\[1.2ex]
\text{subject to} \quad & \operatorname{trace}(X) = 1, \\
& \operatorname{diag}(X) = x, \\
& \frac{1}{\left\lfloor \frac{n}{2} \right\rfloor} \leq \rho \leq 1, \\
& \langle J_n, X \rangle \leq \left\lfloor \frac{n}{2} \right\rfloor, \\
& \begin{pmatrix}
    X & x \\
    x^\top & \rho \\
\end{pmatrix} \in \mathcal{S}_+^{n+1}, \quad \rho \in \mathbb{R}, \\
\end{array}
\end{equation}
and a strengthened relaxation is obtained by adding nonnegativity constraints on $x$ and $X$, resulting in a DNN relaxation; see~\cite{hrga2024connectivity}.

All instances of the edge expansion problem considered here were previously used in~\cite{gupte2026edge}. These include graphs of grlex and grevlex polytopes from~\cite{gupte2018ondantzig}, instances from the tenth DIMACS implementation challenge~\cite{DIMACS1992}, and network graphs from the online network repository~\cite{peixoto2020thenetzschleuder}.

Tables~\ref{table:edge:expansion_basic} and~\ref{table:edge_expansion_dnn} show that the Augmented Mixing Method successfully solves almost all instances of both the basic and the strengthened DNN relaxations for the edge expansion problem. Notably, the largest instance `celegans\_metabolic' with $n=454$ and more than~\num{100000} inequality constraints is solved in about \num{4552}~seconds. The primal-dual errors are typically below $10^{-10}$; two exceptions from the DNN relaxation are solved with reduced accuracy of about $10^{-6}$ due to time or iteration limits.

\begin{figure}[ht]
    \centering
    \includegraphics[width=\textwidth]{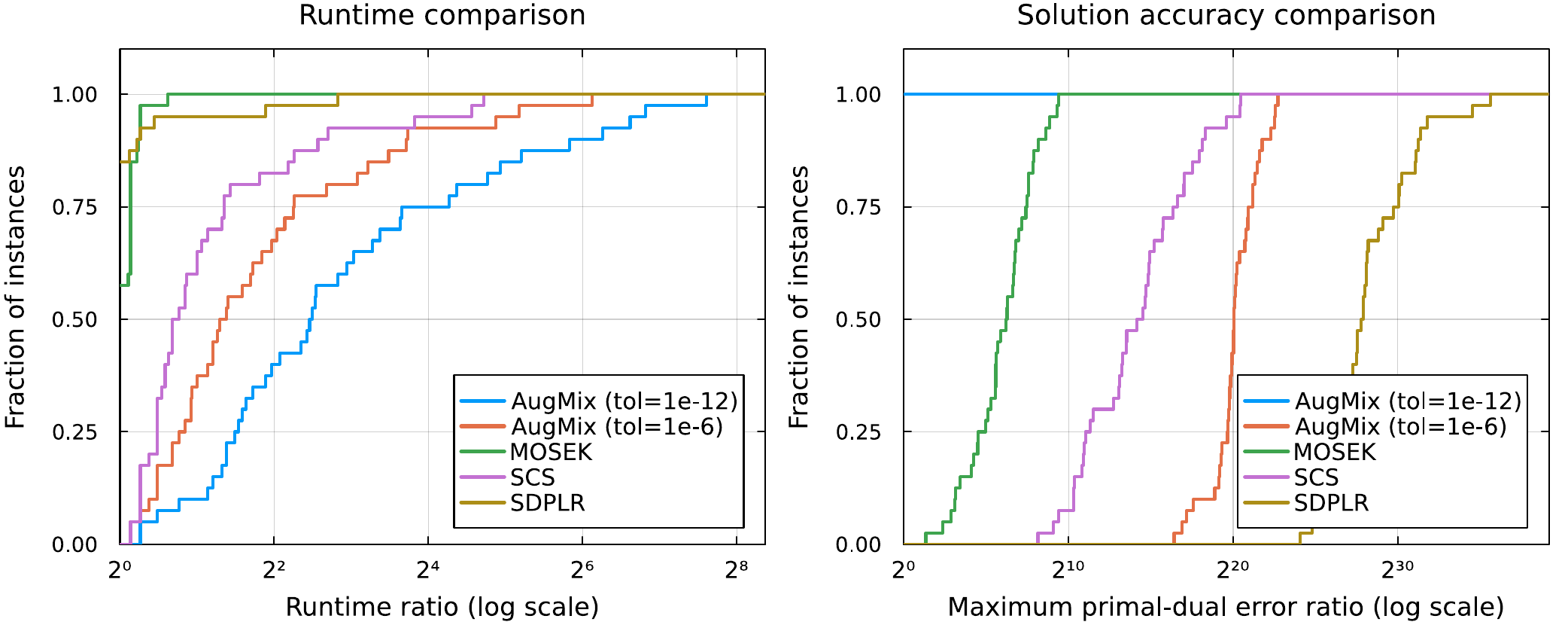}
    \caption{Performance profiles comparing the Augmented Mixing Method (with \texttt{tol=1e-12} and \texttt{tol=1e-6}) against \texttt{MOSEK}, \texttt{SCS}, and \texttt{SDPLR} on the basic SDP relaxation~\eqref{eq:edge_expansion_basic} of the edge expansion problem. Left: Runtime comparison. Right: Maximum primal-dual error comparison.}
    \label{fig:perfplot_edgeExpansionBasic}
\end{figure}

\begin{figure}[ht]
    \centering
    \includegraphics[width=\textwidth]{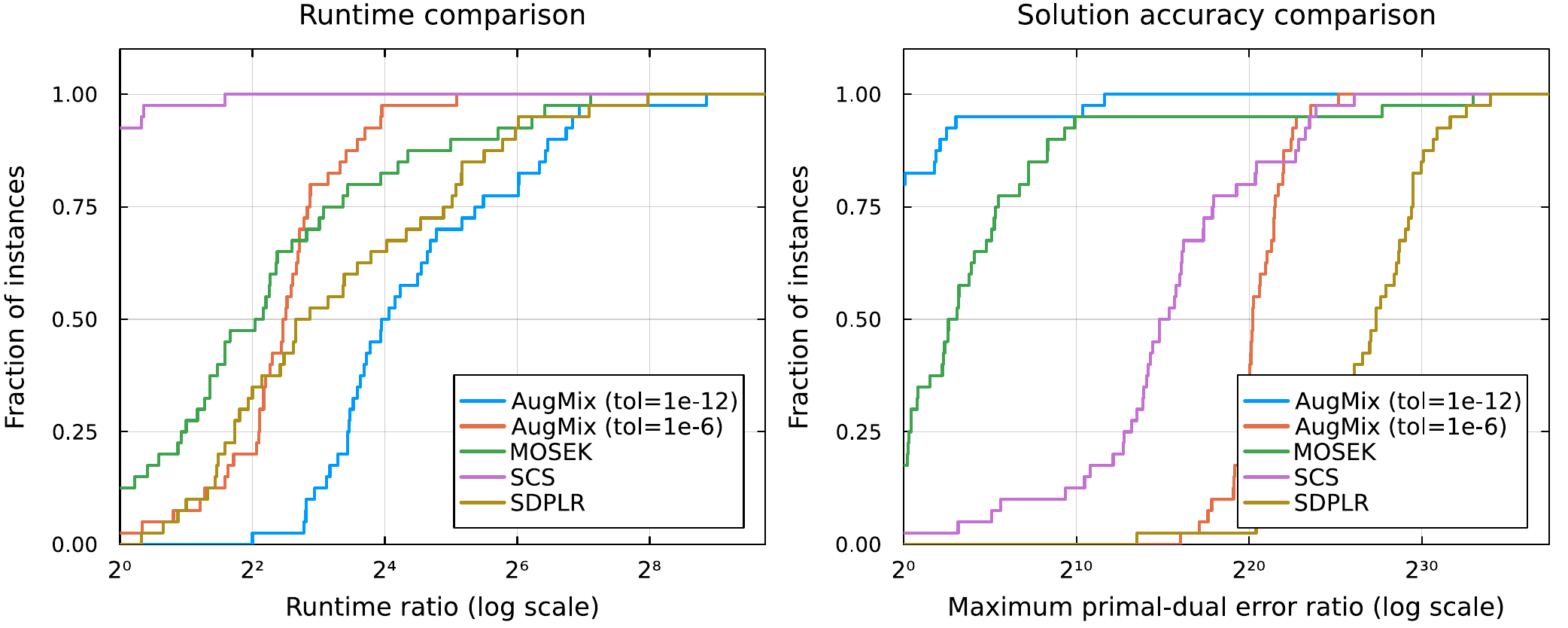}
    \caption{Performance profiles comparing the Augmented Mixing Method (with \texttt{tol=1e-12} and \texttt{tol=1e-6}) against \texttt{MOSEK}, \texttt{SCS}, and \texttt{SDPLR} on the DNN relaxation for bounding the edge expansion of a graph. Left: Runtime comparison. Right: Maximum primal-dual error comparison.}
    \label{fig:perfplot_edgeExpansionDnn}
\end{figure}

Figures~\ref{fig:perfplot_edgeExpansionBasic} and~\ref{fig:perfplot_edgeExpansionDnn} show performance profiles comparing the selected SDP solvers on both the basic SDP relaxation and the DNN relaxation. We note that \texttt{SDPLR} could not be run on the DNN relaxation of the instance `celegans\_metabolic' because the modeling language \texttt{JuMP} failed to pass the problem to \texttt{SDPLR}. For this instance, we conservatively set the runtime of \texttt{SDPLR} to the time limit of two hours and assigned all primal-dual errors a value of 1. Doing so does not influence the overall conclusions drawn from the performance profile in Figure~\ref{fig:perfplot_edgeExpansionDnn}.

The plots for the basic SDP relaxation~\eqref{eq:edge_expansion_basic} in Figure~\ref{fig:perfplot_edgeExpansionBasic} show that \texttt{MOSEK} is the fastest solver on most instances, closely followed by \texttt{SDPLR}. Both versions of the Augmented Mixing Method are the slowest solvers overall. However, when considering solution accuracy, the Augmented Mixing Method with \texttt{tol=1e-12} consistently produces the most accurate solutions, followed by \texttt{MOSEK}. The other solvers produce significantly less accurate solutions.

Regarding the DNN relaxation covered in Figure~\ref{fig:perfplot_edgeExpansionDnn}, \texttt{SCS} is the fastest solver overall, followed by \texttt{MOSEK} and the Augmented Mixing Method with \texttt{tol=1e-6}. However, the version with \texttt{tol=1e-12} produces the most accurate solutions on more than 75\% of all instances. Only \texttt{MOSEK} remains competitive in accuracy, but it was unable to solve the largest instance `celegans\_metabolic' due to the time limit.

\begin{table}[ht]
\centering
\resizebox{1.0\textwidth}{!}{
\setlength{\tabcolsep}{7pt}
\renewcommand{\arraystretch}{1.1}
\begin{filecontents*}{edgeExpansionBasic.csv}
instance,type,seed,statusCode,time,iter,tol,timeLimit,maxIters,itersZ,scaling,warmStart,k,n,numBlocks,avgBlockSize,minBlockSize,maxBlockSize,m,equations,inequalities,nnzA,densityAPercent,timeScaling,timeData,primalObj,dualObj,gap,pinf,dinf,compl,complNotPsdAbs,mu,evals,nnzY,gapOriginal,pinfOriginal,dinfOriginal,complOriginal
grevlex-8,edgeExpansionBasic,584644499,tol,1.6,2600,1.0e-12,7200.0,100000,50,true,false,10,38,1,38.0,38,38,41,38,3,1519,2.5657050199310856,0.0,8.416175842285156e-5,1.8522447099731312,1.8522447099409443,4.4e-13,2.6e-14,5.7e-15,4.6e-13,4.6e-13,8.0e+00,983341,40,6.8e-12,1.8e-14,2.7e-14,7.2e-12
grevlex-9,edgeExpansionBasic,224256750,tol,1.9,2250,1.0e-12,7200.0,100000,50,true,false,10,47,1,47.0,47,47,50,47,3,2302,2.0842009959257584,0.0,0.0001270771026611328,1.9777646304264718,1.9777646304573966,3.4e-13,2.5e-14,1.8e-14,6.1e-14,3.5e-13,7.1e+00,1066799,49,6.2e-12,1.5e-14,1.1e-13,1.1e-12
grevlex-10,edgeExpansionBasic,375645274,tol,8.7,8250,1.0e-12,7200.0,100000,50,true,false,11,57,1,57.0,57,57,60,57,3,3362,1.7246332204780959,0.0,0.0001621246337890625,2.156261836337601,2.1562618363258608,1.0e-13,1.8e-14,3.1e-14,3.3e-13,1.2e-13,8.2e+00,4720251,59,2.2e-12,9.7e-15,1.9e-13,6.9e-12
grevlex-11,edgeExpansionBasic,569158246,tol,11.4,8700,1.0e-12,7200.0,100000,50,true,false,12,68,1,68.0,68,68,71,68,3,4759,1.4495711291973292,0.0,0.00023508071899414062,2.3724659966626724,2.3724659967534314,6.7e-13,2.4e-14,1.1e-13,1.5e-13,6.6e-13,6.7e+00,6131075,70,1.6e-11,1.8e-14,6.9e-13,3.5e-12
grevlex-12,edgeExpansionBasic,366202531,tol,35.9,24100,1.0e-12,7200.0,100000,50,true,false,13,80,1,80.0,80,80,83,80,3,6559,1.2347515060240963,0.0,0.0002739429473876953,2.5522314860903217,2.5522314860277806,3.9e-13,2.3e-14,2.7e-14,6.4e-13,4.1e-13,8.4e+00,15980799,82,1.0e-11,1.1e-14,1.9e-13,1.7e-11
grevlex-13,edgeExpansionBasic,318165265,tol,21.7,9550,1.0e-12,7200.0,100000,50,true,false,14,93,1,93.0,93,93,96,93,3,8834,1.0639476625428759,0.0,0.0004210472106933594,2.704325690439811,2.7043256904351756,2.5e-14,2.0e-14,5.1e-14,5.1e-13,3.7e-14,6.6e+00,9509753,95,7.3e-13,8.7e-15,3.8e-13,1.5e-11
grevlex-14,edgeExpansionBasic,683481647,tol,97.1,47450,1.0e-12,7200.0,100000,50,true,false,15,107,1,107.0,107,107,110,107,3,11662,0.9260038590111086,0.0,0.00046706199645996094,2.889734749903278,2.889734750087006,8.5e-13,4.1e-14,3.3e-14,7.9e-13,8.7e-13,9.2e+00,31384729,109,2.7e-11,1.7e-14,3.0e-13,2.5e-11
grlex-8,edgeExpansionBasic,707290643,tol,0.2,350,1.0e-12,7200.0,100000,50,true,false,10,38,1,38.0,38,38,41,38,3,1519,2.5657050199310856,0.0,8.702278137207031e-5,0.47036416630440386,0.47036416635082057,6.1e-13,2.2e-13,4.4e-14,6.2e-18,6.5e-13,2.1e+00,143362,40,2.4e-11,2.4e-13,1.4e-13,2.9e-15
grlex-9,edgeExpansionBasic,376510247,tol,0.7,850,1.0e-12,7200.0,100000,50,true,false,10,47,1,47.0,47,47,50,47,3,2302,2.0842009959257584,0.0,0.000125885009765625,0.4335491530965476,0.4335491531090901,1.3e-13,1.3e-13,7.6e-15,5.2e-17,1.4e-13,2.5e+00,395430,49,6.7e-12,1.0e-13,2.5e-14,5.4e-15
grlex-10,edgeExpansionBasic,891081383,tol,1.2,1050,1.0e-12,7200.0,100000,50,true,false,11,57,1,57.0,57,57,60,57,3,3362,1.7246332204780959,0.0,0.0001628398895263672,0.41460679313133203,0.4146067931767048,3.7e-13,4.1e-13,2.1e-14,1.0e-17,3.3e-13,3.0e+00,597143,59,2.5e-11,2.2e-13,6.8e-14,6.6e-15
grlex-11,edgeExpansionBasic,577149841,tol,1.5,1150,1.0e-12,7200.0,100000,50,true,false,12,68,1,68.0,68,68,71,68,3,4759,1.4495711291973292,0.0,0.00022101402282714844,0.4054700562390484,0.4054700561881211,3.4e-13,2.9e-13,2.0e-15,3.7e-13,3.7e-13,3.5e+00,788067,70,2.8e-11,1.5e-13,1.5e-14,3.1e-11
grlex-12,edgeExpansionBasic,228141071,tol,2.0,1150,1.0e-12,7200.0,100000,50,true,false,13,80,1,80.0,80,80,83,80,3,6559,1.2347515060240963,0.0,0.00029778480529785156,0.39201522324840266,0.3920152233896529,7.8e-13,2.1e-14,5.0e-14,5.5e-16,7.8e-13,3.1e+00,929494,82,7.9e-11,9.7e-15,1.7e-13,6.0e-14
grlex-13,edgeExpansionBasic,265797155,tol,4.8,2150,1.0e-12,7200.0,100000,50,true,false,14,93,1,93.0,93,93,96,93,3,8834,1.0639476625428759,0.0,0.0003941059112548828,0.37678712258274183,0.376787122669002,4.0e-13,1.5e-13,5.7e-14,2.4e-13,4.1e-13,3.7e+00,2033733,95,4.9e-11,2.0e-13,1.9e-13,2.9e-11
grlex-14,edgeExpansionBasic,317678400,tol,4.5,1550,1.0e-12,7200.0,100000,50,true,false,15,107,1,107.0,107,107,110,107,3,11662,0.9260038590111086,0.0,0.0004570484161376953,0.36811044670487725,0.3681104466660297,1.5e-13,1.5e-13,1.1e-14,8.9e-13,1.4e-13,3.5e+00,1724819,109,2.2e-11,5.4e-13,4.8e-14,1.3e-10
grlex-15,edgeExpansionBasic,13464812,tol,6.1,1850,1.0e-12,7200.0,100000,50,true,false,16,122,1,122.0,122,122,125,122,3,15127,0.8130610051061543,0.0,0.0006101131439208984,0.36369831862195096,0.36369831870619446,2.8e-13,2.7e-13,1.8e-14,4.9e-13,2.7e-13,3.4e+00,2343755,124,4.9e-11,1.1e-12,6.3e-14,8.4e-11
grlex-16,edgeExpansionBasic,302061719,tol,11.6,2600,1.0e-12,7200.0,100000,50,true,false,17,138,1,138.0,138,138,141,138,3,19319,0.7194611657065907,0.0,0.0007519721984863281,0.3568967367156418,0.35689673662729515,2.6e-13,4.0e-13,5.4e-15,5.4e-13,2.8e-13,4.1e+00,3868212,140,5.2e-11,1.4e-13,2.1e-14,1.1e-10
grlex-17,edgeExpansionBasic,659944485,tol,18.3,3250,1.0e-12,7200.0,100000,50,true,false,18,155,1,155.0,155,155,158,155,3,24334,0.6410516471502522,0.0,0.0009450912475585938,0.3487105586048422,0.34871055865437467,1.3e-13,4.7e-14,2.0e-15,6.5e-14,1.2e-13,2.9e+00,5661783,157,2.9e-11,1.6e-14,1.3e-14,1.5e-11
adjnoun,edgeExpansionBasic,102870606,tol,6.0,2150,1.0e-12,7200.0,100000,50,true,false,16,113,1,113.0,113,113,116,113,3,12994,0.8772593106688883,0.0,0.0005350112915039062,0.36782876471675585,0.367828764621727,7.5e-13,6.9e-15,2.4e-14,7.9e-13,7.5e-13,1.5e+00,2392461,115,5.5e-11,8.7e-16,7.8e-14,5.8e-11
blumenau_drug,edgeExpansionBasic,636802559,tol,2.7,1750,1.0e-12,7200.0,100000,50,true,false,13,76,1,76.0,76,76,79,76,3,5927,1.298914758581998,0.0,0.00028705596923828125,0.1899395766679449,0.18993957666103822,1.0e-13,2.7e-15,2.0e-14,1.3e-13,1.0e-13,2.7e+00,1304860,78,5.0e-12,9.6e-16,6.7e-14,6.0e-12
celegans_metabolic,edgeExpansionBasic,40580129,tol,237.8,5500,1.0e-12,7200.0,100000,50,true,false,31,454,1,454.0,454,454,457,454,3,207023,0.21978127674106565,0.0,0.007581949234008789,0.1323256515677879,0.13232565138757788,3.9e-13,4.3e-15,6.8e-15,4.2e-13,3.9e-13,9.0e-01,29077153,456,1.4e-10,8.6e-16,1.9e-14,1.5e-10
celegansneural,edgeExpansionBasic,586253899,tol,389.3,23200,1.0e-12,7200.0,100000,50,true,false,25,298,1,298.0,298,298,301,298,3,89399,0.33445187662523357,0.0,0.003638029098510742,0.4321227100968949,0.4321227103950177,7.8e-13,1.9e-14,6.1e-14,5.7e-14,7.8e-13,2.5e+00,66992567,300,1.6e-10,4.6e-15,2.1e-13,1.2e-11
chesapeake,edgeExpansionBasic,833654923,tol,0.4,650,1.0e-12,7200.0,100000,50,true,false,10,40,1,40.0,40,40,43,40,3,1679,2.440406976744186,0.0,9.012222290039062e-5,0.8712397236403536,0.8712397236208479,2.5e-13,2.5e-14,0.0e+00,2.5e-13,2.5e-13,1.8e+00,268097,42,7.1e-12,1.6e-14,2.2e-15,7.3e-12
cintestinalis,edgeExpansionBasic,579052920,tol,134.2,19450,1.0e-12,7200.0,100000,50,true,false,21,206,1,206.0,206,206,209,206,3,42847,0.48310295357241595,0.0,0.003876924514770508,0.46186863751545115,0.4618686377444574,4.8e-13,4.4e-14,2.1e-14,2.6e-15,4.8e-13,1.7e+00,29958473,208,1.2e-10,1.8e-13,1.5e-13,6.7e-13
dolphins,edgeExpansionBasic,492059676,tol,3.2,2650,1.0e-12,7200.0,100000,50,true,false,12,63,1,63.0,63,63,66,63,3,4094,1.5628698168380708,0.0,0.00023102760314941406,0.08755910598407658,0.08755910601141181,4.9e-13,4.7e-15,1.6e-14,2.6e-14,4.9e-13,4.4e+00,1751915,65,2.3e-11,2.5e-15,8.0e-14,1.3e-12
foodweb_little_rock,edgeExpansionBasic,686181882,tol,32.6,5350,1.0e-12,7200.0,100000,50,true,false,20,184,1,184.0,184,184,187,184,3,34223,0.5405561648959292,0.0,0.0013039112091064453,0.49889721355951844,0.4988972132261619,6.3e-13,1.3e-15,1.8e-14,6.3e-13,6.3e-13,9.5e-01,9330610,186,1.7e-10,7.1e-17,9.6e-14,1.7e-10
football,edgeExpansionBasic,844759941,tol,2.9,950,1.0e-12,7200.0,100000,50,true,false,16,116,1,116.0,116,116,119,116,3,13687,0.8547622378321126,0.0,0.0005629062652587891,0.7388841746988184,0.7388841747100832,8.3e-14,1.1e-13,5.3e-14,6.0e-14,8.6e-14,1.1e+01,1119517,118,4.5e-12,4.6e-13,5.4e-13,3.3e-12
game_thrones,edgeExpansionBasic,742204460,tol,10.4,3750,1.0e-12,7200.0,100000,50,true,false,15,108,1,108.0,108,108,111,108,3,11879,0.9175070131860256,0.0,0.0004589557647705078,0.12225384242274805,0.12225384251465388,8.2e-13,1.7e-14,1.4e-14,5.4e-14,8.2e-13,1.7e+00,3986999,110,7.4e-11,7.0e-15,5.5e-14,4.9e-12
highschool,edgeExpansionBasic,105048298,tol,1.6,1050,1.0e-12,7200.0,100000,50,true,false,13,71,1,71.0,71,71,74,71,3,5182,1.3891495145214645,0.0,0.00026106834411621094,0.40182764094724027,0.4018276409139885,3.9e-13,1.4e-14,1.1e-14,4.7e-13,3.9e-13,5.7e+00,759295,73,1.8e-11,5.4e-14,5.4e-14,2.2e-11
jazz,edgeExpansionBasic,606046100,tol,9.6,1200,1.0e-12,7200.0,100000,50,true,false,21,199,1,199.0,199,199,202,199,3,39998,0.5000123759251004,0.0,0.0014190673828125,0.29927800375968966,0.29927800415518524,7.6e-13,3.5e-14,8.0e-15,5.7e-16,7.6e-13,2.5e+00,2362584,201,2.5e-10,1.4e-13,6.8e-14,2.1e-13
karate,edgeExpansionBasic,419007989,tol,0.2,450,1.0e-12,7200.0,100000,50,true,false,9,35,1,35.0,35,35,38,35,3,1294,2.7798066595059074,0.0,7.796287536621094e-5,0.24004040965970394,0.2400404096435704,3.8e-13,2.2e-15,2.2e-16,3.8e-13,3.8e-13,4.7e+00,159938,37,1.1e-11,9.7e-15,3.0e-15,1.1e-11
lesmis,edgeExpansionBasic,961839584,tol,4.6,2650,1.0e-12,7200.0,100000,50,true,false,13,78,1,78.0,78,78,81,78,3,6239,1.26602056801487,0.0,0.0005381107330322266,0.10545076255733905,0.10545076252056455,4.0e-13,9.7e-15,6.8e-15,4.4e-13,4.0e-13,2.1e+00,2202438,80,3.0e-11,4.5e-15,2.2e-14,3.3e-11
malaria_genes_HVR_1,edgeExpansionBasic,715556285,tol,113.7,5650,1.0e-12,7200.0,100000,50,true,false,25,308,1,308.0,308,308,311,308,3,95479,0.3236279630504377,0.0,0.0038390159606933594,0.10375119003251848,0.10375119012739752,2.2e-13,5.0e-14,5.6e-15,6.5e-13,2.2e-13,3.3e+00,20151609,310,7.9e-11,2.2e-13,7.3e-14,2.3e-10
moviegalaxies-52,edgeExpansionBasic,332653262,tol,1.8,1600,1.0e-12,7200.0,100000,50,true,false,12,60,1,60.0,60,60,63,60,3,3719,1.6397707231040564,0.0,0.00021719932556152344,0.18317120284974298,0.18317120287352398,3.8e-13,8.8e-15,1.2e-13,2.0e-13,3.8e-13,5.9e+00,1020838,62,1.7e-11,4.9e-15,3.2e-13,9.3e-12
moviegalaxies-567,edgeExpansionBasic,227033841,tol,2.1,2150,1.0e-12,7200.0,100000,50,true,false,11,53,1,53.0,53,53,56,53,3,2914,1.8524640187153536,0.0,0.0001518726348876953,0.16897153310324678,0.16897153307650856,4.1e-13,6.4e-15,3.4e-16,4.1e-13,4.1e-13,4.0e+00,1206508,55,2.0e-11,1.9e-14,4.3e-15,2.0e-11
polbooks,edgeExpansionBasic,294515227,tol,5.2,1800,1.0e-12,7200.0,100000,50,true,false,15,106,1,106.0,106,106,109,106,3,11447,0.9346595641140372,0.0,0.0004730224609375,0.1642911548994508,0.1642911549148306,1.3e-13,1.5e-13,1.8e-14,8.7e-13,1.2e-13,5.4e+00,2017161,108,1.2e-11,2.1e-13,1.2e-13,7.9e-11
revolution,edgeExpansionBasic,21858395,tol,31.6,6350,1.0e-12,7200.0,100000,50,true,false,18,142,1,142.0,142,142,145,142,3,20447,0.6993344232466191,0.0,0.0008320808410644531,0.038588654451966714,0.03858865443889787,1.3e-13,2.3e-15,3.6e-16,1.3e-13,1.3e-13,4.5e+00,10596676,144,1.2e-11,8.4e-16,3.9e-15,1.2e-11
sp_office,edgeExpansionBasic,238321863,tol,6.7,2950,1.0e-12,7200.0,100000,50,true,false,14,93,1,93.0,93,93,96,93,3,8834,1.0639476625428759,0.0,0.00039196014404296875,2.4994806768441937,2.4994806741357976,7.5e-13,3.4e-15,0.0e+00,7.5e-13,7.5e-13,5.0e+00,2975783,95,4.5e-10,1.5e-15,4.9e-15,4.5e-10
swingers,edgeExpansionBasic,781423324,tol,4.9,2100,1.0e-12,7200.0,100000,50,true,false,15,97,1,97.0,97,97,100,97,3,9602,1.0205122754809226,0.0,0.0004220008850097656,0.13674471836345928,0.13674471832837162,4.3e-13,4.7e-15,1.1e-14,4.5e-13,4.3e-13,1.8e+00,2025462,99,2.8e-11,2.0e-15,3.3e-14,2.9e-11
terrorists-911,edgeExpansionBasic,995370216,tol,1.3,1100,1.0e-12,7200.0,100000,50,true,false,12,63,1,63.0,63,63,66,63,3,4094,1.5628698168380708,0.0,0.00022220611572265625,0.09273504221797156,0.0927350421938498,4.0e-13,9.7e-13,3.0e-14,8.9e-13,4.6e-13,3.2e+00,647844,65,2.0e-11,5.1e-13,9.9e-14,4.5e-11
train_terrorists,edgeExpansionBasic,273949639,tol,2.3,1800,1.0e-12,7200.0,100000,50,true,false,12,65,1,65.0,65,65,68,65,3,4354,1.515489035851027,0.0,0.00020694732666015625,0.16813402689387105,0.16813402682781897,7.2e-13,1.0e-14,1.5e-16,7.2e-13,7.2e-13,2.5e+00,1245639,67,4.9e-11,5.3e-15,5.0e-15,4.9e-11
\end{filecontents*}
\pgfplotstabletypeset[
    col sep=comma,
    string type,
    columns={instance,n,equations,inequalities,statusCode,time,iter,gapOriginal,pinfOriginal,dinfOriginal,complOriginal},
    columns/instance/.style={column name=\textbf{Instance}, column type=l},
    columns/n/.style={column name=$n$, column type=r},
    columns/equations/.style={column name=$m_a$, column type=r},
    columns/inequalities/.style={column name=$m_b$, column type=r},
    columns/statusCode/.style={column name=\textbf{Status}, column type=c},
    columns/time/.style={column name=\textbf{Time~[s]}, column type=r},
    columns/iter/.style={column name=\textbf{Iter}, column type=r},
    columns/gapOriginal/.style={column name=\texttt{gap}, column type=c},
    columns/pinfOriginal/.style={column name=\texttt{pinf}, column type=c},
    columns/dinfOriginal/.style={column name=\texttt{dinf}, column type=c},
    columns/complOriginal/.style={column name=\texttt{compl}, column type=c},
    every head row/.style={before row=\toprule, after row=\midrule},
    every last row/.style={after row=\bottomrule}
]{edgeExpansionBasic.csv}
}
\caption{Results of the Augmented Mixing Method with \texttt{tol=1e-12} on selected instances of the basic SDP relaxation~\eqref{eq:edge_expansion_basic} for computing bounds on the edge expansion of a graph.}
\label{table:edge:expansion_basic}
\end{table}

\begin{table}[ht]
\centering
\resizebox{1.0\textwidth}{!}{
\setlength{\tabcolsep}{7pt}
\renewcommand{\arraystretch}{1.1}
\begin{filecontents*}{edgeExpansionDnn.csv}
instance,type,seed,statusCode,time,iter,tol,timeLimit,maxIters,itersZ,scaling,warmStart,k,n,numBlocks,avgBlockSize,minBlockSize,maxBlockSize,m,equations,inequalities,nnzA,densityAPercent,timeScaling,timeData,primalObj,dualObj,gap,pinf,dinf,compl,complNotPsdAbs,mu,evals,nnzY,gapOriginal,pinfOriginal,dinfOriginal,complOriginal
grevlex-8,edgeExpansionDnn,681068216,tol,2.5,2200,1.0e-12,7200.0,100000,50,true,false,38,38,1,38.0,38,38,744,38,706,2925,0.27226119202930926,0.0,0.0004100799560546875,1.8969724665228647,1.896972466465673,7.8e-13,1.5e-14,6.2e-15,7.9e-13,7.8e-13,2.3e+00,779771,53,1.2e-11,7.5e-15,2.9e-14,1.2e-11
grevlex-9,edgeExpansionDnn,338050489,tol,5.0,2800,1.0e-12,7200.0,100000,50,true,false,47,47,1,47.0,47,47,1131,47,1084,4464,0.1786758534233597,0.0,0.0006630420684814453,2.022724317615344,2.0227243175431937,7.8e-13,1.6e-14,8.2e-15,7.9e-13,7.8e-13,1.6e+00,1257855,63,1.4e-11,1.7e-14,4.5e-14,1.4e-11
grevlex-10,edgeExpansionDnn,73344184,tol,7.9,3000,1.0e-12,7200.0,100000,50,true,false,57,57,1,57.0,57,57,1656,57,1599,6554,0.12181377250227866,0.0,0.0011532306671142578,2.2053070663928143,2.2053070662913172,9.0e-13,3.3e-14,7.7e-15,9.0e-13,8.8e-13,1.3e+00,1653538,76,1.9e-11,1.8e-14,4.6e-14,1.9e-11
grevlex-11,edgeExpansionDnn,76487558,tol,11.5,3450,1.0e-12,7200.0,100000,50,true,false,68,68,1,68.0,68,68,2349,68,2281,9315,0.08575945591218231,0.0,0.0017168521881103516,2.428343942143872,2.428343942241098,7.2e-13,3.8e-14,9.3e-14,8.8e-15,7.2e-13,1.3e+00,2316279,93,1.7e-11,6.1e-14,6.0e-13,2.0e-13
grevlex-12,edgeExpansionDnn,410322911,tol,17.7,3850,1.0e-12,7200.0,100000,50,true,false,80,80,1,80.0,80,80,3243,80,3163,12879,0.062051919518963924,0.0,0.0026998519897460938,2.611926503609558,2.6119265037410373,8.2e-13,2.8e-14,9.0e-14,7.2e-15,8.2e-13,1.4e+00,3081603,107,2.1e-11,1.5e-14,6.2e-13,1.8e-13
grevlex-13,edgeExpansionDnn,75828636,tol,32.0,4600,1.0e-12,7200.0,100000,50,true,false,93,93,1,93.0,93,93,4374,93,4281,17390,0.04596792564858523,0.0,0.0041158199310302734,2.7655598104434804,2.7655598106118666,9.0e-13,1.2e-14,8.1e-14,4.4e-15,8.9e-13,1.5e+00,4404285,122,2.6e-11,4.8e-14,6.0e-13,1.3e-13
grevlex-14,edgeExpansionDnn,310923696,tol,51.0,5650,1.0e-12,7200.0,100000,50,true,false,107,107,1,107.0,107,107,5781,107,5674,23004,0.034756243738448295,0.0,0.005890846252441406,2.9549775784163907,2.954977578586627,7.9e-13,1.9e-14,6.0e-14,4.1e-15,7.8e-13,1.6e+00,6256677,140,2.5e-11,4.6e-14,4.8e-13,1.3e-13
grlex-8,edgeExpansionDnn,581256958,tol,1.2,1000,1.0e-12,7200.0,100000,50,true,false,38,38,1,38.0,38,38,744,38,706,2925,0.27226119202930926,0.0,0.00038313865661621094,0.4703641663018175,0.47036416629069294,1.5e-13,3.8e-15,0.0e+00,1.5e-13,1.5e-13,2.3e+00,385191,40,5.7e-12,2.6e-15,2.6e-15,5.8e-12
grlex-9,edgeExpansionDnn,301226142,tol,1.8,1000,1.0e-12,7200.0,100000,50,true,false,47,47,1,47.0,47,47,1131,47,1084,4464,0.1786758534233597,0.0,0.0006608963012695312,0.4335491530949363,0.4335491530734974,2.2e-13,3.6e-15,0.0e+00,2.2e-13,2.2e-13,2.0e+00,485100,49,1.1e-11,2.3e-15,4.7e-15,1.2e-11
grlex-10,edgeExpansionDnn,126546846,tol,4.6,1700,1.0e-12,7200.0,100000,50,true,false,57,57,1,57.0,57,57,1656,57,1599,6554,0.12181377250227866,0.0,0.001123189926147461,0.4146067931363388,0.4146067930927457,3.6e-13,1.8e-14,0.0e+00,3.6e-13,3.6e-13,2.0e+00,1001972,59,2.4e-11,1.0e-14,3.5e-15,2.4e-11
grlex-11,edgeExpansionDnn,638220867,tol,5.2,1500,1.0e-12,7200.0,100000,50,true,false,68,68,1,68.0,68,68,2349,68,2281,9315,0.08575945591218231,0.0,0.0017249584197998047,0.4054700562426312,0.4054700561865046,3.8e-13,3.4e-14,0.0e+00,3.8e-13,3.8e-13,1.8e+00,1093529,70,3.1e-11,1.7e-14,8.5e-15,3.1e-11
grlex-12,edgeExpansionDnn,778703852,tol,7.7,1600,1.0e-12,7200.0,100000,50,true,false,80,80,1,80.0,80,80,3243,80,3163,12879,0.062051919518963924,0.0,0.0027179718017578125,0.39201522324805843,0.3920152230752625,9.6e-13,4.3e-14,1.7e-15,9.6e-13,9.6e-13,2.6e+00,1393631,82,9.7e-11,2.1e-14,5.8e-15,9.7e-11
grlex-13,edgeExpansionDnn,929010428,tol,16.6,2300,1.0e-12,7200.0,100000,50,true,false,93,93,1,93.0,93,93,4374,93,4281,17390,0.04596792564858523,0.0,0.004099845886230469,0.3767871225804945,0.37678712253206,2.2e-13,3.6e-15,0.0e+00,2.2e-13,2.2e-13,1.7e+00,2340211,95,2.8e-11,1.6e-15,1.3e-14,2.8e-11
grlex-14,edgeExpansionDnn,798487695,tol,23.4,2500,1.0e-12,7200.0,100000,50,true,false,107,107,1,107.0,107,107,5781,107,5674,23004,0.034756243738448295,0.0,0.005883932113647461,0.368110446702852,0.3681104465087175,7.6e-13,1.2e-14,3.8e-15,7.7e-13,7.6e-13,1.8e+00,2969556,109,1.1e-10,4.8e-15,1.3e-14,1.1e-10
grlex-15,edgeExpansionDnn,975125365,tol,31.9,2550,1.0e-12,7200.0,100000,50,true,false,122,122,1,122.0,122,122,7506,122,7384,29889,0.026753657541582966,0.0,0.010777950286865234,0.36369831862498864,0.36369831840173517,7.5e-13,8.2e-15,1.4e-14,7.7e-13,7.5e-13,2.3e+00,3469501,124,1.3e-10,3.1e-15,4.7e-14,1.3e-10
grlex-16,edgeExpansionDnn,977283631,tol,47.7,2950,1.0e-12,7200.0,100000,50,true,false,138,138,1,138.0,138,138,9594,138,9456,38225,0.02092134528699915,0.0,0.0180208683013916,0.35689673672212846,0.35689673661381244,3.2e-13,2.8e-14,2.3e-15,3.2e-13,3.2e-13,1.4e+00,4737050,140,6.3e-11,5.2e-14,8.2e-15,6.4e-11
grlex-17,edgeExpansionDnn,774141164,tol,77.1,3450,1.0e-12,7200.0,100000,50,true,false,155,155,1,155.0,155,155,12093,155,11938,48204,0.016591499128373214,0.0,0.02851390838623047,0.34871055860542866,0.3487105584541644,3.8e-13,1.3e-14,4.0e-15,3.9e-13,3.8e-13,1.4e+00,6342905,157,8.9e-11,3.4e-14,1.4e-14,9.0e-11
adjnoun,edgeExpansionDnn,504263236,tol,84.4,10200,1.0e-12,7200.0,100000,50,true,false,113,113,1,113.0,113,113,6444,113,6331,25650,0.03117273809518601,0.0,0.007311820983886719,0.3707890460874795,0.3707890460567417,2.4e-13,3.5e-13,6.5e-13,9.7e-13,2.5e-13,7.2e-01,10195954,122,1.8e-11,2.7e-14,2.1e-12,7.1e-11
blumenau_drug,edgeExpansionDnn,216394955,tol,6.2,1650,1.0e-12,7200.0,100000,50,true,false,76,76,1,76.0,76,76,2929,76,2853,11627,0.0687260076661979,0.0,0.0022020339965820312,0.18993957666799502,0.18993957662066066,7.1e-13,1.5e-14,1.6e-13,8.6e-13,7.1e-13,2.8e+00,1257075,78,3.4e-11,1.0e-15,5.3e-13,4.2e-11
celegans_metabolic,edgeExpansionDnn,521781560,tol,4552.3,17050,1.0e-12,7200.0,100000,50,true,false,454,454,1,454.0,454,454,103288,454,102834,412685,0.00193846118768382,0.5,1.1363399028778076,0.13280342646860588,0.13280342662657388,3.5e-13,7.4e-13,4.9e-14,1.6e-13,3.4e-13,1.8e+00,73784603,918,1.2e-10,8.6e-14,1.3e-13,5.8e-11
celegansneural,edgeExpansionDnn,284075511,time,7200.0,93657,1.0e-12,7200.0,100000,50,true,false,298,298,1,298.0,298,298,44554,298,44256,177905,0.0044964412711405635,0.0,0.15144801139831543,0.43298388047040776,0.4329862967003864,6.3e-09,1.1e-12,Inf,Inf,6.3e-09,2.1e+01,207663645,3457,1.3e-06,2.5e-12,2.5e-09,3.0e-09
chesapeake,edgeExpansionDnn,544222124,tol,0.8,600,1.0e-12,7200.0,100000,50,true,false,40,40,1,40.0,40,40,823,40,783,3239,0.24597509113001215,0.0,0.0004680156707763672,0.8712397236407498,0.8712397237173907,9.6e-13,3.2e-14,5.0e-14,7.8e-14,9.6e-13,2.5e+00,251626,42,2.8e-11,1.3e-13,1.5e-13,2.3e-12
cintestinalis,edgeExpansionDnn,599961352,tol,518.9,15750,1.0e-12,7200.0,100000,50,true,false,206,206,1,206.0,206,206,21324,206,21118,85077,0.009401757662888343,0.1,0.07901191711425781,0.461868637515318,0.4618686373918645,2.6e-13,8.4e-15,2.1e-15,2.6e-13,2.6e-13,1.6e+00,28955932,208,6.4e-11,2.5e-15,6.7e-14,6.5e-11
dolphins,edgeExpansionDnn,182145680,tol,79.7,25200,1.0e-12,7200.0,100000,50,true,false,63,63,1,63.0,63,63,2019,63,1956,8000,0.09983264305300202,0.0,0.0014541149139404297,0.08755910598397346,0.08755910602925758,8.0e-13,3.9e-14,1.4e-13,9.9e-13,8.0e-13,1.6e+00,15660079,66,3.9e-11,2.0e-14,6.6e-13,4.7e-11
foodweb_little_rock,edgeExpansionDnn,469145144,tol,395.3,19050,1.0e-12,7200.0,100000,50,true,false,184,184,1,184.0,184,184,17023,184,16839,67895,0.011780563036754343,0.0,0.041359901428222656,0.4988972135596007,0.4988972130357946,9.9e-13,7.6e-15,0.0e+00,9.9e-13,9.9e-13,3.9e+00,25296250,186,2.6e-10,2.1e-16,1.0e-14,2.6e-10
football,edgeExpansionDnn,677830017,tol,63.8,7800,1.0e-12,7200.0,100000,50,true,false,116,116,1,116.0,116,116,6789,116,6673,27027,0.029585305249153655,0.0,0.007961034774780273,0.757077493081755,0.7570774931697146,6.5e-13,7.9e-15,2.6e-14,2.7e-15,6.5e-13,7.6e+00,6685923,2458,3.5e-11,3.1e-15,2.6e-13,1.7e-13
game_thrones,edgeExpansionDnn,758074638,tol,82.3,10750,1.0e-12,7200.0,100000,50,true,false,108,108,1,108.0,108,108,5889,108,5781,23435,0.03411739727249497,0.0,0.006398916244506836,0.12309842023795554,0.12309842023669534,1.1e-14,9.1e-13,2.4e-13,5.8e-13,8.4e-15,5.1e-01,10182822,399,1.0e-12,4.4e-14,8.7e-13,5.2e-11
highschool,edgeExpansionDnn,624847530,tol,37.0,10600,1.0e-12,7200.0,100000,50,true,false,71,71,1,71.0,71,71,2559,71,2488,10152,0.07869816856989567,0.0,0.0018658638000488281,0.40523336003339966,0.40523336010339306,8.2e-13,2.3e-14,6.7e-14,5.0e-14,8.2e-13,1.0e+01,6111549,341,3.9e-11,9.2e-14,3.2e-13,2.4e-12
jazz,edgeExpansionDnn,756200486,tol,619.8,20900,1.0e-12,7200.0,100000,50,true,false,199,199,1,199.0,199,199,19903,199,19704,79400,0.010073857577955897,0.1,0.05311703681945801,0.2994629188872777,0.2994629191391877,4.8e-13,8.9e-13,2.1e-13,5.3e-13,4.8e-13,5.8e-01,37474505,202,1.6e-10,2.3e-14,1.1e-12,1.7e-10
karate,edgeExpansionDnn,857422156,tol,6.5,7600,1.0e-12,7200.0,100000,50,true,false,35,35,1,35.0,35,35,633,35,598,2484,0.3203404584582648,0.0,0.0003190040588378906,0.24020930447765876,0.24020930448365377,1.4e-13,9.6e-13,1.8e-13,2.3e-13,1.5e-13,5.2e-01,1958443,64,4.0e-12,3.0e-13,5.3e-13,6.7e-12
lesmis,edgeExpansionDnn,172000776,iter,464.3,100000,1.0e-12,7200.0,100000,50,true,false,78,78,1,78.0,78,78,3084,78,3006,12245,0.06526122397119105,0.0,0.0024459362030029297,0.10545077198040291,0.10545282897948875,2.2e-08,1.4e-10,Inf,Inf,2.2e-08,2.9e+00,79467506,87,1.7e-06,6.5e-11,5.1e-09,4.4e-07
malaria_genes_HVR_1,edgeExpansionDnn,135556971,tol,4552.0,49750,1.0e-12,7200.0,100000,50,true,false,308,308,1,308.0,308,308,47589,308,47281,190035,0.004209452209208237,0.0,0.41091299057006836,0.11026171115017612,0.11026171150661013,8.3e-13,9.4e-13,9.9e-13,6.6e-13,8.3e-13,5.4e-01,150368402,3576,2.9e-10,6.8e-13,7.7e-12,2.3e-10
moviegalaxies-52,edgeExpansionDnn,717200605,tol,14.1,5600,1.0e-12,7200.0,100000,50,true,false,60,60,1,60.0,60,60,1833,60,1773,7259,0.1100048493665515,0.0,0.0012290477752685547,0.1846950142484011,0.18469501425827098,1.6e-13,7.5e-15,1.6e-14,1.5e-14,1.6e-13,6.5e+00,3002534,793,7.2e-12,3.2e-14,4.2e-14,6.9e-13
moviegalaxies-567,edgeExpansionDnn,414755637,tol,74.7,26550,1.0e-12,7200.0,100000,50,true,false,53,53,1,53.0,53,53,1434,53,1381,5670,0.1407609432323777,0.0,0.0009319782257080078,0.1800125222262946,0.1800125222464909,3.1e-13,2.7e-13,1.9e-13,6.7e-13,3.1e-13,2.4e+02,20097338,77,1.5e-11,1.5e-14,5.1e-13,3.2e-11
polbooks,edgeExpansionDnn,917979893,tol,92.4,12100,1.0e-12,7200.0,100000,50,true,false,106,106,1,106.0,106,106,5674,106,5568,22577,0.03541319990518416,0.0,0.005876064300537109,0.182829096475146,0.18282909648639814,9.4e-14,7.9e-13,2.8e-13,4.8e-13,9.2e-14,5.5e-01,11937478,1616,8.2e-12,2.2e-13,1.4e-12,4.2e-11
revolution,edgeExpansionDnn,96094795,tol,136.7,7350,1.0e-12,7200.0,100000,50,true,false,142,142,1,142.0,142,142,10156,142,10014,40469,0.01976164494078841,0.0,0.019397974014282227,0.038588654451936884,0.0385886544694425,1.8e-13,1.8e-14,2.0e-15,3.3e-19,1.8e-13,3.0e+00,14159510,144,1.6e-11,6.5e-15,8.2e-15,1.2e-15
sp_office,edgeExpansionDnn,310650509,tol,21.4,3200,1.0e-12,7200.0,100000,50,true,false,93,93,1,93.0,93,93,4374,93,4281,17390,0.04596792564858523,0.0,0.004199028015136719,2.499480676844131,2.4994806762250663,1.7e-13,2.0e-15,0.0e+00,1.7e-13,1.7e-13,3.2e+00,3151132,95,1.0e-10,8.8e-16,5.6e-15,1.0e-10
swingers,edgeExpansionDnn,969151401,tol,34.6,5450,1.0e-12,7200.0,100000,50,true,false,97,97,1,97.0,97,97,4756,97,4659,18914,0.04226667361502118,0.0,0.004759073257446289,0.13739902000940407,0.137399020010279,1.1e-14,5.2e-13,5.8e-15,1.9e-14,1.7e-14,5.9e-01,4852963,102,6.8e-13,2.8e-14,1.8e-14,1.2e-12
terrorists-911,edgeExpansionDnn,374926724,tol,70.4,29350,1.0e-12,7200.0,100000,50,true,false,63,63,1,63.0,63,63,2019,63,1956,8000,0.09983264305300202,0.0,0.0014090538024902344,0.09279272676294743,0.0927927267575527,9.1e-14,6.9e-13,9.4e-14,3.4e-13,8.8e-14,5.2e-01,13366594,86,4.6e-12,5.6e-14,3.1e-13,1.7e-11
train_terrorists,edgeExpansionDnn,890715037,tol,73.6,27950,1.0e-12,7200.0,100000,50,true,false,65,65,1,65.0,65,65,2148,65,2083,8514,0.09381508049320683,0.0,0.0015521049499511719,0.16814655822906752,0.16814655822095542,8.8e-14,8.9e-13,7.3e-14,2.9e-13,8.8e-14,5.0e-01,14609628,104,6.1e-12,7.0e-14,2.7e-13,2.0e-11
\end{filecontents*}
\pgfplotstabletypeset[
    col sep=comma,
    string type,
    columns={instance,n,equations,inequalities,statusCode,time,iter,gapOriginal,pinfOriginal,dinfOriginal,complOriginal},
    columns/instance/.style={column name=\textbf{Instance}, column type=l},
    columns/n/.style={column name=$n$, column type=r},
    columns/equations/.style={column name=$m_a$, column type=r},
    columns/inequalities/.style={column name=$m_b$, column type=r},
    columns/statusCode/.style={column name=\textbf{Status}, column type=c},
    columns/time/.style={column name=\textbf{Time~[s]}, column type=r},
    columns/iter/.style={column name=\textbf{Iter}, column type=r},
    columns/gapOriginal/.style={column name=\texttt{gap}, column type=c},
    columns/pinfOriginal/.style={column name=\texttt{pinf}, column type=c},
    columns/dinfOriginal/.style={column name=\texttt{dinf}, column type=c},
    columns/complOriginal/.style={column name=\texttt{compl}, column type=c},
    every head row/.style={before row=\toprule, after row=\midrule},
    every last row/.style={after row=\bottomrule}
]{edgeExpansionDnn.csv}
}
\caption{Results of the Augmented Mixing Method with \texttt{tol=1e-12} on selected instances of the DNN relaxation for bounding the edge expansion of a graph.}
\label{table:edge_expansion_dnn}
\end{table}

\subsection{Randomly generated SDPs}\label{sec:rand}

We also test the Augmented Mixing Method on randomly generated SDPs, including instances with multiple positive semidefinite matrix variables. For simplicity, we only consider SDPs with equality constraints. We generate these SDPs following the approach in~\cite{yamashita2003sdpa}, with an additional parameter $d \in (0,1]$ to specify the density of the data matrices.

For one-block SDPs, we first fix the order of the matrix variable $n$ and the number of equality constraints $m_a$. We set $A_1$ to be the $n \times n$ identity matrix $I_n$. For the remaining matrices $C, A_2, \ldots, A_{m_a}$, a fraction $d$ of the entries is selected uniformly at random and then filled with values drawn from the uniform distribution $[-1,1]$. Finally, we set the entries of the right-hand side vector $a$ to $a_i = \operatorname{trace}(A_i)$ for $i = 1,\ldots,m_a$. By construction, the primal and dual SDPs always satisfy Slater's condition as assumed throughout the paper. Multi-block SDPs with block sizes $n_1,\ldots,n_q$ are constructed analogously.

We generate three different types of SDPs with one or more matrix variables and different densities: `rand\_$n$\_$m_a$\_$d$', `rand\_$q$x$n$\_$m_a$\_$d$', and `rand\_$p$\_$q$x$n$\_$m_a$\_$d$'. Here, $m_a$ always determines the number of equality constraints and $d$ determines the density of data matrices. SDPs of type `rand\_$n$\_$m_a$\_$d$' are one-block SDPs with a matrix variable of order $n$, SDPs of type `rand\_$q$x$n$\_$m_a$\_$d$' have $q$ blocks of equal size $n$, and SDPs of type `rand\_$p$\_$q$x$n$\_$m_a$\_$d$' have one block of size $p$ and $q$ blocks of size $n$.

Table~\ref{table:rand} shows the results of the Augmented Mixing Method with \texttt{tol=1e-12} on a wide range of randomly generated SDPs. The method consistently solves the instances with primal-dual errors below $10^{-11}$, although the time limit is reached on two instances. In general, the runtime and number of outer iterations increase with problem size, but the method consistently delivers reliable results across all tested SDPs.

\begin{figure}[ht]
    \centering
    \includegraphics[width=\textwidth]{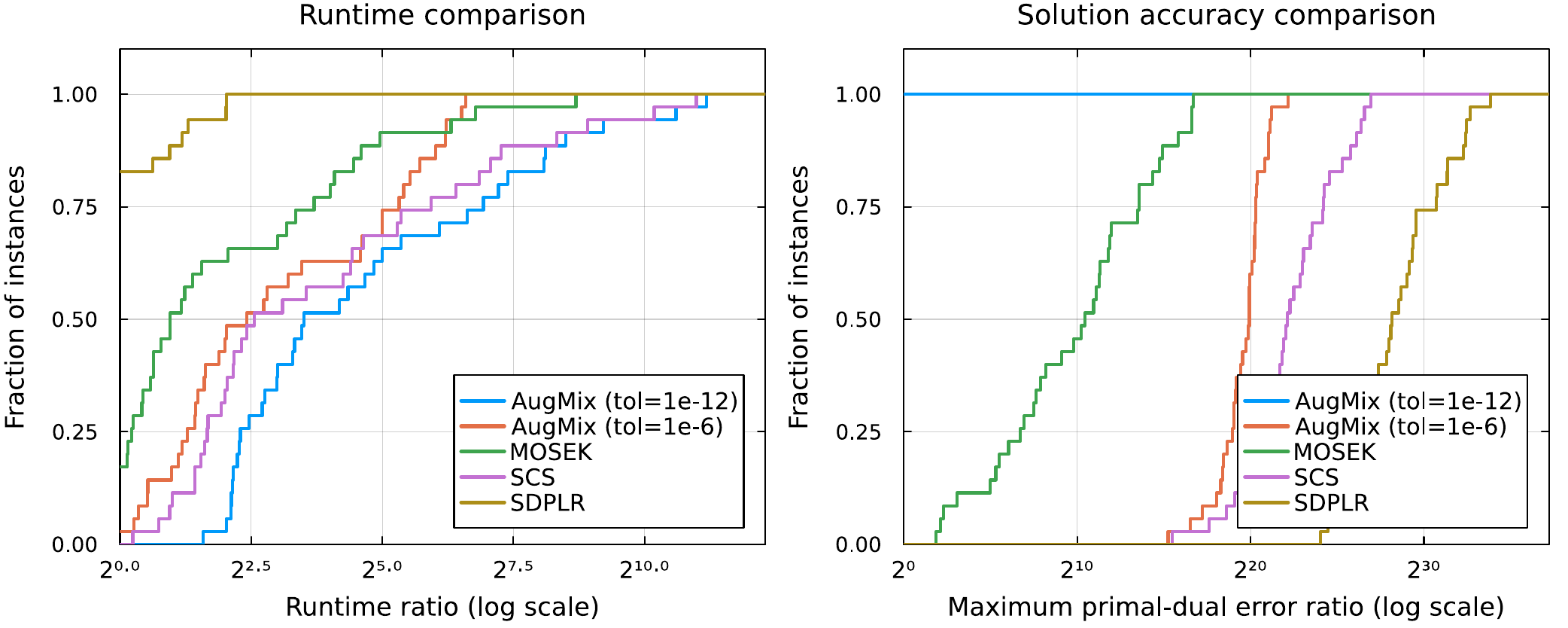}
    \caption{Performance profiles comparing the Augmented Mixing Method (with \texttt{tol=1e-12} and \texttt{tol=1e-6}) against \texttt{MOSEK}, \texttt{SCS}, and \texttt{SDPLR} on a wide range of randomly generated SDPs. Left: Runtime comparison. Right: Maximum primal-dual error comparison.}
    \label{fig:perfplot_rand}
\end{figure}

Figure~\ref{fig:perfplot_rand} reveals that \texttt{SDPLR} is the fastest of all solvers on more than 75\% of the instances tested. However, it also produces the least accurate solutions by a wide margin. In terms of speed, \texttt{MOSEK} is the second-fastest solver and the Augmented Mixing Method with \texttt{tol=1e-6} is third. These are followed by \texttt{SCS}, and then by the Augmented Mixing Method with \texttt{tol=1e-12}. The latter produces the most accurate solutions on all instances considered. \texttt{MOSEK} ranks second in terms of solution quality, followed by the Augmented Mixing Method with \texttt{tol=1e-6} and \texttt{SCS}.

\begin{table}[ht]
\centering
\resizebox{1.0\textwidth}{!}{
\setlength{\tabcolsep}{7pt}
\renewcommand{\arraystretch}{1.1}
\begin{filecontents*}{rand.csv}
instance,type,seed,statusCode,time,iter,tol,timeLimit,maxIters,itersZ,scaling,warmStart,k,n,numBlocks,avgBlockSize,minBlockSize,maxBlockSize,m,equations,inequalities,nnzA,densityAPercent,timeScaling,timeData,primalObj,dualObj,gap,pinf,dinf,compl,complNotPsdAbs,mu,evals,nnzY,gapOriginal,pinfOriginal,dinfOriginal,complOriginal
rand_50_100_1.0,rand,372791249,tol,0.8,350,1.0e-12,7200.0,100000,50,true,false,15,50,1,50.0,50,50,100,100,0,247550,99.02,0.0,0.0013298988342285156,-287.7703862754443,-287.7703862754088,4.5e-14,6.4e-14,1.6e-13,2.4e-13,2.8e-14,2.4e+00,228267,100,6.1e-14,5.2e-13,2.3e-12,3.2e-13
rand_50_200_1.0,rand,61473150,tol,3.3,800,1.0e-12,7200.0,100000,50,true,false,20,50,1,50.0,50,50,200,200,0,497550,99.51,0.0,0.009387969970703125,-261.166090756184,-261.16609075608534,1.3e-13,1.1e-13,1.5e-13,3.2e-13,6.6e-14,2.2e+00,594466,200,1.9e-13,9.1e-13,2.4e-12,4.6e-13
rand_50_400_1.0,rand,583750979,tol,9.1,1100,1.0e-12,7200.0,100000,50,true,false,29,50,1,50.0,50,50,400,400,0,997550,99.755,0.0,0.006529092788696289,-193.50301676736987,-193.50301676724334,2.1e-13,2.0e-13,2.3e-13,4.3e-13,2.0e-13,1.9e+00,969124,400,3.3e-13,1.7e-12,3.4e-12,6.7e-13
rand_50_800_1.0,rand,662188238,tol,5.7,300,1.0e-12,7200.0,100000,50,true,false,40,50,1,50.0,50,50,800,800,0,1997550,99.8775,0.0,0.03650307655334473,-108.37697429035077,-108.37697429034921,3.5e-15,4.3e-15,1.2e-14,2.6e-14,3.7e-15,3.2e+00,310703,800,7.4e-15,3.6e-14,1.9e-13,5.5e-14
rand_100_200_1.0,rand,190913164,tol,5.5,500,1.0e-12,7200.0,100000,50,true,false,20,100,1,100.0,100,100,200,200,0,1990100,99.505,0.0,0.046556949615478516,-908.8239864291731,-908.8239864288636,1.3e-13,1.1e-13,1.6e-13,4.9e-13,8.6e-14,3.0e+00,618068,200,1.7e-13,1.3e-12,4.6e-12,6.5e-13
rand_100_400_1.0,rand,596971048,tol,10.8,450,1.0e-12,7200.0,100000,50,true,false,29,100,1,100.0,100,100,400,400,0,3990100,99.7525,0.0,0.23722505569458008,-860.05844649828,-860.0584464981631,5.1e-14,1.3e-13,1.5e-14,4.7e-14,6.4e-15,8.6e-01,638088,400,6.8e-14,1.5e-12,4.3e-13,6.2e-14
rand_100_800_0.001,rand,368254007,tol,6.4,1400,1.0e-12,7200.0,100000,50,true,false,40,100,1,100.0,100,100,800,800,0,8019,0.1002375,0.0,0.0015301704406738281,-44.745168023541446,-44.74516802351852,2.1e-13,7.6e-13,1.3e-13,2.9e-13,7.0e-15,2.3e+00,1477072,800,2.5e-13,2.2e-13,1.7e-13,3.4e-13
rand_100_800_1.0,rand,919640344,tol,45.2,950,1.0e-12,7200.0,100000,50,true,false,40,100,1,100.0,100,100,800,800,0,7990100,99.87625,0.0,0.2734990119934082,-705.4966497895135,-705.4966497888029,3.5e-13,1.7e-13,1.5e-13,5.3e-13,2.6e-13,2.3e+00,1480133,800,5.0e-13,2.0e-12,4.6e-12,7.5e-13
rand_100_1600_0.001,rand,320726000,tol,4.1,400,1.0e-12,7200.0,100000,50,true,false,57,100,1,100.0,100,100,1600,1600,0,15935,0.09959375,0.0,0.015582799911499023,-25.25669448229579,-25.25669448229945,4.9e-14,1.2e-13,1.6e-13,7.4e-14,1.4e-15,1.7e+00,756656,1600,6.8e-14,5.5e-14,2.6e-13,1.1e-13
rand_100_1600_1.0,rand,86330892,tol,75.7,750,1.0e-12,7200.0,100000,50,true,false,57,100,1,100.0,100,100,1600,1600,0,15990100,99.938125,0.2,0.33797597885131836,-547.6054296375016,-547.6054296377561,1.5e-13,9.2e-13,2.0e-13,8.7e-13,2.7e-14,6.0e-01,1140308,1600,2.3e-13,1.1e-11,5.8e-12,1.4e-12
rand_100_3200_0.001,rand,464616117,tol,6.8,450,1.0e-12,7200.0,100000,50,true,false,80,100,1,100.0,100,100,3200,3200,0,31770,0.09928125,0.0,0.018033981323242188,-5.7035694727104715,-5.70356947270698,1.3e-13,1.8e-13,4.7e-14,2.4e-13,2.3e-15,1.6e+00,920379,3200,2.8e-13,8.5e-14,5.9e-14,5.4e-13
rand_100_3200_1.0,rand,823100754,tol,156.7,850,1.0e-12,7200.0,100000,50,true,false,80,100,1,100.0,100,100,3200,3200,0,31990100,99.9690625,0.1,0.6419739723205566,-292.792174842191,-292.7921748420395,1.2e-13,3.0e-13,1.9e-13,7.9e-13,2.3e-14,6.8e-01,1348597,3200,2.6e-13,3.7e-12,5.6e-12,1.7e-12
rand_200_400_1.0,rand,312434,tol,43.5,600,1.0e-12,7200.0,100000,50,true,false,29,200,1,200.0,200,200,400,400,0,15960200,99.75125,0.2,0.37706494331359863,-2768.166215037435,-2768.1662150393886,2.7e-13,2.3e-13,8.2e-14,9.1e-13,2.9e-13,2.3e+00,1269771,400,3.5e-13,3.8e-12,4.7e-12,1.2e-12
rand_200_800_0.001,rand,944356629,tol,10.2,1100,1.0e-12,7200.0,100000,50,true,false,40,200,1,200.0,200,200,800,800,0,31962,0.09988125,0.0,0.006971120834350586,-177.22248318656386,-177.22248318655681,1.7e-14,9.5e-13,3.6e-14,1.5e-13,4.5e-15,6.7e-01,1283659,800,2.0e-14,5.1e-13,8.0e-14,1.7e-13
rand_200_800_1.0,rand,978601445,tol,168.4,1200,1.0e-12,7200.0,100000,50,true,false,40,200,1,200.0,200,200,800,800,0,31960200,99.875625,0.1,0.46750688552856445,-2540.997506923309,-2540.9975069223174,1.5e-13,1.2e-13,1.9e-13,6.2e-13,1.4e-13,2.8e+00,2899239,800,1.9e-13,2.0e-12,1.1e-11,8.2e-13
rand_200_1600_0.001,rand,337628751,tol,9.8,550,1.0e-12,7200.0,100000,50,true,false,57,200,1,200.0,200,200,1600,1600,0,63830,0.099734375,0.0,0.006855964660644531,-123.30178212594343,-123.30178212594855,1.7e-14,8.9e-13,4.3e-14,1.5e-13,1.2e-14,6.0e-01,952740,1600,2.0e-14,4.6e-13,9.6e-14,1.8e-13
rand_200_1600_1.0,rand,294144527,tol,199.5,800,1.0e-12,7200.0,100000,50,true,false,57,200,1,200.0,200,200,1600,1600,0,63960200,99.9378125,0.2,0.8947029113769531,-2363.625792797503,-2363.6257927985066,1.6e-13,3.3e-13,6.0e-14,3.6e-13,8.7e-14,9.6e-01,1636473,1600,2.1e-13,5.3e-12,3.5e-12,4.8e-13
rand_200_3200_0.001,rand,74981947,tol,55.1,1500,1.0e-12,7200.0,100000,50,true,false,80,200,1,200.0,200,200,3200,3200,0,127556,0.099653125,0.0,0.016987085342407227,-94.96980623664624,-94.96980623657933,2.7e-13,4.7e-13,2.2e-13,6.6e-13,2.3e-15,1.7e+00,4380750,3200,3.5e-13,2.7e-13,5.4e-13,8.5e-13
rand_200_3200_1.0,rand,473171981,tol,583.3,800,1.0e-12,7200.0,100000,50,true,false,80,200,1,200.0,200,200,3200,3200,0,127960200,99.96890625,0.5,1.7353968620300293,-1986.0936354160915,-1986.0936354146686,2.5e-13,1.0e-13,2.2e-13,7.0e-13,1.6e-13,2.3e+00,2616670,3200,3.6e-13,1.7e-12,1.3e-11,1.0e-12
rand_200_6400_0.001,rand,199619049,tol,67.1,1050,1.0e-12,7200.0,100000,50,true,false,114,200,1,200.0,200,200,6400,6400,0,254878,0.09956171875,0.0,0.04164290428161621,-61.227868694936355,-61.227868694926926,5.2e-14,4.2e-13,1.4e-13,6.7e-13,2.0e-14,6.7e-01,2985042,6400,7.5e-14,2.5e-13,3.4e-13,9.8e-13
rand_200_12800_0.001,rand,940160956,tol,73.0,400,1.0e-12,7200.0,100000,50,true,false,160,200,1,200.0,200,200,12800,12800,0,509683,0.0995474609375,0.0,0.06021523475646973,-23.295865753982483,-23.295865753980838,1.7e-14,8.2e-15,1.2e-14,8.5e-14,3.1e-15,3.1e+00,2022544,12800,3.4e-14,4.6e-15,2.5e-14,1.7e-13
rand_400_800_1.0,rand,85150322,tol,410.1,950,1.0e-12,7200.0,100000,50,true,false,40,400,1,400.0,400,400,800,800,0,127840400,99.8753125,0.4,1.1382498741149902,-8081.041732476969,-8081.041732479744,1.3e-13,9.9e-14,3.1e-14,5.5e-13,1.6e-13,3.0e+00,2986677,800,1.7e-13,2.3e-12,3.8e-12,7.0e-13
rand_400_1600_0.0001,rand,971366104,tol,81.4,2550,1.0e-12,7200.0,100000,50,true,false,57,400,1,400.0,400,400,1600,1600,0,25924,0.0101265625,0.0,0.008107900619506836,-264.67897066472875,-264.67897066472887,1.8e-16,7.4e-14,2.1e-13,7.6e-13,8.9e-16,2.0e+01,7206929,1600,6.4e-16,1.3e-14,3.4e-13,8.1e-13
rand_400_1600_1.0,rand,732760939,tol,848.2,1200,1.0e-12,7200.0,100000,50,true,false,57,400,1,400.0,400,400,1600,1600,0,255840400,99.93765625,0.5,2.44583797454834,-7812.495442911056,-7812.495442911343,1.4e-14,1.5e-13,5.5e-14,5.4e-13,5.6e-14,2.3e+00,4632325,1600,1.8e-14,3.5e-12,6.2e-12,6.9e-13
rand_400_3200_0.0001,rand,492529223,tol,134.4,3850,1.0e-12,7200.0,100000,50,true,false,80,400,1,400.0,400,400,3200,3200,0,51467,0.0100521484375,0.0,0.0363309383392334,-342.7240541326808,-342.7240541326816,1.1e-15,2.1e-15,8.6e-13,4.3e-16,1.0e-16,9.6e+00,8887958,3200,8.3e-17,6.5e-16,1.2e-12,6.3e-16
rand_400_6400_0.0001,rand,916383521,tol,98.8,1800,1.0e-12,7200.0,100000,50,true,false,114,400,1,400.0,400,400,6400,6400,0,102520,0.01001171875,0.0,0.04269695281982422,-223.38759737082862,-223.38759737081503,2.8e-14,9.2e-13,4.8e-14,1.6e-13,1.2e-16,7.5e-01,4132971,6400,3.1e-14,1.9e-13,6.8e-14,1.7e-13
rand_400_12800_0.0001,rand,975528625,tol,266.8,2950,1.0e-12,7200.0,100000,50,true,false,160,400,1,400.0,400,400,12800,12800,0,204649,0.009992626953125,0.0,0.20035910606384277,-205.3198924780041,-205.31989247792663,1.7e-13,5.6e-13,1.5e-13,3.1e-13,1.7e-15,1.2e+00,10378345,12800,1.9e-13,1.6e-13,2.4e-13,3.5e-13
rand_400_25600_0.0001,rand,773156656,time,7200.1,40659,1.0e-12,7200.0,100000,50,true,false,227,400,1,400.0,400,400,25600,25600,0,408909,0.0099831298828125,0.1,0.29340100288391113,-86.12471232966305,-86.12471232968754,1.0e-13,9.1e-12,Inf,Inf,8.2e-15,5.9e-01,154408768,25600,1.4e-13,2.0e-12,1.9e-12,5.3e-12
rand_400_51200_0.0001,rand,434330530,tol,772.8,2000,1.0e-12,7200.0,100000,50,true,false,320,400,1,400.0,400,400,51200,51200,0,817571,0.00998011474609375,0.2,0.27559494972229004,-21.218555226688274,-21.21855522666451,2.9e-13,2.8e-13,1.2e-13,7.7e-13,1.0e-15,8.5e-01,11788442,51200,5.5e-13,7.5e-14,1.7e-13,1.4e-12
rand_200_20x10_2500_1.0,rand,84223329,tol,554.3,800,1.0e-12,7200.0,100000,50,true,false,71,400,21,19.047619047619047,10,200,2500,2500,0,104958400,99.96038095238096,0.5,1.5552260875701904,-4232.143654740144,-4232.143654737248,2.6e-13,1.5e-13,9.9e-14,7.8e-13,1.8e-13,1.9e+00,3113245,2500,3.4e-13,1.8e-12,5.9e-12,1.0e-12
rand_5x200_5000_0.001,rand,350918615,tol,351.1,1950,1.0e-12,7200.0,100000,50,true,false,100,1000,5,200.0,200,200,5000,5000,0,995721,0.0995721,0.1,0.23369097709655762,-661.2218233706086,-661.2218233706327,1.4e-14,8.3e-13,6.2e-14,2.3e-13,1.4e-14,6.0e-01,17838579,5000,1.8e-14,4.2e-13,2.8e-13,3.1e-13
rand_50x50_5000_0.01,rand,545101698,tol,1518.4,2400,1.0e-12,7200.0,100000,50,true,false,50,2500,50,50.0,50,50,5000,5000,0,6374063,1.01985008,0.4,0.6506049633026123,-1941.0184182051848,-1941.0184182045705,7.0e-14,2.2e-13,7.9e-13,2.9e-13,2.9e-14,6.6e-01,60392060,5000,1.6e-13,1.8e-13,8.6e-12,6.7e-13
rand_100_100x10_5000_1.0,rand,909351052,tol,3173.9,2400,1.0e-12,7200.0,100000,50,true,false,100,1100,101,10.891089108910892,10,100,5000,5000,0,99981100,99.9811,0.7,2.389385938644409,-2199.0399220420904,-2199.0399220421446,6.3e-15,3.6e-14,9.0e-13,1.2e-13,7.1e-16,7.0e-01,32620936,5000,1.2e-14,1.9e-13,3.7e-11,2.3e-13
rand_20x100_10000_0.001,rand,901859648,tol,3875.9,10400,1.0e-12,7200.0,100000,50,true,false,100,2000,20,100.0,100,100,10000,10000,0,1982094,0.0991047,0.3,0.4751760959625244,-824.095036109392,-824.0950361089228,1.8e-13,9.9e-14,3.4e-13,9.6e-13,1.7e-13,1.9e+00,204305163,10000,2.8e-13,3.8e-14,1.6e-12,1.5e-12
rand_400_10x40_10000_0.01,rand,270384289,time,7200.1,8868,1.0e-12,7200.0,100000,50,true,false,142,800,11,72.72727272727273,40,400,10000,10000,0,17560296,0.9977440909090909,0.4,0.6128931045532227,-1218.1113307192809,-1218.1113307167127,8.2e-13,4.9e-13,Inf,Inf,4.5e-13,1.8e+00,45779273,10000,1.1e-12,8.5e-13,5.7e-12,4.8e-12
\end{filecontents*}
\pgfplotstabletypeset[
    col sep=comma,
    string type,
    columns={instance,statusCode,time,iter,gapOriginal,pinfOriginal,dinfOriginal,complOriginal},
    columns/instance/.style={column name=\textbf{Instance}, column type=l},
    columns/statusCode/.style={column name=\textbf{Status}, column type=c},
    columns/time/.style={column name=\textbf{Time~[s]}, column type=r},
    columns/iter/.style={column name=\textbf{Iter}, column type=r},
    columns/gapOriginal/.style={column name=\texttt{gap}, column type=c},
    columns/pinfOriginal/.style={column name=\texttt{pinf}, column type=c},
    columns/dinfOriginal/.style={column name=\texttt{dinf}, column type=c},
    columns/complOriginal/.style={column name=\texttt{compl}, column type=c},
    every head row/.style={before row=\toprule, after row=\midrule},
    every last row/.style={after row=\bottomrule}
]{rand.csv}
}
\caption{Results of the Augmented Mixing Method with \texttt{tol=1e-12} on a wide range of randomly generated SDPs.}
\label{table:rand}
\end{table}

\subsection{Computations with extended floating-point arithmetic}\label{sec:arbitrary_precision}

To conclude our computational experiments, we evaluate the two-phase approach of the Augmented Mixing Method, as described in Section~\ref{sec:warmstart}, using extended floating-point arithmetic. The test instances are generated as in Section~\ref{sec:rand}, but are generally smaller to account for the increased computational cost associated with higher-precision arithmetic. We use the \texttt{Double64} data type from the \texttt{DoubleFloats.jl} package~\cite{sarnoff2022doublefloats}. Table~\ref{table:high} presents results comparing the Augmented Mixing Method to the interior-point solver \texttt{Hypatia}. Both solvers target a solution accuracy of $10^{-20}$.

The results in Table~\ref{table:high} show that the Augmented Mixing Method, guided by its warm-start strategy, can efficiently compute high-accuracy solutions in extended precision. With only a few exceptions around $10^{-20}$, the primal-dual errors are typically on the order of $10^{-21}$ or even smaller. The total number of iterations remains modest, except for the last instance, which requires~\num{17400} iterations to reach the desired precision. The iteration counts in parentheses show that both phases typically require a similar number of iterations, indicating steady progress throughout the solution process. Accordingly, the runtimes in parentheses suggest that switching from \texttt{Float64} to \texttt{Double64} increases the per-iteration cost by a factor typically between~$5$ and~$20$.

In terms of runtime, the Augmented Mixing Method outperforms \texttt{Hypatia} on all but the two smallest instances. On the three largest instances, \texttt{Hypatia} even reaches the time limit. Moreover, \texttt{Hypatia} did not attain the target precision of $10^{-20}$ on any instance. According to solver logs, it often stalls and terminates prematurely due to stagnating progress before reaching the desired tolerance.

\begin{table}[ht]
\centering
\resizebox{1.0\textwidth}{!}{
\setlength{\tabcolsep}{7pt}
\renewcommand{\arraystretch}{1.1}
\begin{filecontents*}{high.csv}
instance,type,seed,statusCode,time,iter,tol,timeLimit,maxIters,itersZ,scaling,warmStart,k,n,numBlocks,avgBlockSize,minBlockSize,maxBlockSize,m,equations,inequalities,nnzA,densityAPercent,timeScaling,timeData,primalObj,dualObj,gap,pinf,dinf,compl,complNotPsdAbs,mu,evals,nnzY,gapOriginal,pinfOriginal,dinfOriginal,complOriginal,timeHypatia,maxErrorHypatia,time1,iter1,time2,iter2
rand_10x10_100_1.0,high,235442785,tol,7.9 (1.6),1300 (850),1.0e-20,7198.4342629909515,99150,50,true,true,10,100,10,10.0,10,10,100,100,0,99100,99.1,0.0,0.0013930797576904297,-2.27875830087178170923579440544911751e+02,-2.27875830087178170923599575364688908e+02,2.1e-23,3.5e-23,4.5e-22,2.6e-22,2.2e-23,1.7e+00,292690,100,4.4e-23,1.3e-22,4.4e-21,5.5e-22,2.9,1.3e-16,1.6,850,6.3,450
rand_50_100_1.0,high,48535243,tol,24.4 (1.7),2400 (1150),1.0e-20,7198.274051904678,98850,50,true,true,15,50,1,50.0,50,50,100,100,0,247550,99.02,0.0,0.00766301155090332,-2.84003183857616948092642698368484396e+02,-2.84003183857616948087423437012214134e+02,6.7e-21,9.8e-21,1.2e-21,5.1e-21,7.8e-22,9.5e-01,250000,100,9.2e-21,8.0e-20,1.7e-20,7.0e-21,17.3,3.4e-18,1.7,1150,22.7,1250
rand_50_200_0.01,high,292392538,tol,3.5 (0.5),500 (300),1.0e-20,7199.479796886444,99700,50,true,true,20,50,1,50.0,50,50,200,200,0,5135,1.027,0.0,0.0006909370422363281,-4.00170997867386554691453950566565711e+01,-4.00170997867386554691486205809305815e+01,3.2e-23,3.7e-22,8.7e-22,1.3e-22,1.6e-23,2.2e+00,103394,200,4.0e-23,2.8e-22,1.7e-21,1.6e-22,31.9,1.3e-17,0.5,300,3.0,200
rand_100_10x10_100_1.0,high,998984545,tol,76.2 (3.8),1050 (550),1.0e-20,7196.211987018585,99450,50,true,true,15,200,11,18.181818181818183,10,100,100,100,0,1089200,99.01818181818182,0.1,0.057897090911865234,-1.99103928478341061806712199929105576e+03,-1.99103928478341061806315343774463102e+03,8.1e-22,5.2e-22,2.6e-22,1.6e-21,7.9e-22,2.6e+00,608150,100,1.0e-21,4.4e-21,8.1e-21,1.9e-21,102.6,8.4e-18,3.8,550,72.4,500
rand_100_400_0.01,high,765923170,tol,26.7 (2.1),1050 (450),1.0e-20,7197.889054059982,99550,50,true,true,29,100,1,100.0,100,100,400,400,0,39621,0.990525,0.0,0.004875898361206055,-8.95653664803259048568907011915813556e+01,-8.95653664803259048567396538489393637e+01,6.6e-22,9.0e-21,4.2e-22,1.5e-21,1.3e-22,5.5e-01,240836,400,8.4e-22,1.1e-20,1.2e-21,1.9e-21,463.9,1.9e-18,2.1,450,24.5,600
rand_100_400_1.0,high,958176910,tol,94.0 (10.0),800 (450),1.0e-20,7189.964840173721,99550,50,true,true,29,100,1,100.0,100,100,400,400,0,3990100,99.7525,0.4,0.23005986213684082,-8.35847650239686479373145703003425464e+02,-8.35847650239686479374379672705315221e+02,5.5e-22,6.2e-21,8.7e-22,3.2e-21,1.0e-21,8.9e-01,140000,400,7.4e-22,7.2e-20,2.6e-20,4.4e-21,480.0,1.8e-18,10.0,450,84.0,350
rand_2x100_200_0.01,high,612244091,tol,17.2 (3.1),600 (400),1.0e-20,7196.872037887573,99600,50,true,true,20,200,2,100.0,100,100,200,200,0,39582,0.98955,0.0,0.005712032318115234,-2.28974470223766972082864148215254006e+02,-2.28974470223766972081656120218400384e+02,2.0e-21,3.3e-21,2.4e-21,2.7e-21,2.7e-21,1.6e+00,236559,200,2.6e-21,3.9e-21,1.0e-20,3.5e-21,375.0,2.9e-17,3.1,400,14.1,200
rand_20x20_200_0.01,high,378414885,tol,36.5 (4.8),1700 (800),1.0e-20,7195.206212043762,99200,50,true,true,20,400,20,20.0,20,20,200,200,0,15608,0.9755,0.0,0.0022940635681152344,-2.84571163954726907988681505427758262e+02,-2.84571163954726907989861943442503484e+02,1.5e-21,5.8e-22,7.9e-21,4.9e-21,1.4e-21,8.1e+00,1972039,200,2.1e-21,2.6e-22,2.5e-20,7.1e-21,81.3,2.6e-16,4.8,800,31.8,900
rand_200_20x20_200_1.0,high,258599152,tol,732.4 (35.1),1450 (750),1.0e-20,7164.886979103088,99250,50,true,true,20,600,21,28.571428571428573,20,200,200,200,0,9552600,99.50625,0.7,0.29662299156188965,-8.64551092305115663594097465000486908e+03,-8.64551092305115663592962347649822512e+03,5.5e-22,3.6e-22,2.1e-22,2.4e-21,4.9e-22,3.0e+00,1793174,200,6.6e-22,3.8e-21,1.3e-20,2.9e-21,1300.4,2.5e-18,35.1,750,697.3,700
rand_200_800_0.001,high,485994070,tol,460.1 (14.3),4100 (1550),1.0e-20,7185.747787952423,98450,50,true,true,40,200,1,200.0,200,200,800,800,0,32001,0.100003125,0.0,0.008514881134033203,-1.82899462527429238403264078061329472e+02,-1.82899462527429238403317625372022007e+02,1.3e-22,4.6e-21,1.4e-22,6.0e-22,1.7e-23,5.8e-01,2279920,800,1.5e-22,2.7e-21,3.5e-22,6.9e-22,7284.6,4.3e-17,14.3,1550,445.9,2550
rand_300_600_0.001,high,977849511,tol,240.2 (11.4),1250 (550),1.0e-20,7188.631647109985,99450,50,true,true,35,300,1,300.0,300,300,600,600,0,54051,0.10009444444444444,0.0,0.013062000274658203,-3.04369559775001535897321692747869939e+02,-3.04369559775001535897453282370267441e+02,1.9e-22,7.5e-21,1.8e-23,1.1e-22,2.0e-23,5.3e-01,850996,600,2.2e-22,4.5e-21,5.5e-23,1.3e-22,7520.0,9.9e-11,11.4,550,228.9,700
rand_4x150_600_0.001,high,170430419,tol,4512.9 (145.5),17400 (6700),1.0e-20,7054.516135931015,93300,50,true,true,35,600,4,150.0,150,150,600,600,0,52986,0.09812222222222222,0.0,0.011495828628540039,-5.64139857794180449630527496825689817e+02,-5.6413985779418044962981588796251128e+02,5.3e-22,9.4e-21,1.0e-20,5.3e-22,9.2e-23,1.6e+00,35429271,600,6.3e-22,4.2e-21,3.3e-20,6.3e-22,7297.1,2.4e-13,145.5,6700,4367.4,10700
\end{filecontents*}
\pgfplotstabletypeset[
    col sep=comma,
    string type,
    columns={instance,statusCode,time,iter,gapOriginal,pinfOriginal,dinfOriginal,complOriginal,timeHypatia,maxErrorHypatia},
    columns/instance/.style={column name=\textbf{Instance}, column type=l},
    columns/statusCode/.style={column name=\textbf{Status}, column type=|c},
    columns/time/.style={column name=\textbf{Time~[s]}, column type=r},
    columns/iter/.style={column name=\textbf{Iter\phantom{abc}}, column type=r},
    columns/gapOriginal/.style={column name=\texttt{gap}, column type=c},
    columns/pinfOriginal/.style={column name=\texttt{pinf}, column type=c},
    columns/dinfOriginal/.style={column name=\texttt{dinf}, column type=c},
    columns/complOriginal/.style={column name=\texttt{compl}, column type=c|},
    columns/timeHypatia/.style={column name=\textbf{Time~[s]}, column type=r},
    columns/maxErrorHypatia/.style={column name=\textbf{Max. error}, column type=c},
    every head row/.style={
        before row={
            \toprule
            & \multicolumn{7}{c|}{\textbf{Augmented Mixing Method}} & \multicolumn{2}{c}{\textbf{Hypatia}} \\
        },
        after row=\midrule
    },
    every last row/.style={after row=\bottomrule}
]{high.csv}
}
\caption{Results of the Augmented Mixing Method and the interior-point solver \texttt{Hypatia} on randomly generated SDPs using the \texttt{Double64} data type from the \texttt{DoubleFloats.jl} package~\cite{sarnoff2022doublefloats}. The target relative precision for both solvers is set to $10^{-20}$. For the Augmented Mixing Method, the two-phase warm-starting approach described in Section~\ref{sec:warmstart} is applied. The reported runtime and iteration counts include both phases, with the values in parentheses indicating Phase~1 (\texttt{Float64}).}
\label{table:high}
\end{table}

\section{Conclusion} \label{sec:conclusion}

In this paper, we have presented the \emph{Augmented Mixing Method}, a novel SDP solver that combines an inexact augmented Lagrangian method with a block coordinate descent approach. The computational results in Section~\ref{sec:computational_experiments} demonstrate the strong practical performance of the Augmented Mixing Method using its default settings across various types of SDPs. It frequently yields solutions of significantly higher accuracy than those produced by state-of-the-art interior-point methods. Accurate solutions are even obtained for SDPs with numbers of affine constraints that exceed the capabilities of interior-point methods and most first-order methods.

A key feature of the Augmented Mixing Method is the direct handling of inequality constraints without maintaining slack variables in the algorithm. Another important component for achieving high-accuracy primal-dual solutions is the dynamic penalty parameter update scheme described in Section~\ref{sec:update}. Compared to other solvers based on the low-rank factorization framework, this update scheme makes the initial choice of the penalty parameter less critical and leads to balanced progress of both the primal and dual variables. In contrast, other solvers may produce approximate primal feasible solutions more quickly, but often yield dual solutions of only moderate accuracy. This advantage is reflected in the success of the Augmented Mixing Method in Section~\ref{sec:arbitrary_precision} when solving SDPs using extended floating-point arithmetic. Apart from the external L-BFGS solver used for subproblems, most of the computational time is spent performing (sparse) matrix-vector multiplications, making the method well-suited for use with extended precision arithmetic.

\paragraph{Future research.} Despite its strong practical performance, further investigation into the convergence theory of the Augmented Mixing Method is warranted. Section~\ref{sec:fail} shows that without additional assumptions convergence cannot be guaranteed in general. Moreover, we believe that a more sophisticated scaling strategy than the one presented in Section~\ref{sec:scaling} could further improve practical performance. We have also refrained from implementing a dynamic update strategy for the maximum rank $k$ in the Burer--Monteiro factorization, which is commonly used by other low-rank solvers to improve scalability, especially with respect to the matrix dimension $n$. In our experience, the Augmented Mixing Method often becomes trapped in local minima if the rank $k$ is not chosen as in~\eqref{eq:our_k}, even when the optimal solution of the SDP has significantly lower rank. Finally, a potential enhancement would be to implement an active-set style strategy in the Augmented Mixing Method to ignore inactive inequality constraints during iterations and reconsider them only when they become relevant. This could be particularly useful in the context of cutting-plane approaches based on SDP relaxations for combinatorial optimization problems.

\section*{Acknowledgements}
This research was funded in part by the Austrian Science Fund (FWF) [10.55776/DOC78]. For open access purposes, the author has applied a CC BY public copyright license to any author-accepted manuscript version arising from this submission. We thank two anonymous reviewers for their constructive comments, which helped improve this paper.

\paragraph{Disclosure statement.}
The authors report there are no competing interests to declare. The linguistic revision of this paper was carried out with the support of ChatGPT.

\bibliographystyle{plain}
\bibliography{references}

\end{document}